\documentclass[11pt,english]{article}

\usepackage{stocmds}
\usepackage{chemarr,bm}
\usepackage{mathtools}
\usepackage{algorithm2e}

\usepackage{amsmath,amssymb}

  \newtheorem{proof}{Proof}
  \newtheorem{example}{Example}

 \makeatletter
\newcommand*{\rom}[1]{\expandafter\@slowromancap\romannumeral #1@}
\makeatother

\usepackage{tikz}
\newcommand*\circled[1]{\tikz[]{ 
            \node[shape=circle,draw,inner sep=1pt,anchor=base] (char)
            {\scriptsize #1};}}

\usepackage{subcaption}
\usepackage{etoolbox}
\let\bbordermatrix\bordermatrix
\patchcmd{\bbordermatrix}{8.75}{4.75}{}{}
\patchcmd{\bbordermatrix}{\left(}{\left[}{}{}
\patchcmd{\bbordermatrix}{\right)}{\right]}{}{}
\usepackage{pifont}
\usepackage{blkarray}
\usepackage{multirow}

\usepackage{hyperref}

\newcommand{\harr}{\stackrel{\displaystyle\longrightarrow}{\longleftarrow}}

\newcommand{\duarr}{\mathbin{\hspace{.25em}\downarrow\hspace{-.1em}\uparrow}}
\begin{document}

\title
{ A Multi-Time-Scale Analysis of Chemical Reaction Networks : II. Stochastic Systems
\thanks{Supported in part by NSF Grants DMS \# 9517884
    and 131974 and NIH Grant \# GM 29123 to H. G. Othmer  and by National Research Foundation of Korea (2014R1A1A2054976) to CH. Lee.} 
    }
 \author{Xingye Kan$^{1}$, Chang Hyeong Lee$^{2}$, Hans G. Othmer$^1$\\
 \\
1. School of Mathematics  \\
              University of Minnesota\\
              Minneapolis, MN 55455\\
   E-mail: xkan@umn.edu, \;\; othmer@math.umn.edu   \\\\
2. Ulsan National Institute of Science and Technology\\
Ulsan Metropolitan City 698-798\\
 South Korea\\
 E-mail: chlee@unist.ac.kr
 \thanks{All authors contributed equally to this work.}}

\date{\today}

\maketitle

\begin{abstract}
  We consider stochastic descriptions of chemical reaction networks in which there are
  both fast and slow reactions, and for which the time scales are widely
  separated.  We develop a computational algorithm that produces  the generator
  of the full chemical master equation for arbitrary systems, and show how to obtain
  a reduced equation that governs the evolution on the slow time scale. This is
  done by applying a state space decomposition to the full equation that 
  leads to the reduced dynamics in terms of certain projections and the
  invariant distributions of the fast system. The rates or propensities of the reduced
  system are shown to be the rates of the slow reactions conditioned on the
  expectations of fast steps.  We also show that the generator of the
  reduced system is a Markov generator, and we present an efficient stochastic simulation
  algorithm for the slow time scale dynamics. We illustrate the numerical
  accuracy of the approximation by simulating several examples. Graph-theoretic
  techniques are used throughout to describe the structure of the reaction
  network and the state-space transitions accessible under the dynamics.

\end{abstract}

 \section{Introduction and background }

 While singular perturbations techniques and the quasi-steady-state
 approximation (QSSA) have a long history of use in deterministic descriptions
 of chemical reaction kinetics (\cf  ~\cite{Lee:2009:MTS} for a review), it was
 apparently not applied to discrete stochastic descriptions of chemical kinetics
 until Janssen proposed a method for adiabatic elimination of fast variables in
 stochastic chemical reaction networks
 \cite{Janssen:1989:EFVII,Janssen:1989:EFVIII}.  Janssen began with a master
 equation description and, using projection techniques related to volume
 expansions, obtained a reduced master equation and showed in some examples that
 an intermediate chemical species can be eliminated in a network if a reaction
 occurs much faster than the others. Moreover, he also gave examples in which
 the master equation cannot be reduced via these techniques
 \cite{Janssen:1989:EFVIII}. It is sometimes assumed that the results of a
 deterministic reduction of a reaction network yields correct results for the
 stochastic description, but examples show that this is not correct
 \cite{Thomas:2011:CLS}.

 The Gillespie algorithm \cite{Gillespie:2007:SSC} is the most widely-used
 algorithm for simulating stochastic reactions but it can be very inefficient
 when there are multiple time scales in the reaction dynamics.  Stochastic
 simulation algorithms of multiscale reaction networks were developed more or
 less simultaneously by Haseltine and Rawlings (\cite{Haseltine:2002:ASC})
 and Rao and Arkin (\cite{Rao:2003:SCK}). The former authors formulated the
 changes in the species numbers $n$ in terms of extents $\eta$ by defining
 $n(t)=n(0)+\nu\cE\eta(t)$\footnote{The notation used is defined later.}, and
 divided the reaction extents into those of slow reactions and fast
 reactions.  Rao and Arkin (\cite{Rao:2003:SCK}) divided the set of
 \emph{species} into `primary' $(y)$ and `intermediate or ephemeral' $(z)$, and
 assumed two conditions on the conditional variable $z|y$, (i) the conditional
 variable $z|y$ is Markovian and (ii) $z|y$ is at quasi-steady-state on the slow
 time scale.  However, the first assumption was not justified and one cannot
 classify species as slow and fast -- rather it is reactions that are slow or
 fast and species can participate in both.  Mastny {\it et al.}
 (\cite{Mastny:2007}) applied singular perturbation analysis to remove the
 QSS species for a number of model systems having small populations, and for
 networks where non-QSS species have large populations, they also utilized the
 $\Omega$-expansion to reduce the master equation.  More recent work deals with
 the role played by the reduction on the level of noise in the solution
 \cite{Srivastava:2011:SQS}.

Several others have proposed and implemented variations of these two approaches
for hybrid simulations of stochastic systems. These include partitioning the
system dynamically \cite{Salis:2005:AHS} and solving the fast reactions using
an SDE approximation. The effective reaction rate expressions for Hill-type
kinetics derived from singular perturbation analysis of the deterministic system
of equations have been used in stochastic simulations and the results compared
to simulations for the full system \cite{Bundschuh:2003:FSV}.  There is also a
partitioning method based on the variance of species, which led to a hybrid
method coupling deterministic (small variance) and stochastic (large variance)
by Hellander and L{\"o}tstedt (\cite{Hellander:2007}).  Goutsias showed that
fast reactions can be eliminated when probabilities of slow reactions depend at
mostly linearly on the states of fast reactions
\cite{Goutsias:2005:QAF}. Utilizing the degree of advancement or extent, he
followed Haseltine and Rawling's approach in order to separate variables into
fast and slow variables. By utilizing a Taylor series expansion he showed that
when the slow transition rates or propensity functions depend linearly on fast
extents, the fast reaction kinetics of a stochastic biochemical system can be
approximated by the conditional mean of extents for the fast kinetics. Although his
approach is rigorous, it has some drawbacks. In many nonlinear systems, slow
reactions depend nonlinearly on fast reactions. For example, if a slow reaction
subsystem consists of bimolecular reactions with two different reactants which
are affected by fast reactions, his method cannot be applied due to the
nonlinear dependence of the slow reactions on fast reactions.  Peles {\it et al}
applied singular perturbation theory and finite state projection method to
obtain an approximate master equation for certain linear reaction systems
\cite{Samad:2006:SHS}, where clusters of fast-transitioned states have been
identified. A rigorous mathematical framework is lacking in their
derivation, which we provide herein.

Computational methods have also been developed by many groups. Cao {\it et al.},
(\cite{Cao:2005:SSS}) proposed a slow scale stochastic simulation algorithm
(SSA). They first identified fast and slow reaction channels and defined fast
and slow species according as the species are changed by fast reactions or not.
Assuming that the fast processes are stable, they approximated the propensity
functions on a slow time scale, and showed that computation of the effective
propensity functions requires first and second moments of the steady-state
distribution of the partitioned fast reaction subsystem. They illustrated
numerical results for a system with one independent fast species.  Using a
partitioning approach similar to that in 
\cite{Haseltine:2002:ASC,Rao:2003:SCK,Cao:2005:SSS}, Chevalier and EI-Samad
(\cite{Chevalier:2009}) derived evolution equations for both fast and slow
reactions that led to an SSA with slow reaction trajectories in which the SSA is
run until the first occurrence of a slow reaction.  E {\it et
  al.}\,(\cite{E:2005:NSS}) proposed a nested stochastic simulation
algorithm with inner loops for the fast reactions and outer loops for the slow
reactions. Strong convergence of the nested SSA has been proved recently by
Huang and Liu \cite{Huang:2014} and a speed up of the the algorithm is proposed
by using the tau-leaping method as the inner solver.  More recently, Kim et al
(\cite{Kim:2014}) investigated the validity of the stochastic QSSA, where the
propensity functions resulted from their deterministic counterpart. Under the
moment closure assumption, they have shown that the stochastic QSSA is accurate
for a two-dimensional monostable system, if the deterministic QSS solution is
not very sensitive.  Goutsias and Jenkinson (\cite{Goutsias:2013}) provide an
excellent review article covering recently developed techniques on approximating
and solving master equations for complex networks. Other recent papers include \cite{Cotter:2015:ARX, Smith:2015:ARX}. 

As described above, many computational results and some algorithms have been
reported to date, but a rigorous analysis of stochastic reaction networks that
involve multiple time scales based on singular perturbation has heretofore not been done. Preliminary work on
this by one of the authors \cite{Lee:2009} proposed a reduction method based on
a singular perturbation analysis for networks with two or more time scales under
the assumption that the sub-graph of fast steps in the network is strongly
connected.  Our objectives here are to develop a rigorous analytic framework for
the reduction of general stochastic networks with two widely-separated time
scales, to prove that the generator of the reduced system is Markovian, and to
illustrate the numerical accuracy and speedup of the reduction method for
several biological models.

\section{The general formulation for reacting systems}

\subsection{Deterministic evolution equations }
\label{Sec:2.1}

To set notation for the stochastic analysis, we first recall some notation used
in our previous analysis of deterministic systems \cite{Lee:2009:MTS}, hereafter
referred to as I. Further details can be found there.

Suppose that the reacting mixture contains the set $\cM$ of $m$ chemical species
$\cM_{i}$ that participate in a total of $r$ reactions. Let $ \nu_{i\ell}$ be the
stoichiometric coefficient  of
the $i^{th}$ species in the $\ell^{th}$ reaction. The $  \nu_{i\ell}$ are non-negative
integers that represent the normalized molar proportions or stoichiometric
coefficients of the species in a reaction.  Each reaction is written in the form
\begin{equation}
{\displaystyle
{\sum_{i}}^{reac.}  \nu_{i\ell}^{reac}\,\cM_{i} \rightarrow {\sum_{i}}^{prod}
  \nu_{i\ell}^{prod}\,\cM_{i}\hspace{.5in}\ell=1,\,\ldots\,r,} 
\label{graphrxn}
\end{equation}
 where the sums are over reactants and products, respectively in the $\ell^{th}$
reaction. In this formulation, the forward and reverse reaction of a reversible
pair are considered separately, as two irreversible reactions. Once the
reactants and products for each reaction are specified, the significant entities
so far as the network topology is concerned are not the species themselves, but
rather the linear combinations of species that appear as reactants or products
in the various elementary steps. These linear combinations of species are {\em
complexes} \cite{Horn:1972:GMA}, and we suppose that there are $p$ of them.  A
species may also be a complex as is the case for first-order reactions.  Once
the complexes are fixed, their composition is specified unambiguously, and we
let $\nu$ denote the $m\times p$ matrix whose $j^{th}$ column encodes the
stoichiometric amounts of the species in the $j^{th}$ complex. 

The set of reactions gives rise to a directed graph $\cG$ as follows.  Each
 complex is identified with a vertex $V_{j}$ in $\cG$ and a directed edge
 $E_{\ell}$ is introduced into $\cG$ for each reaction. 
 The topology of $\cG$ is encoded in its  vertex-edge {\em incidence matrix}
$\cE$, which is defined as follows.
\begin{equation}
\cE_{j\ell}=\left\{ \begin{array}{ll}
+1 & \mbox{ if }E_{\ell}\;\mbox{is incident at }V_{j}\mbox{ and is
directed toward it}  \\
-1 & \mbox{ if }E_{\ell}\:\mbox{is incident at }V_{j}\mbox{ and is
directed away from it}\\
\phantom{-}0 & \mbox{ otherwise}\end{array}\right.
\label{incid}
\end{equation}
 Since there are $p$ complexes and $r$ reactions,
 $\cE$ has $p$ rows and $r$ columns, and every column has exactly one $+1$ and
 one $-1$. Each edge carries a nonnegative weight $R_{\ell}(c)$ given by the intrinsic rate
 of the corresponding reaction.  
For example, the following table gives four classes of first-order reactions
 studied in Gadgil, {\it et al.}, (\cite{Gadgil:2005:SAF}) and two additional bimolecular reaction types.
\begin{table}[h]
\begin{center}
\begin{tabular}{|c|l|c|c|}
\hline 
Label&
~~~Type of reaction&
Reaction&
Rate\tabularnewline
\hline
\hline 
I&
Production from a source&
$\phi\rightarrow\cM_{i}$&
$k_{i}^{s}$\tabularnewline
\hline 
II&
Degradation&
$\cM_{i}\rightarrow\phi$&
$k_{i}^{d_1}n_{i}$\tabularnewline
\hline 
III&
Conversion&
$\cM_{j}\rightarrow \nu_i \cM_{i}$&
$k_{ij}^{con_1}n_{j}$\tabularnewline
\hline 
IV&
Catalytic production from source &
$\phi\stackrel{\cM_{j}}{\longrightarrow}\cM_{i}$&
$k_{ij}^{cat}n_{j}$\tabularnewline
\hline
V &  Bimolecular  degradation           & $\cM_j+\cM_k\rightarrow\phi$ & $k_{jk}^{d_2}n_jn_k$ \tabularnewline
\hline
VI        &  Bimolecular conversion & $\cM_j+\cM_k\rightarrow \nu_i\cM_i$ & $k_{ijk}^{con_2}n_jn_k$ \tabularnewline
\hline
\end{tabular}\end{center}
\caption{The four classes of first-order reactions considered in
\protect\cite{Gadgil:2005:SAF} and two types of bimolecular reactions. The relationship
between the deterministic and stochastic rates of these reactions 
  are discussed later. Rates are given in terms of number of molecules ($n$ will be introduced later).}
\label{rxntypes}
\end{table}
For either type III or VI reactions there may also be different types of
 products, {\it e.g.}, $A \rightarrow B + C$ may represent the decomposition of a
 complex. Inclusion of such types poses no difficulties, but if such reactions
 are reversible we restrict the type to uni- or bimolecular reactions. 

If at least one reaction of type I, II, IV, or V  is present the stoichiometric
matrix is
\begin{equation}
\nu=[~~\bI~~|~~{\textbf{0}}~~].
\label{nu0}
\end{equation}
wherein the column of zeroes in the complex matrix represents the null
 complex $\phi$, which by definition contains no time-varying species.
If the system is closed and the reactions are all first order 
\[
\nu=[~~\bI~~].
\]
If all species are also complexes, which can occur when there are inputs ($\emptyset\rightarrow \mathcal M_i$),
outputs ($\mathcal M_i\rightarrow\emptyset$), and first-order decay or conversion reactions ($\mathcal M_j \rightarrow \mathcal M_i$), the stoichiometric
matrix has the form 
\begin{equation}
\nu=[~~\bI~~| ~~\nu_1~~|~~{\textbf{0}}~~].
\label{nu2}
\end{equation}
wherein $\nu_1$ defines the stoichiometry of the higher-order complexes. 
In this case the  corresponding incidence matrix $\cE$ can be written as follows.
\begin{equation}
\cE=\left[~~\cE_1~~|~~\cE_2~~| ~~\cE_o \right] 
\end{equation}
where $\cE_1$ represents first-order reactions, $\cE_2$ represents second-order 
reactions, and $\cE_o$ represents input and output steps, all of the appropriate 
dimensions.
 An alternate form of the complex and incidence matrices
arises if the inputs or outputs are only of type I or II, for then  the null
complex $\phi$ can be omitted from $\nu$, the $\pm{1}$'s omitted from $\cE$,
and the inputs or outputs represented by a separate vector in the evolution
equations given below. In either case  the stoichiometry of the reactions and the
topology of the network are easily encoded in $\nu$ and $\cE$,
respectively. 

In this notation the  evolution of the composition of a reacting mixture is governed by
\begin{equation}
\frac{dc}{dt}=\nu\cE R (c), \ \ c(0)=c_0
\label{dd}
\end{equation}
where the $j^{th}$ column of $\nu$ gives the composition of the $j^{th}$ complex
and $R_\ell(c)$ is the rate of the $\ell^{th}$ reaction, or equivalently, the
flow on the $\ell^{th}$ edge of $\cG$.  The matrix $ \hat{\nu} \equiv \nu\cE$ is
called the stoichiometric matrix when the composition of complexes and the
topology of $\cG$ are not encoded separately, as we do here
\cite{Aris:1965:PRA}. One can interpret the factored form in (\ref{dd}) as
follows: the vector $R$ gives the flows on edges due to reactions of the
complexes, the incidence matrix maps this flow to the sum of all flows entering
and leaving a given node (a complex), and the matrix $\nu$ converts the net
change in a complex to the appropriate change in the molecular species. 

 A {\em component} is a connected subgraph $\cG_1 \subset \cG^0$ that is maximal
 with respect to the inclusion of edges, {\em i.e.}, if $\cG_2$ is a connected
 subgraph and $\cG_1\subset \cG_2 \subset \cG^0,$ then $\cG_1 = \cG_2$. An
 isolated vertex is a component and every vertex is contained in one and only
 one component. A directed graph $\cG$ is {\em strongly connected} if for every
 pair of vertices $(V_i,V_j)$ , $V_i$ is reachable from $V_j$ and vice-versa.
 Thus a directed graph is strongly connected if and only if there exists a
 closed, directed edge sequence that contains all the edges in the graph. A {\em
   strongly-connected component} of $\cG$ (a strong component or SCC for short)
 is a strongly-connected subgraph of a directed graph $\cG$ that is maximal with
 respect to inclusion of edges. An isolated vertex in a directed graph is a
 strong component, and every vertex is contained in one and only one component.
 Strong components in the directed graph $\cG$ are classified into three
 distinct types: sources, internal strong components and absorbing strong
 components.  A source is a subgraph which has outgoing edges to other strong
 components and has no incoming edges from other strong components. An internal
 strong component is a strong component in which edges from other strong
 components terminate and from which edges to other strong components
 originate. An absorbing strong component is a strong component from which no
 edges to other strong components originate. If $\cG$ has $p$ vertices and $q$
 strong components then it is easily shown that the rank of $\cE$ is $\rho(\cE)
 = p-q$ \cite{Chen:1971:AGT}.

For ideal mass-action kinetics, which we consider here,   the flow on the $\ell^{th}$ edge, which
originates at the $j^{th}$ vertex, depends only on the species in the $j^{th}$
complex, and the rate can be written as
\begin{equation}
R_{\ell}(c)=k_{\ell j} P_{j}(c)  \quad \mbox{where} \quad  P_{j}(c) =\prod_{i=1}^{m}(c_{i})^{\nu_{ij}}
\label{mak_1}
\end{equation}
for every reaction that involves the $j^{th}$ complex as the
reactant.\footnote{This form also includes  non-ideal mass action rate laws, but 
the concentrations in (\ref{mak_1}) are then replaced by the 
activities of the species in the reactant complex, and as a result  the flow on a edge may
depend on all species in the system.} Thus the rate vector can be
written 
\begin{equation}
R(c)={\cal K} P(c)\end{equation}
 where $\cal K$ is an $r\times p$ matrix with $k_{\ell j}>0$ if and
only if the $\ell^{th}$ edge leaves the $j^{th}$ vertex, and $k_{\ell j}=0$
otherwise. Since each row of $\cK$ has one and only one positive entry, we now denote 
the only positive entry in $\ell^{th}$ row by $k_{\ell}$.

The topology of the underlying graph $\cG$ enters into $\cal K$ as
follows. Define the exit matrix $\cE_{e}$ of $\cG$ by replacing all $1$'s in
$\cE$ by zeroes, and changing the sign of the resulting matrix. Let $\hat{K}$ be
the $r\times r$ diagonal matrix with the $k_{\ell}$'s, $\ell=1,\,\ldots r$,
along the diagonal. Then it is easy to see that ${\cal K}=\hat{K}\cE_{e}^{T}$
and therefore
\begin{equation}
\dfrac{dc}{dt}=\nu\cE\hat{K}\cE_{e}^{T} P(c)
\label{facteqn}
\end{equation}
 It follows from the definitions that (i) the $(p,q)^{th}$ entry,
$p\neq q$, of $\cE\hat{K}\cE_{e}^{T}$ is nonzero (and positive) if and only if
there is a directed edge $(q,p)\in\cG$, (ii) each diagonal entry of
$\cE\hat{K}\cE_{e}^{T}$ is minus the sum of the $k$'s for all edges that leave
the $j^{th}$ vertex, and (iii) the columns of $\cE\hat{K} \cE_{e}^{T}$ all sum
to zero, and so the rank of $\cE\hat{K}\cE_{e}^{T}$ is $\leq p-1$.

If one separates the inputs, which are constants, one can write this as 
\begin{equation} 
\dfrac{dc}{dt}=\nu\cE\hat{K}\cE_{e}^{T} P(c) + \Phi,
\label{facteqn1} 
\end{equation} 
where $\Phi$ is the constant input and both $\nu$ and $\cE$ are modified
appropriately. Herein we use the evolution equations in the form (\ref{facteqn})
unless stated otherwise.

One can also describe the evolution of a reacting system in terms of the number
of molecules present for each species. Let $n=(n_{1},n_{2},\ldots,n_{m})$ denote
the discrete composition vector whose $i^{th}$ component $n_{i}$ is the number
of molecules of species $\cM_{i}$ present in the volume $V$. This is related to
the composition vector $c$ by $n=\cN_{A}Vc$, where $\cN_{A}$ is Avagadro's
number, and although the $n_i$ take discrete values, they are regarded as
continuous when large numbers are present. From (\ref{dd}) we obtain the
deterministic evolution for $n$ as
\begin{equation}
\dfrac{dn}{dt}=\nu\cE\tilde{\cR}(n)
\label{evol_n}
\end{equation}
 where $\tilde{\cR}(n)\equiv\cN_{A}V R(n/\cN_{A} V)$. In particular, for 
 ideal mass-action kinetics 
\begin{eqnarray}
\tilde{\cR_{\ell}}(n)&=&\cN_{A}Vk_{\ell}{P}_{j}(n/\cN_{A} V)\\
&=&\cN_A Vk_{\ell}\prod_{i=1}^{m}\left(\dfrac{n_{i}}{\cN_{A}
  V}\right)^{\nu_{ij}}=\dfrac{k_{\ell}}{(\cN_{A}V)^{\sum_{i}\nu_{ij}-1}}\prod_{i=1}^{m}(n_{i})^{\nu_{ij}}=\hat{k}_{\ell}\prod_{i=1}^{m}(n_{i})^{\nu_{ij}}.    
\label{mak_2}       
\end{eqnarray}

\subsection{Invariants of reaction networks}

The kinematic invariants and the kinetic invariant manifolds in a deterministic
description of reactions in a constant-volume system are discussed in detail in
\cite{Othmer:1979:GTA}.  In general the concentration  space has the
decomposition 
\begin{equation}
{\bf R}_m = \cN[(\nu\cE)^T] \oplus  \cR[\nu\cE],
\label{decomp}
\end{equation} 
 where $\cN(A)$ denotes the null space of $A$, and $\cR(A)$ denotes the range of $A$.
The solution of  (\ref{dd}), which defines a curve in ${\bf R}_m$
through an initial point $c_0$,  can be written
$$
c(t) = c_0 + \nu\cE\int_0^t R(c(\tau))d \tau.
$$
This  shows that $c(t)-c_0 \in \cR(\nu\cE)$, and the intersection of the
translate of $ \cR(\nu\cE)$ by $c_0$, which formally is a coset of
$\cR(\nu\cE)$ with the non-negative cone ${\bf C}^+_m$ of ${\bf R}_m$, defines the
reaction simplex $\Omega(c_0)$. While this terminology has a long history
\cite{Aris:1965:PRA}, $\Omega$ is a simplex in the mathematical sense only if
the intersection of the coset with ${\bf C}^+_m$ is compact, which occurs if and
only if there is a vector $y >0 \in \cN[(\nu\cE)^T]$
\cite{Othmer:1979:GTA}. This is only guaranteed in closed systems, where the
total mass is conserved and $y$ comprises the molecular weights of the species,
and therefore in general we should call $\Omega$ the kinetic manifold, but we
retain the standard terminology.

First suppose that the system is closed -- the case of an open system is discussed
later. A vector $a\in {\bf R}_m$ defines an invariant linear combination of
concentrations if
\begin{eqnarray}
\langle a,\nu\cE R(c)\rangle   = 0, 
\label{inva}
\end{eqnarray}
and these  are called kinematic invariants if $a \in \cN[(\nu\cE)^T]$
\cite{Othmer:1979:GTA}. The following important properties of these invariants
are proven in \cite{Lee:2009:MTS}. 
\begin{itemize}
\item [{\bf P1}]  One can choose a basis for $\cN[(\nu \cE)]^T$ of vectors with
integer entries. 

\item [{\bf P2}]
  If the reaction simplex $\Omega(c_0)$ is compact, then there is a
  basis for $\cN[(\nu\cE)^T]$ for which all basis vectors have
  nonnegative integer entries.
\end{itemize}

 If the reactions are partitioned into subsets of fast and slow steps we can write 
\begin{equation} 
\nu\cE = \nu \left[\cE^f  \mid \cE^s ~\right]
\label{fsdecomp}
\end{equation}
and it follows that any invariant of the full system is simultaneously an
invariant of the fast and slow systems. We assume throughout that the slow and
fast reactions are independent, and therefore $  \cR(\nu\cE) = \cR(\nu\cE^{f})
\cup  \ \cR(\nu\cE^{s})$. Furthermore, one has that 
\begin{align}
\label{rinclu}
 \cR[\nu\cE^{(f,s)}] & \subseteq  \cR[\nu\cE]  \\
\label{nullinclu}
\cN[(\nu\cE)^T] & \subseteq  \cN[(\nu\cE^{(f,s)})^T].
\end{align}
that the ranges of the slow and fast subsets are no larger than that of the full
system, and the corresponding null spaces are no smaller, but the properties
{\bf P1, P2} of the full system do not necessarily carry over to the
subsystems. However, it follows from (\ref{nullinclu}) that one can define a map
$P^f:{\bf R}_m \rightarrow {\bf R}_{m-r_f} $ (where $r_f$ = dim$
\cR[\nu\cE^{f}]$) for the fast subsystem that represents a vector in
$\cN[(\nu\cE^f)^T]$ in terms of intrinsic coordinates on
$\cN[(\nu\cE^f)^T]$. The associated matrix $\cP^f$ has rows given by basis
vectors with integer components of $\cN[(\nu\cE^f)^T]$. It follows that the
reaction simplex for the fast subsystem is given by
$$\Omega_f(c_0) \equiv \{c \,: c \in c_0 + \cR[\nu\cE^f]\} \cap \bar{\bf R}_m^+ = \{c : \cP^f
c=\cP^fc_0 \equiv \tilde c \in {\bf R}_{m-r_f} \}\cap \bar{\bf R}_m^+.
$$
Here  $\tilde c$ represents a conserved quantity for the
fast subsystem, but it may vary as slow reactions occur.

If the system is open then one can regard the effect of inputs and outputs as
moving the dynamics between simplexes of fixed mass, as dictated by the 
evolution equations given at (\ref{facteqn}). This point of view will be useful
in the stochastic formulation. 

\section{The stochastic formulation} 

\subsection{The master equation} 

In a stochastic description the number of molecules of a species is too small
to be treated as a continuous variable -- they are random variables. We define
$N(t)=(n_1(t),n_2(t),$ $\dots,n_m(t))$, where $n_i(t)$ is as before, but now
$N(t)$ is a random vector. Under the assumption that the process is Markovian,
which is appropriate if all the relevant species are taken into account, the
evolution of $N(t)$ is governed by a continuous-time Markov process with
discrete states, and we denote the probability that $\{N(t)=n\} $ by $P(n,t)$.
The governing equation for the evolution of $P(n,t)$ is called the chemical
master equation, and is given as
\begin{equation}
\dfrac{d}{dt}P(n,t)=\sum_{\ell}\cR_{\ell}(n-\nu\cE_{(\ell)})\cdot
P(n-\nu\cE_{(\ell)},t)-\sum_{\ell}\cR_{\ell}(n)\cdot
P(n,t)
\label{prob_evol1}
\end{equation}
where  $\cE_{(\ell)}$, the $\ell^{th}$ column of $\cE$, denotes $\ell^{th}$
reaction and the stochastic rates have the form
\begin{equation}
\cR_{\ell}=c_{\ell}h_{j(\ell)}(n)\label{trans_1}.
\end{equation}
 Here $c_{\ell}$ is the probability per unit time that the molecular species in
the $j^{th}$ complex reacts, $j(\ell)$ denotes the reactant complex for the
$\ell^{th}$ reaction, and
$h_{j(\ell)}(n)$ is the number of independent combinations of the molecular
components in this complex \cite{Gadgil:2005:SAF}.  One sees in
(\ref{prob_evol1}) how the reaction graph $\cG$ determines the state transition
graph $\cG_s$ that defines the steps in the master equation.

In the stochastic analysis
all reactions are assumed to follow mass-action kinetics, and thus $c_{\ell} =\hat{
{k}}_{\ell}$, and the combinatorial coefficient is given by \footnote{This
formulation applies only to ideal solutions -- in nonideal solutions the number
of molecules must be replaced by an appropriate measure of its activity in the
solution. In particular, this involves a suitable description of diffusion when
the solution is not ideal \protect\cite{Othmer:1976:NEC,Schnell:2004:RKI}.}
\begin{equation}
h_{j(\ell)}=\prod_{i}{{n_{i} \choose \nu_{ij(\ell)}}}.\label{hpee}
\end{equation}
If $ \nu_{ij(\ell)} =1$ the stochastic rate reduces to (\ref{mak_2}) but for a
bimolecular reaction the rate of a step in the stochastic
framework is always smaller than in the deterministic framework.

The master equation (\ref{prob_evol1}) can be expressed more explicitly for
unimolecular and bimolecular mass action kinetic mechanisms. The general form
for the uni- and bimolecular molecular reactions given in Table \ref{rxntypes} is 
\begin{align*}
  \dfrac{dP(n,t)}{dt} &=
  \sum_{i=1}^s\bigg[K_{ii}^s(S_i^{-1}-1)P(n,t)+\sum_{j=1}^s\bigg(K_{ij}^{con_1}(S_i^{-\nu_i}\cS_{j}^{+1}-1)+K_{ij}^{cat}(S_i^{-1}-1)+K_{ii}^{d1}(S_i^{+1}-1)\\ 
 &+\sum_{k=1}^s\bigg(K_{ijk}^{con_2}(S_i^{-\nu_i}S_j^{+1}S_k^{+1}-1)+K_{jk}^{d2}(S_j^{+1}S_k^{+1}-1)\bigg)n_k    
 \bigg)(n_jP(n,t))\bigg].
\end{align*}
where $\cS_{i}^{k}$ is the shift operator that increases the $i^{th}$
component of $n$ by an integer amount $k$, and $S_{i}^{k}(n_{i}P(n,t))=S_{i}^{k}n_{i}\cdot P(S_{i}^{k}n,t)$.
In this expression the matrices $K^{s}, K^{d1}, K^{cat},K^{d2},$ $K^{con_1}$, and $K^{con_2}$ correspond
to the six types of reactions given in Table \ref{rxntypes}. %
\\
A compact representation of (\ref{prob_evol1}) is obtained by defining an
ordering of the accessible states, of which we first assume there finitely many
($n_s$) in number and which lie in the non-negative cone ${\bf C}^+_m \subset
{\bf R}_{n_s}$\footnote{The case in which there are infinitely many states will
  be discussed in Section \ref{Subsec:bound}.}.  The evolution of the vector
$p(t)$ of the probabilities of these states is governed by the matrix Kolmogorov
equation
\begin{eqnarray}
\dfrac{dp(t)}{dt} = K p(t), \label{kolmo}
\end{eqnarray}
where $K$ is the matrix of transition rates between the states, the entries of
which are defined as follows.  Let $n^i =(n_1^i,n_2^i, \dots,n_m^i)$  and
$n^j=(n_1^j,n_2^j, \dots,n_m^j)$ denote the $i^{th}$ and $j^{th}$ states of
the system and denote the $i^{th}$ and
$j^{th}$ entries of the vector $p$ as $p_i=P(n^i,t)$ and $p_j=P(n^j,t)$, resp.
Then the $(i,j)^{th}$ entry of $K$ is given by 
\begin{equation}
K_{ij}=\left\{\begin{array}{cc} \cR_{\ell}(n^j) &  \textrm{ if $n^i=n^j +
\nu\cE_{(\ell)}$ for some $\ell=1,\dots,r$ }\\             0 & \textrm{
otherwise } 
\end{array} 
\right.
\label{stoK}
\end{equation}
and the following pseudo-code shows how to
generate $K$ and the corresponding state transition graph $\cG_s$ for any reacting system. By starting with a given initial state vector, the algorithm checks all the reactions and adds resulting new state vectors into state space while updating the transition matrix, and then repeats the same procedure on the newly added state vector until no new state vector can be added.

\vspace*{10pt}
\small
\begin{algorithm}[H]\label{alg}
\KwData{initial state vector $V_1$ and reactions $\{R_\ell\}_{\ell=1}^r$}
\KwResult{transition matrix $K_{n_s\times n_s}$}
Initialization: set current state index $C_s=1$, set accessible state vector space $\{V\}$ to be $\{V_1\}$, set accessible state vector space size $n_s=1$, set transition rate matrix $K$  to be \{0\}\;
\While{current state index $\leq n_s$  }{ 
   set current reaction index to one: $C_r = 1$\;
   \While{current reaction $\leq r$ }
{
check if current reaction $R_C$ reacts from current state $V_C$\;
    \If{True}{
               get current state index: source $\gets C_s$\;
               get target state $V_T = V_C + \nu\cE_{(C_r)}$\;
               check if target state $V_T$ is already in $\{V\}$\;
                  \If{True}{
                            get the index $i$ of the state in $\{V\}$ that is equal to $V_T$\;
                            check if $K_{i,source} = 0$\;
                               \If{True}{
                                         update $K_{i,source}$\;}{}
               
               }\Else{
                 add $V_T$ to $\{V\}$: $V_{n_s} \gets V_T$\;
                 increase accessible state space size by one: $n_s \gets n_s+1$\;
                 update $K_{n_s,source}$\;
                 
}}{}
    increase current reaction index by one: $C_r \gets C_r + 1$\;}
    increase current state index by one: $C_s \gets C_s + 1$\;
    }
               
update diagonal entries $K_{jj} = -\sum_{i=1,i\neq j}^{n_s}K_{ij}$.              
\vspace*{10pt}                     
               \caption{An algorithm for generating the transition matrix $K$ for a chemical reaction network}
\end{algorithm}
\normalsize

Off-diagonal elements in the $i^{th}$ row of $K$ are transition rates at which
source states reach state $i$ in one step, whereas off-diagonal elements in the
$i^{th}$ column represent the  rates at which the $i^{th}$ state reaches its
target states in one reaction.  Since only uni- and bimolecular reactions are
realistic, the states that reach the current state in one step (the sources)
differ from the current state by at most two molecules, but those that can be
reached in one step (the targets) may involve more than two.  Under certain
orderings  $K$ has a narrow bandwidth when the number of states
$n_{s}$ is small, but expands with increasing $n_{s}$.  The number $n_s$ grows
combinatorially with the number of molecules -- when there are $m$ species and a
total of $N_0$ molecules
$$
n_s = { N_0 + m -1 \choose m -1}.
$$
If, for example, there are 4 species and 50 molecules,  the number of states
is 23,426.  In any case, if all reactive states  are accounted for
the matrix $K$ is the generator of a Markov chain \cite{Norris:1998}, {\em i.e.} $\sum_i K_{ij} =0
$ for each $j$, $ K_{ij} \ge 0$ for each $j \neq i$, and the vector ${\bf 1}^T
\equiv (1,1, \dots, 1)^T$ of length equal to $n_s$  is in $\cN({K^T})$.

The formal solution of (\ref{kolmo}) is 
$$p(t)=e^{Kt}p(0),$$
but there are well-known difficulties in computing the exponential if the state  space
is large, and various approximate methods have been used
\cite{Kazeev:2014:DSC,Menz:2012:HSD,Deuflhard:2008:ADG}. An alternate approach
is to do direct stochastic simulations, using the Gillespie algorithm or one of
its many variants. However the computational time can be extremely long if there
are many species or the system is spatially distributed \cite{Hu:2013:SAR}, and
thus it is advantageous to reduce the system if possible.

Reduction of deterministic systems is frequently done by taking advantage of the
presence of multiple time scales in the evolution, and reducing the number of
variables by invoking the QSSH, as was done in  \cite{Lee:2009:MTS}.  In the
stochastic system one can determine whether a reaction is fast or not according
to the magnitude of the transition rate of a reaction, the so-called propensity
function in the chemical literature. A larger transition rate implies that the
corresponding step occurs more frequently, but since the transition rate of a
reaction depends on the numbers of reactant molecules as well as the 
rate constant of the reaction, a reaction that is fast in the interior of
${\bf C}^+_m$ may be slow near the boundary of the cone.  The following example
shows that a reaction with a large  rate constant may have to be considered as a
slow reaction.
\begin{Example}
\label{eg:simple}
Consider the following reaction network
 $$ A ~\ssang{k_1}{k_{-1}} ~ B 
$$
 If $k_1= 0.01$, $k_{-1} = 0.1$, and initially $n_A(0) = 100$, $n_B(0) = 1$,
 then the rate of $A \rightarrow B$ is $k_1n_A(0) = 1$, and from $B
 \rightarrow A$ is $k_{-1}n_B(0) = 0.1$. Thus even though $k_1 < k_{-1}$,
 initially $ A \xrightarrow{} B$ is  fast compared with $ B \rightarrow  A$. Of
 course when $A$ is small  $k_1n_A(t)<k_{-1}n_B(t)$, and which is fast and which
 is slow is  interchanged.
\end{Example}
This issue must be addressed on a
case-by-case basis, but here we simply assume that all steps can be identified
as occurring on either an $O(\frac{1}{\ep})$ scale or on an $O(1)$ scale {\it a
priori}. Since we consider a small system with small numbers of molecular
species involved in at most second-order reactions, the separation assumption
can be justified by stipulating that 
$$k_f \approx O(\frac{1}{\ep}) \textrm{ and } N_0^2 k_s \approx O(1),
$$
where $N_0$ is the maximum number of molecules of reactants involved in slow
reactions and $k_f$ and $k_s$ are characteristic rate constants of fast and slow
reactions, respectively. Thus we can rewrite the matrix $K$ as
$$
K= \frac{1}{\ep} K^f + K^s,
$$
where $K^f$ and $K^s$ are the fast and slow transition matrices whose entries
are the transition rates of fast and slow reactions, respectively. Then 
 on the $O(1)$ timescale (\ref{kolmo}) can be written 
\begin{eqnarray}
\frac{d p(t)}{dt}= (\frac{1}{\ep} K^f + K^s) p(t). 
\label{twome}
\end{eqnarray}
As before, conservation of probability must be satisfied, but now it must hold
separately on the slow and fast scales, \ie, ${\bf 1}^T \in \cN(({K^f})^T)~
\textrm {and}~ \cN(({K^s})^T)$. Thus both the transition rate matrices $K^s$ and $K^f$ for slow and fast reactions are
generators of Markov chains. As we see later, the sub-graphs $\cG_s^f$ and
$\cG_s^s$ of $\cG_s$ associated with the fast and slow reactions may have
several disjoint connected components when considered as undirected
graphs. These need not be strong components of $\cG_s$.

Equation (\ref{twome}) is written on a slow time, but on the fast
timescale $\tau = t/\ep$ it reads
\begin{equation}
\frac{d p(t)}{d\tau}= (K^f + \ep K^s) p(t)
\label{slowscal}
\end{equation}
and from this one sees that the slow reactions act as a perturbation of the fast
reactions on this scale. It follows that the perturbation can only have an effect
on the dynamics represented by the eigenvector(s) corresponding to the zero
eigenvalue of $K^f$, and this will be exploited later.

\subsection{Conditions under which the state space is bounded }
\label{Subsec:bound}

A number of technical issues have to be discussed before proceeding
further. Firstly, it is easy to see that the positive cone of $\bZ^{m}$ is
invariant, which therefore preserves non-negativity of the
probabilities. Furthermore, if the reaction simplex is compact there are only
finitely many states in the system, and the generator of the Markov chain is a
finite matrix. In that case the invariant distribution is well-defined and,
given certain additional properties discussed later, it is unique. When the
state space has infinitely many points the generator is an infinite matrix, and
in fact, is unbounded as an operator on $l_{\infty}$\footnote{A vector space whose elements are infinite sequences of real numbers.} because the entries of $K$
depend on $n$. In this case little is known in general about its spectral
properties and invariant measure(s). This is in contrast with random walks on
homogeneous lattices, where the generator is bounded even if the state space is
infinite.  Moreover, the fact that the corresponding deterministic description
has a compact invariant set is of no import in the stochastic description, since
large deviations from the mean dynamics are possible. However, as
the following example shows, the probabilities of very large numbers may in
general be very small under suitable conditions.  

\begin{Example}

Consider the simple process
$\phi
\stackrel{\ds{k_{1}}}{\ds\longrightarrow} A
\stackrel{\ds{k_{2}}}{\ds\longrightarrow} \phi$, where as usual, $\phi$
represents a source and sink, and let $p_{n}(t)$ be the probability of having
$n$ molecules of $A$ in the system at time t. Then one can show that the
generating function $$ G(s,t) = \sum_{n=0}^{\infty} s^{n}p_{n}(t) $$ satisfies
$$
\pdt{G} +k_{2}(s-1) \pds{G} = k_{1}(s-1)G.
$$
The solution of this for the initial condition $p_{1}(0)= 1, p_{n}(0)=0$
otherwise, is 
$$
G(s,t) = \left(1 +
(s-1)e^{-k_{2}t}\right)\cdot exp\left(\dfrac{k_{1}}{k_{2}}(s-1)(1-e^{-k_{2}t})\right)
$$
and from this one finds, by expanding this as a series in $s$,  that 
\begin{equation}
p_{n}(t) =
\dfrac{1}{k_{2}n!}\left(\dfrac{k_{1}}{k_{2}}\right)^{n-1}\left(1-e^{-k_{2}t}\right)^{n-1}\left(k_{1}(1-e^{-k_{2}t})^{2}
+ k_{2} n e^{-k_{2}t}\right) \cdot exp\left(-\dfrac{k_{1}}{k_{2}}(1-e^{-k_{2}t})\right)
\end{equation}
Therefore the stationary distribution is
$$ \ds \lim_{t \rightarrow \infty} p_{n}(t)=
\dfrac{1}{n!}\left(\dfrac{k_{1}}{k_{2}}\right)^{n}exp\left(-\dfrac{k_{1}}{k_{2}}\right),
$$ 
 which is a Poisson distribution with parameter $k_{1}/k_{2}$. Thus $p_{n}(t)$
is non-zero for arbitrarily large $n$ in both the transient and stationary
distributions, but it decays rapidly with $n$. For example, if $k_{1}/k_{2} \sim
O(1)$ and $n = 25$, $p_{n} \sim O(10^{-25})$ in the stationary
distribution. Even if the stationary mean $k_{1}/k_{2} \sim O(10)$, $p_{n} \le
10^{-20}$ for $ n \ge \sim 50$ (one must always choose $n$ greater than the mean
in order that $p_{k} < p_{n}$ for $k >n$). 

\end{Example}

While this is only a heuristic
justification for truncating the state space to a box of the form $\Pi_{i =
1}^{m}[0,M_{i}],$ we shall do this in the remainder of the paper. One could
formalize the truncation by modifying the dynamics so as to make all states that
lie in the positive cone and outside the box transient, and then restricting the
initial data to lie in the box. Further work on the validity of this truncation
is needed, but in any case the following results apply rigorously to closed
systems and those systems whose inputs are terminated when the total number of the
molecules being input reaches a pre-determined level, as in the following
example.

\begin{Example}
  Consider a triangular system in which there is a input of one species, as shown
  below.  
  $$
\phi\xrightarrow{k_0}A\xrightleftharpoons[k_2]{k_1} B\xrightleftharpoons[k_4]{k_3}C\xrightleftharpoons[k_6]{k_5}A.
$$
If the initial state is $(n_A(0), n_B(0),
  n_C(0))=(0,0,0)$ this defines the initial simplex by $n_A(0) + n_B(0) + n_C(0)$. Whenever a molecule of
  $A$ is added the simplex shifts upward by one unit in  the first octant. To
  insure  that the simplex remains bounded we simply turn off the input  $\phi$ when the 
   number of species A achieves a certain threshold, say, $N_A$. At 
  this moment $t=t^*$, we also observe the molecular number of species B and C, 
  which are $n_B(t^*), n_C(t^*)$, and the system evolves on this 
  reaction simplex thereafter, since the  system is closed thereafter.
    
\end{Example}
 
In general we shall assume that the inputs and outputs occur in the slow
dynamics, and thus the fast reactions occur on a fixed finite simplex. The
inputs and outputs can move the system `up' and `down' in the positive cone, but
we assume that this occurs on the slow time scale and if the total inputs are controlled the
state space will remain bounded.

\subsection{The role of invariants}
 \label{subsec:inv}

Any conservation condition that applies to a system defined in terms of
concentrations applies when the system is defined in terms of numbers, and
therefore the invariants that characterize (\ref{dd}) apply in the stochastic
framework as well. In particular, the definition of a kinetic manifold
$\Omega_f$ for the fast subsystem, in which the slow variable $\tilde c$ is
invariant, carries over with only minor modification.  Since the state space is
discrete in a  stochastic description, only the set of integer points in a reaction
simplex $\Omega$ are relevant, and we call the subset of integer points in $\Omega$
the  discrete reaction simplex and denote it by~$\cL$. 

For any $c\in \Omega(c_0)$, we can write $c=c_0+\nu\cE \eta \ge 0$, where $\eta
$ is the extent, and therefore
$$
n= \cN_A V c= \cN_A V c_0 + \cN_A V \nu\cE  \eta \in \Omega( \cN_A V c_0) =:
\Omega(n_0).
$$ 
Since $n$ lies in $ \overline{{\bf Z}^+_m}$, which is the closure of the set
${\bf Z}^+_m$ of $m$-dimensional vectors with positive integer entries, it
follows that the set of all integer points in the reaction simplex is given
$$\cL(n_0) \equiv\Omega(n_0) \cap  \overline{{\bf  Z}^+_m},$$
and we call this the (full) discrete reaction simplex through $n_0$.
Analytically, a discrete simplex $\cL(n_0)$ is defined as a coset of
$\cR(\nu\cE)$, and it can be generated numerically using
Algorithm~\ref{alg}. To illustrate this, consider the following reaction network
for  an open system. 

\begin{figure}[h!]
\begin{center}
$\begin{array}{cc}
{\includegraphics[scale = 0.3]{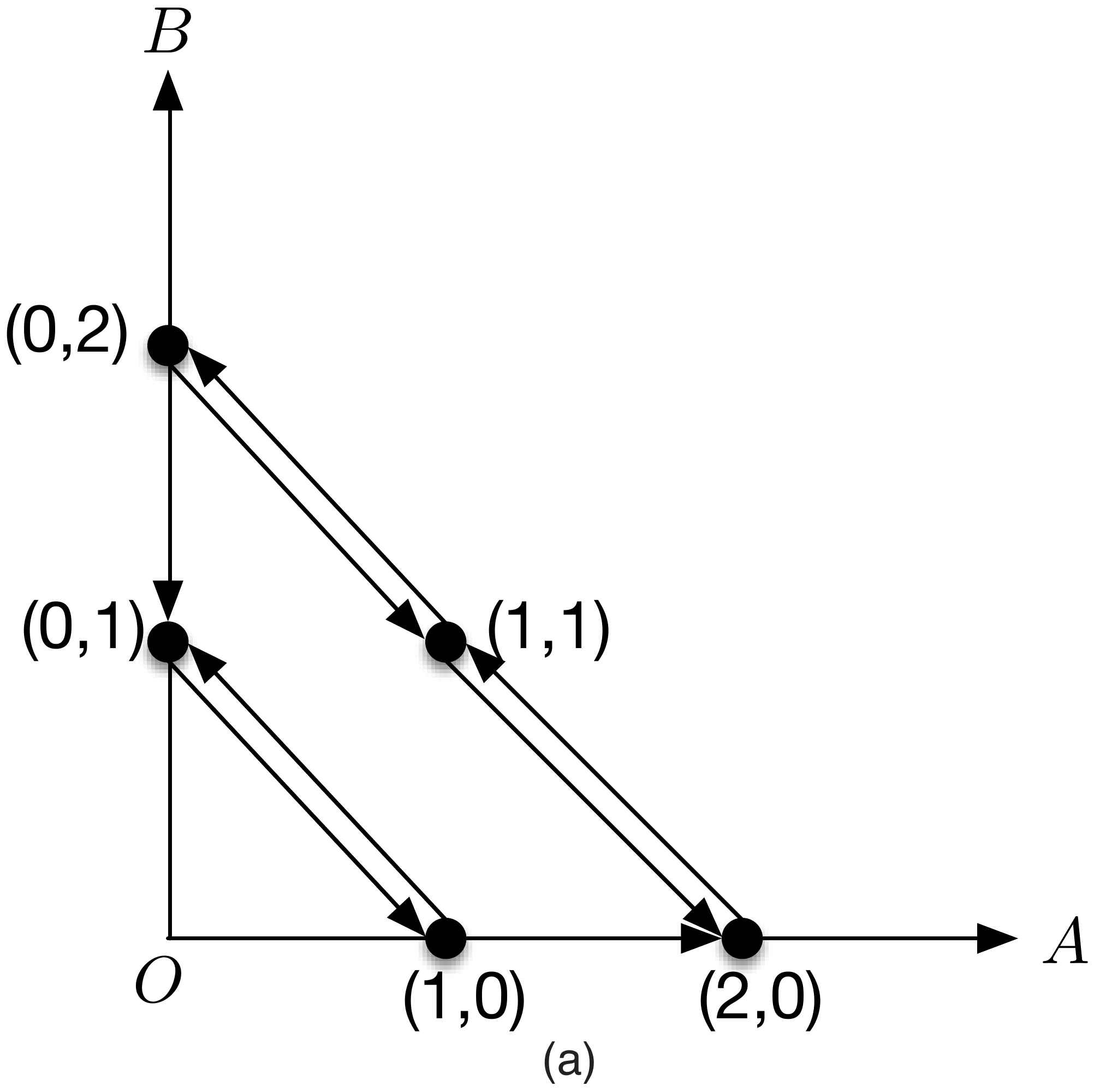}}&\hspace*{20pt}{\includegraphics[scale = 0.3]{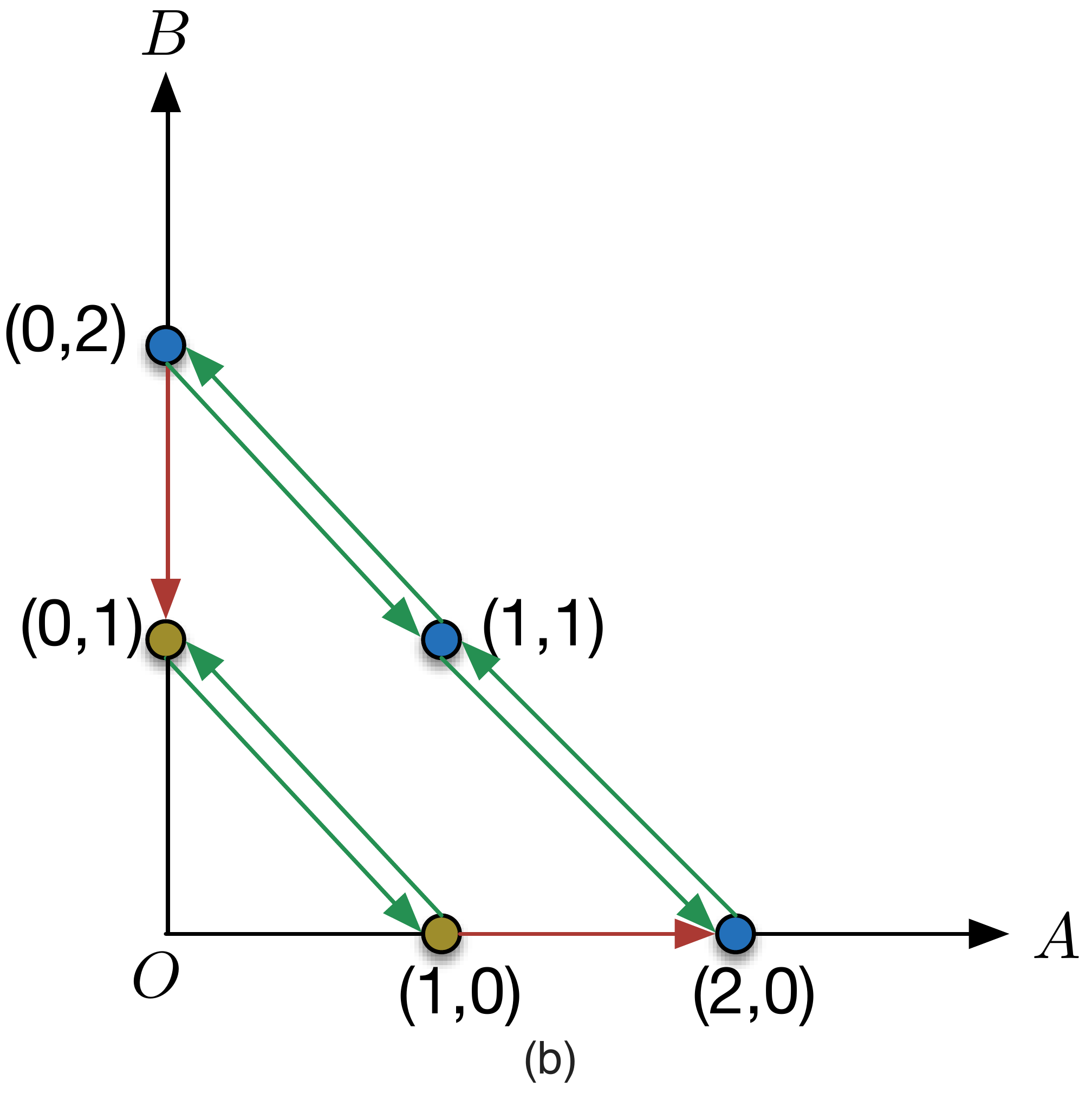}  }
\end{array}
$
\caption {(a) The full  discrete reaction simplex; (b) two fast discrete
reaction simplexes distinguished by yellow and blue colors. }
\label{fig:simplexDef}
\end{center}
\end{figure}

\begin{Example}
\label{eg:simplex}
Consider the network 
$$
\phi \rightarrow A\rightleftharpoons B\rightarrow\phi.
$$
where the input and output reactions are slow steps. One finds that 
$$
\cR(\nu\cE) \cap  \overline{{\bf Z}^+_2}= span\left\{\begin{pmatrix}
-1\\ ~0 \end{pmatrix} \begin{pmatrix} 0\\ 1 \end{pmatrix}  \right\}
\cap  \overline{{\bf    Z}^+_2} 
$$
which is $\overline{{\bf Z}^+_2}$. If we start with the initial state $(1,0)$,
and  stop the input reaction once one molecule of $A$ is added from
the source,  as described in Section~\ref{Subsec:bound}, we can generate the full
discrete reaction simplex as the set of all points denoted `$\bullet$'
Figure~\ref{fig:simplexDef}(a). 

   Since the inputs and outputs occur on the slow time scale, these reactions are denoted by red arrows in
   Figure 3.1 (b), and the remaining reactions are in fast components (green
   arrows).  There are two distinct fast (absorbing) strong components, as shown in
   Figure~\ref{fig:simplexDef}(b),  defined as 
$\cL_1^f =   \{(1,0),(0,1)\}$ and $\cL_2^f = \{(2,0),(1,1),(0,2)\}$. 
Since
   $\cR(\nu\cE^f)=span\{[-1\quad1]^T\}$, $\cN[(\nu\cE^f)^T]=span\{[1\quad 1]\}$,
   we can identify each fast simplex uniquely by
$$
\tilde n = A^fn = [1\quad 1]\left[ \begin{array}{c} n_1 \\
n_2
\end{array}\right] = [n_1+n_2],
$$
where $n_1 (n_2) $ denotes the number of molecules of species $A$ ($B$). The
rows of the matrix $A^f$ comprise a basis for $\cN[(\nu\cE^f)^T]$, and thus
$\cL_1^f$ is defined by $ \tilde n = 1$, and $\cL_2^f$ by $\tilde n = 2$. In
this simple example, we can interpret $\tilde n$ as the total mass,
which is conserved by the fast reactions.
\end{Example}

To define fast simplexes in general we consider only fast reactions, and first identify all
distinct  absorbing strong components in $\cG_s^f$. Then, working
backwards, we identify the sources and internal strong components that feed into
an absorbing component. This identifies a sub-graph of $\cG_s$ comprised of an
absorbing strong component and its `feeders', and leads to the following
definition of a fast simplex.
\begin{Definition}
\label{fastsimplex}
 A fast discrete reaction simplex $\cL^f$ is the set of vertices in an  absorbing strong component and
 its predecessor sources and internal strong components. 
 \end{Definition}
\begin{Remark}
\mbox{}
\
\begin{enumerate} 
\item A fast simplex is uniquely determined by its  absorbing strong
  component. Each fast simplex is characterized by the invariants of the fast
  dynamics. These are 
\begin{equation}
\tilde n = A^f n
\label{fasta}
\end{equation}
where the rows of the matrix $A^f$ comprise a basis for $\cN[(\nu\cE^f)^T]$.

\item The invariant distribution of the fast dynamics on a given fast simplex is
  non-zero  only for the  states in the corresponding absorbing component. 

\item Both source components and internal strong components may belong to more
  than one fast simplex, as is exemplified by a reaction network with a binary  tree
  structure. In such a network the root belongs to all fast simplexes and
  internal nodes at the $k^{th}$ level belong to  $2^{n-k}$ simplexes, where $n$
  is the depth of the tree.

\item If there are no slow reactions in the network there is only one fast
  simplex, and this case is not of interest here. 

\item  In general the disjoint union of all fast simplexes is the full fast simplex, i.e.
$$
\dot \cup\cL^f \equiv \Omega^f(n_0) \cap \overline{{\bf Z}^+_m}.
$$

\end{enumerate}
\end{Remark}

Although the following reaction network due to Wilhelm (\cite{Wilhelm:2009:SCR})
involves unrealistic trimolecular steps, we use it to illustrate how to define
the fast simplexes when the fast subgraph $\cG_s^f$ has multiple absorbing
components whose predecessors  are not disjoint.

\begin{Example}
\label{eg:MAK}
\
\begin{center}
{\includegraphics[scale=0.6]{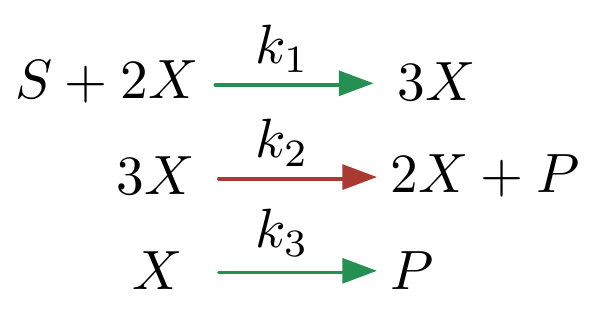}  }
\end{center}

 The species are  S, X, P, and the complexes: S+2X, 3X, 2X+P, X, P are labeled
 1-5 in this order.   Then the matrices $\nu$, $\cE$ and $\nu\cE$ are as follows. 
  $$
  \nu = \left[
  \begin{array}{ccccc}
     1 & 0 & 0 & 0& 0\\
    2 & 3 & 2 & 1 & 0\\
    0 & 0 & 1 & 0 & 1
  \end{array}  \right]
  \quad
  \cE = \left[ 
\begin{array}{ccc}
    -1 & 0 & 0\\
    1 & -1 & 0\\
    0 & 1 & 0\\
    0 & 0 & -1\\
    0 & 0 & 1
  \end{array}
\right]
  \quad
  \nu\cE = \left[
  \begin{array}{ccc}
     -1 & 0 & 0\\
    1 & -1 & -1\\
     0 & 1 & 1
  \end{array}
\right].
  $$

  It follows that  $\rho(\cE) = 3, \rho(\nu\cE) = 2$, and the deficiency, which
  is defined as  $\delta =  \rho(\cE) - \rho(\nu\cE) =1$. Clearly $S + X + P$ is
  a constant and equivalently, $\cN[(\nu\cE)^T] = span\left( ( 1,1,1)^T\right)$. 
  If we set the initial state  as $(n_S(0), n_X(0), n_P(0)) = (1, 2, 0)$, then 
  the possible transitions in the system are  as shown in Figure~\ref{fig:simplexDef2}(a).

Assume that only the second reaction is slow and the others are fast.
Considering only the fast transitions indicated by green directed edges, it is
clear that the resulting $\cG^f_s$ has the two absorbing strong components
$(1,0,2)$ and $(0,0,3)$ because there are two possible exit steps from the node
$(1,2,0)$, one with rate $k_1$, the other with rate $k_3$. The two
fast simplexes are
\begin{align*}
\cL_1^f=&\{(1,2,0),(1,1,1),(1,0,2)\}\\
\cL_2^f=&\{(1,2,0),(0,3,0),(0,2,1),(0,1,2),(0,0,3)\}
\end{align*}
 as shown in Figure~\ref{fig:simplexDef2}(b) \& (c). Figuratively, and also in
 the language of graph theory, the simplexes are branches of the tree that
 defines the full simplex in Figure~\ref{fig:simplexDef2}(a). $\cL_1^f$ lies on
 the coset of $\cR(\nu\cE^{f_1}) = span\{[0\quad-1\quad1]^T\}$, thus
 $\cN[(\nu\cE^{f_1})^T] = span\{[1\quad0\quad0],[0\quad1\quad1]\}$, and
 $\cL_1^f$ is identified by $\tilde n_1 = (1,2)$. $\cL_2^f$ lies on the coset of
 $\cR(\nu\cE^{f_2}) = span\{[-1\quad1\quad0]^T, [0\quad-1\quad1]^T\}$, thus
 $\cN[(\nu\cE^{f_2})^T] = span\{[1\quad1\quad1]\}$, and $\cL_2^f$ is identified
 by $\tilde n_2 = 3$.

\begin{figure}[h!]
\begin{center}
$
\begin{array}{ccc}
{\includegraphics[scale = 0.3]{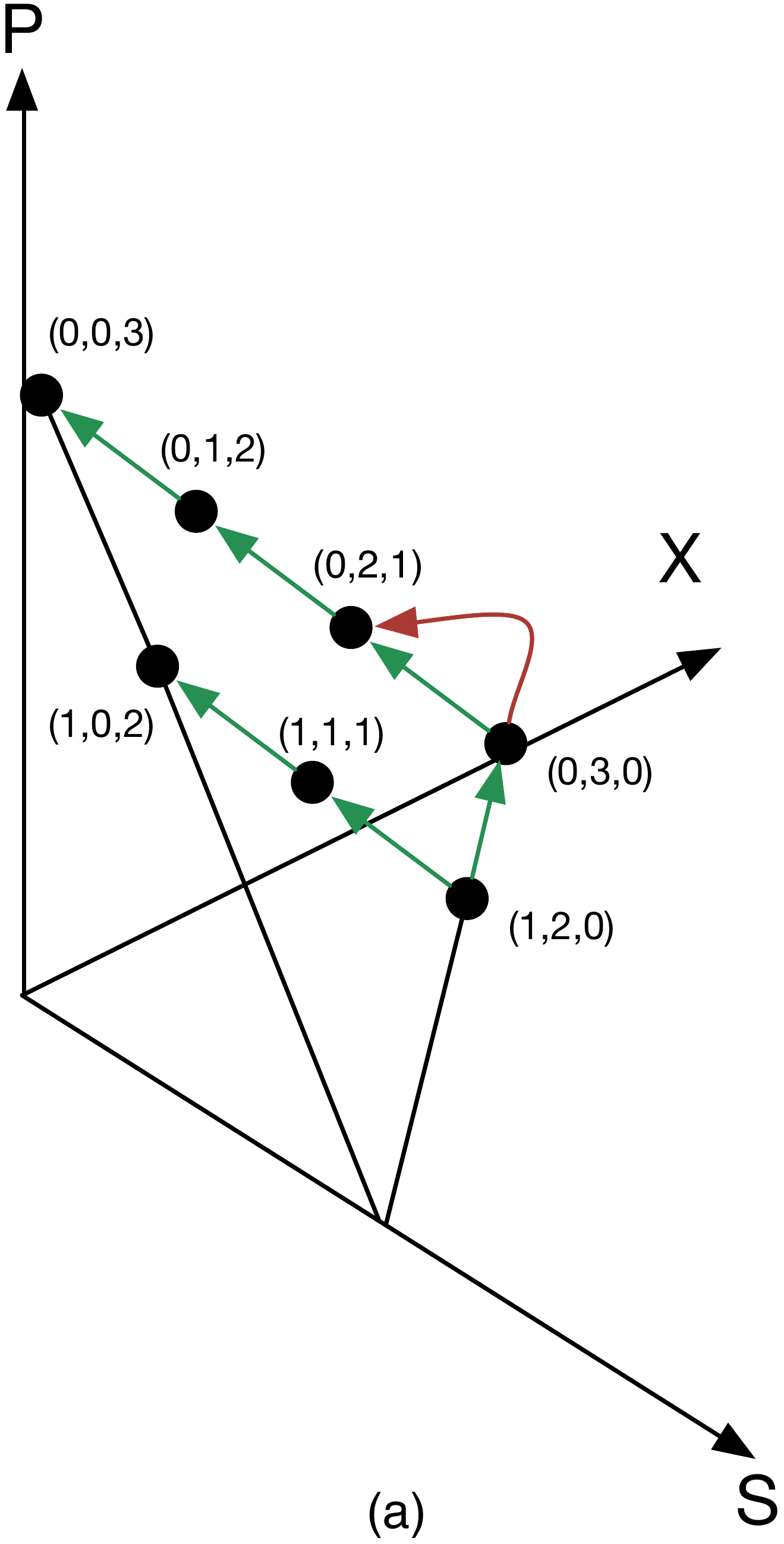}}&{\includegraphics[scale =    0.3]{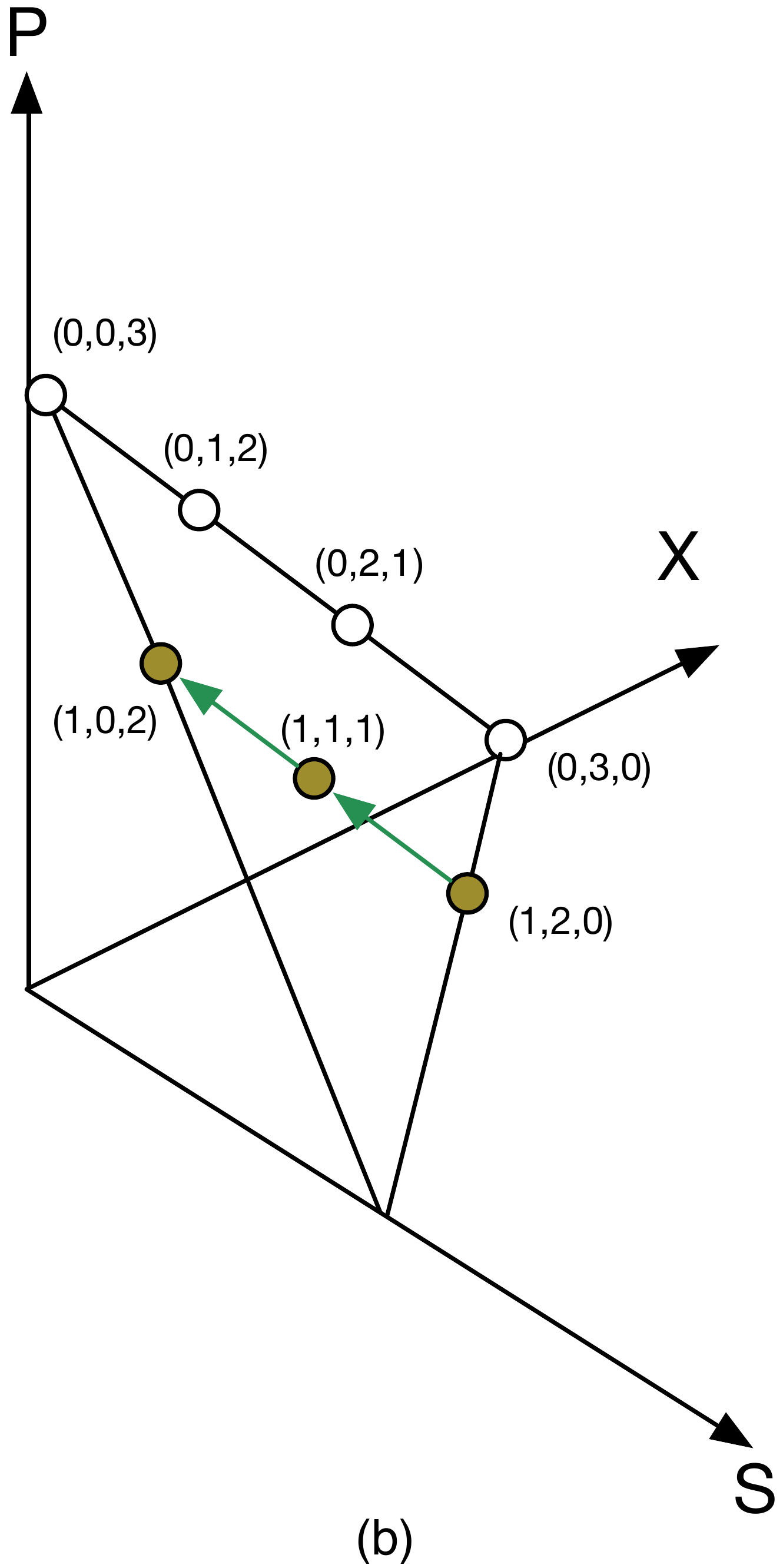}} &{\includegraphics[scale =0.3]{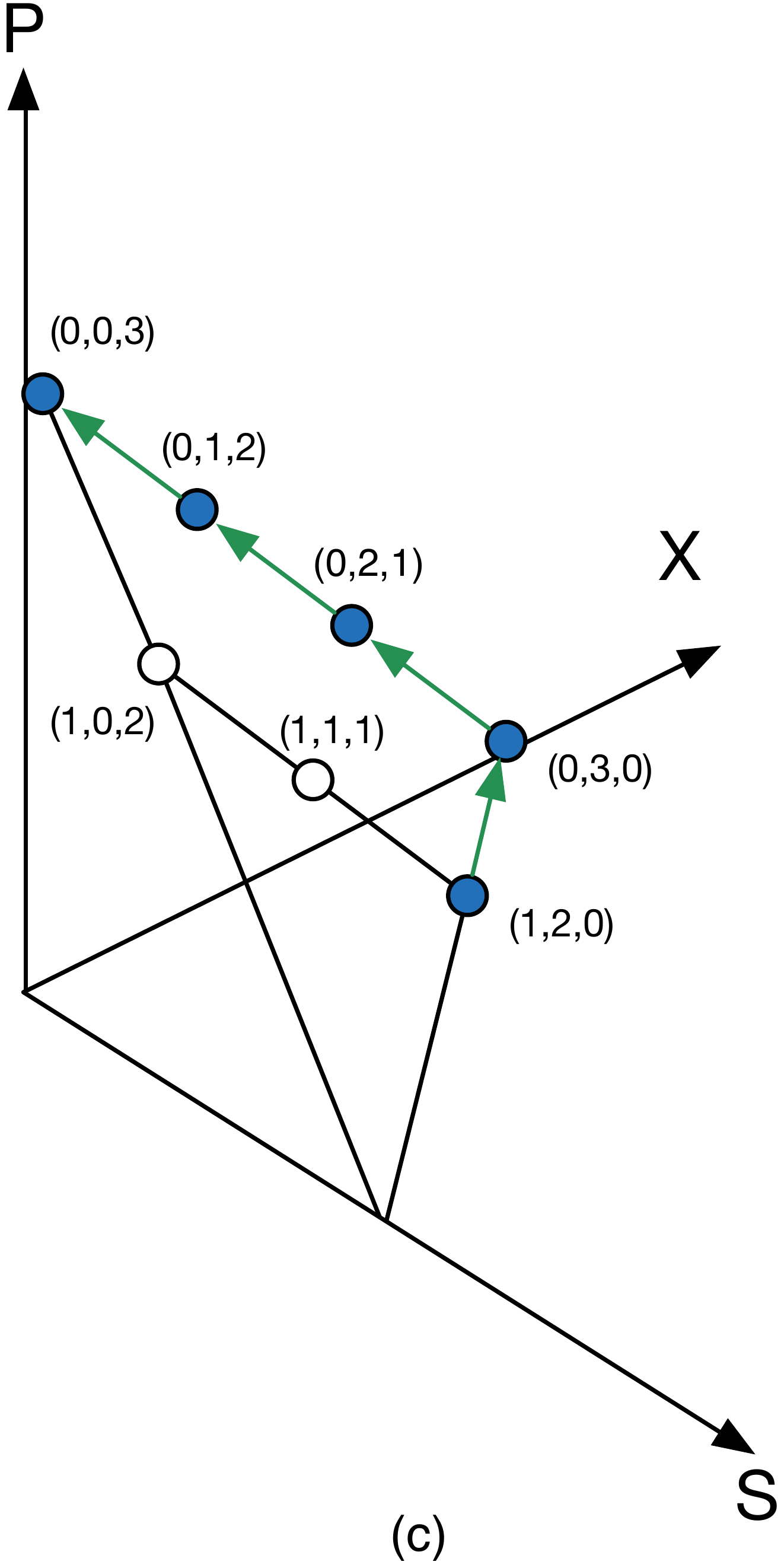}}                                  
\end{array}
$
\caption{\footnotesize  Example~\ref{eg:MAK}: (a) The full discrete reaction  simplex(black circles); (b) \& (c) the fast simplexes. }
\label{fig:simplexDef2}
\end{center}
\end{figure}

\end{Example}
  In general a component of the reaction graph may have more than one absorbing
  component and more than one fast simplex, as in the foregoing example. To
  determine them in general one must first identify the number of strong
  components and their type in every component graph $\cG_{\alpha} \subset
  \cG$ .  Since no vertex in a source is reachable from any vertex
  outside its component and no vertex in a sink is reachable from a vertex in
  any other sink, the relationship of reachability defines a partial order on
  the strong components of $\cG_{\alpha}$.  This in turn leads to the acyclic
  skeleton $\stackrel{\circ}{\cG}_{\alpha}$ of $\cG_{\alpha}$, which is defined
  as follows.  Associate a vertex $\stackrel{\circ}{V_j}$ with each strong
  component of $\cG_{\alpha}$, and introduce a directed edge from
  $\stackrel{\circ}{V_i}$ to $\stackrel{\circ}{V_j}$ if and only if one (and
  hence every) vertex in $V_{\alpha j}$ is reachable from $V_{\alpha i}$.
  $\stackrel{\circ}{\cG}_{\alpha}$ is connected since $\cG_{\alpha}$ is
  connected, but it is acyclic; in fact, it is a directed tree.  One of two
  cases obtains: either $\stackrel{\circ}{\cG}_{\alpha}$ consists of a single
  vertex and no edges, which occurs when $\cG_{\alpha}$ consists of one strong
  component, or it has at least one vertex of in-degree zero and at least one
  vertex of out-degree zero.  One can relabel the strong components if
  necessary, and write the adjacency matrix of $\stackrel{\circ}{\cG}_{\alpha}$
  in the form

  \begin{displaymath}
   \cA = \left[ \begin{array}{ll|ll|ll}
0 &0 &0 &0 &0 &0 \\ 0 &0 &0 &0 &0 &0 \\ \hline
x &x &0 &0 &0 &0 \\ x &x &x &0 &0 &0 \\ \hline
x &x &x &x &0 &0 \\ x &x &x &x &0 &0 \end{array} \right].
  \end{displaymath}
 where the $x$'s represent blocks that may be non-zero.  The diagonal blocks are
 square matrices of dimensions equal to the number of sources, the number of
 internal strong components, and the number of sinks, respectively.  The
 vertices corresponding to internal strong components can always be ordered so
 that the central block is lower triangular, since the strong components are
 maximal with respect to inclusion of edges.  The number of sinks or absorbing
 strong components is the number of zero columns of $\cA$.

\section{The reduction of  a master equation with  two time scales}
\label{master_reduc}

\subsection{The splitting of the evolution equations} 

In this section we show how to obtain the lowest-order  approximation to the
slow dynamics for a system that is described by    
\begin{eqnarray}
\frac{dp}{dt}= (\frac{1}{\ep} K^f + K^s) p. 
\label{evoslow}
\end{eqnarray}
on the $O(1)$ timescale and by 
\begin{eqnarray}
\frac{dp}{d\tau}= (K^f + \ep K^s) p 
\label{evofast}
\end{eqnarray}
on the fast time scale $\tau = t/\ep$.

By a suitable ordering of the states  we can write the  fast transition 
matrix $K^f$ as a block diagonal matrix 
\begin{eqnarray}
K^f=\left[ \begin{array}{cccccc}  K^f_1 & &  &  &  &   \\   & K^f_2 &  &  &  &
\\  &  & K^f_3 &  &   & \\ & & & \ddots &  &  \\ & & & & K^f_{l-1}  &  \\ & & &
&  & K^f_{l} \end{array} \right]
\label{generalblock}
\end{eqnarray}
wherein the number of blocks is equal to the number of fast components in the
state transition graph $\cG_s^f$ of the fast reactions.  We show later in Theorem~\ref{thm:semisimple} that each
block has zero as a semisimple eigenvalue and we analyze the structure of the
blocks in detail, but here we first define the index of a matrix
\cite{Campbell:1991:GIL}.

\begin{Definition}
Let $A: {\bf R}_n \rightarrow {\bf R}_n$ be a linear transformation. The  index of $A$ is
the smallest non-negative integer $k$ such that $\cR(A^{k+1}) = \cR(A^{k})$. 
\end{Definition} 
An equivalent definition is that the index of $A$ is the largest Jordan block
corresponding to the zero eigenvalue of $A$.  It follows from the general theory
in \cite{Campbell:1991:GIL} that a matrix with a semisimple zero eigenvalue has
index 1, and as a result, the following properties hold.

\begin{itemize}
 \item [(i)] ${\bf R}_{n} = \cN(A) \oplus \cR(A)$~~\footnote{This is the direct sum, but is
       generally not the orthogonal direct sum.}

\item [(ii)] $\cR(A^{2}) = \cR(A)$

\end{itemize}

Since we assume that the inputs and outputs only occur on the slow time scale,
the fast dynamics occur on a fixed compact simplex in ${\bf C}^{+}_{m}$ and $K^f$ is
a bounded operator that generates a Markov chain, since ${\bf 1} \in
\cN(K^f)^T$.  We define $\dim(\cN(K^f))=n_f$, and of course it follows that both $\cN(K^{f})$ and
$\cN((K^{f})^{T})$ have a basis of $n_{f}$ linearly-independent eigenvectors.
Let $\Pi$ be the $n_{s}\times n_f$ matrix whose columns are eigenvectors
corresponding to zero eigenvalues of $K^f$, and let $L$ be the $n_f\times n_{s}$
matrix whose rows are eigenvectors corresponding to zero eigenvalues of
$(K^f)^T$.  One can choose these to form a biorthogonal set and therefore $L\Pi
= I_{n_f}$. Furthermore, $\Pi L = P_0,$ where $P_0$ is the eigenprojection
corresponding to the zero eigenvalue of $K^f$, and thus $\cR(P_0)=\cN(K^f)$ and
${\bf R}_{n_s}=\cR(P_0)\oplus \cR(I-P_0)$. 

By property (i)  above ${\bf R}_{n_s}=\cN(K^f)\oplus \cR(K^f)$ and to show that
$\cR(I-P_0)=\cR(K^f)$ we observe that the adjoint $P_{0}^{T}$ leads to the
decomposition ${\bf R}_{n_s}=\cR(P_0{^T})\oplus \cR(I-P_0^{T})$. We have that
$\cR(P_{0}^{T}) \equiv \cN((K^{f})^{T})$ is orthogonal to $\cR(I-P_{0})$, and
from the basic properties of projections and their adjoints \cite{Kato:1966:PTL}
it follows that $\cR(I-P_0)=\cR(K^f)$.

 Since  
$\mathbf{R}_{n_{s}}= \cN(K^f)  \oplus \cR(K^f),$ we can write 
\begin{eqnarray}
p = P_{0}p + (I-P_{0})p = \Pi \tilde p + \Gamma \hat p \label{sepa}. 
\end{eqnarray}
where $\tilde p = L p$ and $\Gamma$ is an $n_{s}\times (n_{s}-n_f)$ matrix whose
columns are basis vectors of $\cR(K^f)$. Then the matrix $[\Pi \ | \ \Gamma]$ is
an $n_{s}\times n_{s}$ matrix whose columns are basis vectors of ${\bf R}_{n_s}$. 
Using (\ref{sepa}) in (\ref{evoslow}) leads to
\begin{eqnarray}
\Pi \frac{d\tilde p}{dt} + \Gamma \frac{d\hat p}{dt} &=& (\frac{1}{\ep} K^f+
K^s)(\Pi \tilde p + \Gamma \hat p). \label{sepa1}
\end{eqnarray}

To obtain the evolution equations for $\tilde p$ and $\hat p$ separately, we define the
$(n_{s}-n_f)\times n_{s}$ matrix $V = [0 \ | \ I_{n_{s}-n_f}][\Pi \ | \ \Gamma
]^{-1}$, which has the property that $V [\Pi \ | \ \Gamma ] = [0 \ | \
I_{m-n_f}].$ In other words, $V\Pi = 0,\ V\Gamma =I_{n_{s}-n_f}$, and by multiplying
(\ref{sepa1}) by $L$ and $V$, respectively, we obtain
\begin{eqnarray}
\frac{d\tilde p}{dt} &=&   L K^s (\Pi\tilde p + \Gamma \hat p),  \label{slowp}\\
\frac{d\hat p}{dt} &=& V ( \frac{1}{\ep} K^f +  K^s) (\Pi\tilde p + \Gamma \hat p) =\frac{1}{\ep} V K^f \Gamma \hat p + V K^s(\Pi \tilde p +\Gamma \hat p). \label{fastp}
\end{eqnarray}
where we used the fact that $L\Gamma = 0$ because $\cN((K^{f})^{T})) \perp
\cR(K^{f})$.  

On the fast time scale  $\tau=t/\ep$  we obtain in the limit $\ep
\rightarrow 0$ that  
\begin{align}
\frac{d\tilde p}{d\tau}&=  0, \\
\frac{d\hat p}{d\tau} &=  V K^f \Gamma \hat p. 
\label{hatp_fast}
 \end{align}
On this time scale the slow variable $\tilde p$ remains constant
while the fast variable $\hat p$ evolves according to  (\ref{hatp_fast}).

On the slow time scale 
\begin{eqnarray}
\frac{d \tilde p}{dt} &=& L K^s (\Pi \tilde p+\Gamma \hat p), \label{slowp1}\\
\ep \frac{d\hat p}{dt} &=&  (V K^f + \ep V K^s) (\Pi\tilde p + \Gamma \hat p).
\label{fastp1} 
\end{eqnarray}
and  the limit $\ep \to 0$ leads to the condition  that 
$$ 
 V K^f \Gamma \hat p = 0.
$$

\begin{Lemma} 
$$V K^f \Gamma \hat p =0$$ 
only if ~  $\Gamma \hat p =0$.
\end{Lemma}

\begin{proof}
By definition $\Gamma \hat p \in \cR(K^{f})$, and since $\cR((K^{f})^{2}) =
\cR(K^{f})$ by property (ii) above, it follows that $K^{f}\Gamma\hat p = \Gamma
q$ for some $q$. Since $V \Gamma = I_{n_{s}-n_{f}}$, $V K^f \Gamma \hat p =0$
only if $q =0$, which implies that  $\Gamma \hat p =0$.
\end{proof}

Since $\Gamma \hat p =0$, it follows from the equation (\ref{slowp1}) that
\begin{eqnarray}
 \frac{d \tilde p}{dt} &=& L K^s \Pi \tilde p \equiv \tilde K \tilde p, 
\label{reduced}
\end{eqnarray}
which shows how  the slow evolution between fast simplexes depends on components
of the fast evolution through $L$ and $\Pi$. The components of  $\tilde p$ represent the
probabilities of aggregated states, and since  $\Gamma \hat p =0$ it follows
from (\ref{decomp}) that $p = {\tilde \Pi}\tilde p$,  from which one obtains the 
probabilities of the individual states. 

\medskip
The following  example  illustrates
the structure underlying the general reduction.
\begin{Example} 
\label{Triex}
Consider a closed triangular reaction system with all
  reversible reactions, in which reactions 1 and 2 are fast, as shown
  in Figure \ref{Tri}. The total number of molecules is conserved in a
  closed system of first order reactions of this type, and the
  transitions in state space are as illustrated in Figure~\ref{Tri}
  when the total number of molecules is two .
  
  \begin{figure}[h!]
\begin{center}
$\begin{array}{cc}
{\includegraphics[scale=0.15]{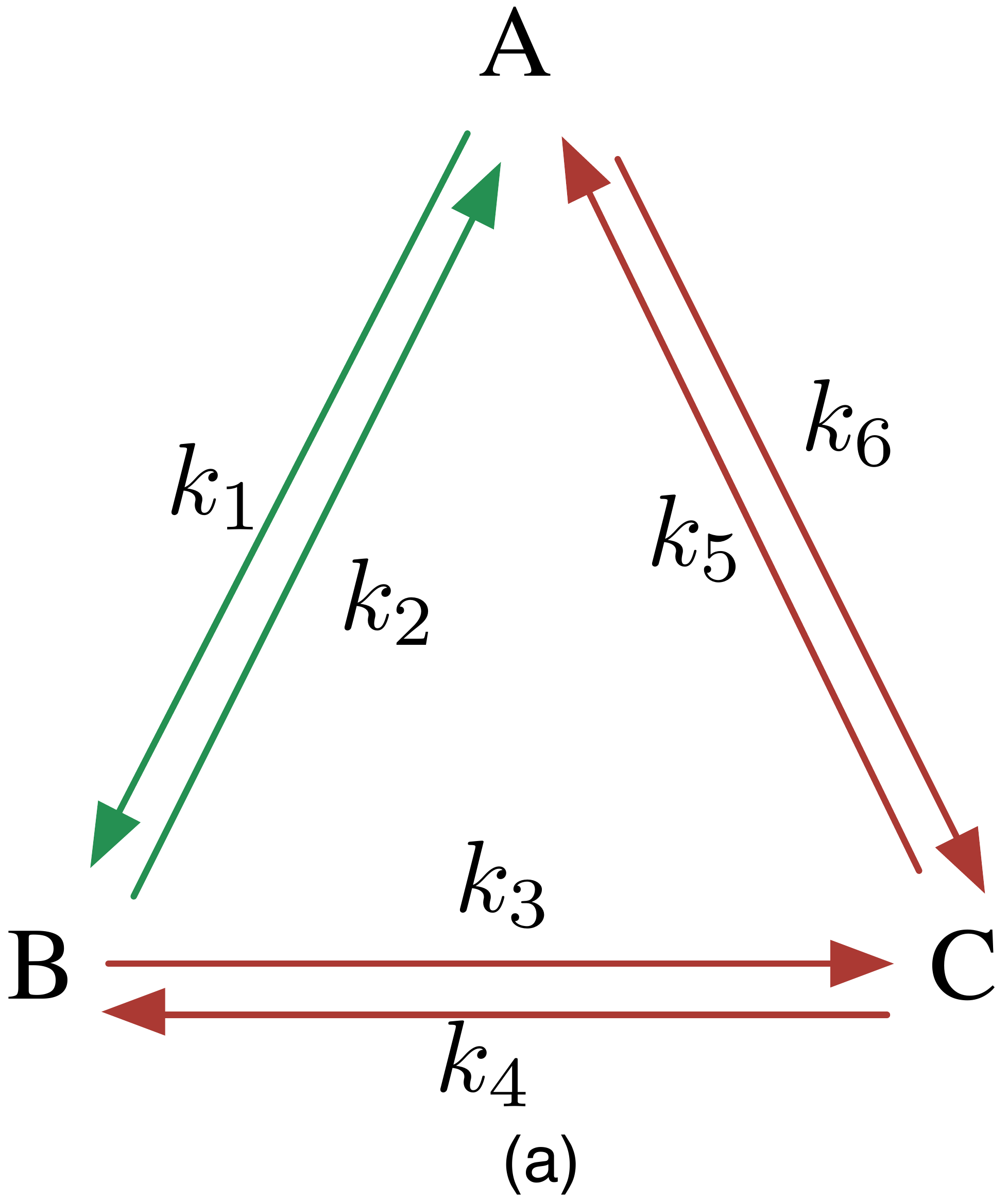}} & \hspace*{.2in}{\includegraphics[scale=0.2]{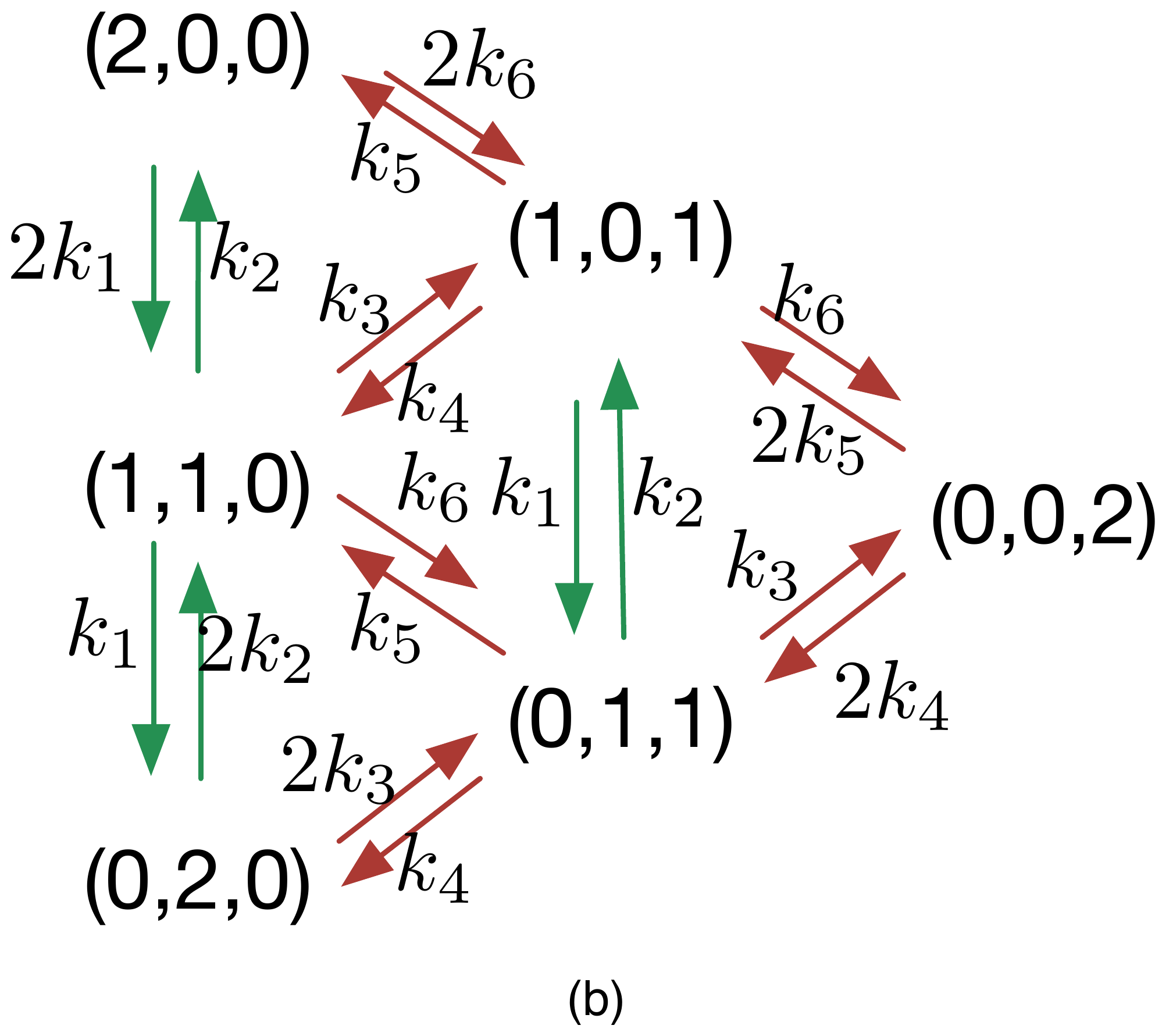}}\\
\end{array}$
 \caption{Triangular reaction: red and green arrows denote  slow and fast
 reactions, respectively. (a) The reaction network; (b) The state transition diagram
 for a total of two molecules. \label{Tri} }
\end{center}
\end{figure}       
  
The matrix $K^f$ is obtained by setting the rates of all slow steps to zero. The
resulting graph has three connected components, 
$(2,0,0) \ssang{}{} (1,1,0)\ssang{}{} (0,2,0)$,~~ $(1,0,1)
\ssang{}{} (0,1,1)$~ and $~~(0,0,2)$, and this leads to the representation 

\begin{equation}
K=\frac{1}{\epsilon}\left[\begin{array}{ccc}  K_1^f &  0 &  0     \\
                                 0 & K_2^f &  0     \\
                                 0   &   0  &  K_3^f   
                                \end{array}\right] + \left[\begin{array}{ccc}  K_1^s &  K_{1,2}^s &  K_{1,3}^s     \\
                                 K_{2,1}^s & K_2^s &  K_{2,3}^s     \\
                                 K_{3,1}^s   &   K_{3,2}^s  &  K_3^s   
                                \end{array}\right],
\label{examp-5}
\end{equation}
where $K_i^f, i = 1, 2, 3$ is the $m_i\times m_i$ matrix of transition rates
 within each fast component, and $K_{i,j}^s$ is an $m_i\times m_j$
 matrix of slow transition rates between fast components. In this
 example, the $K_i^f$ are 
$$
K_1^f=\left[\begin{array}{ccc}  -2k_1 &  k_2 &  0     \\
                                 2k_1 & -(k_1+k_2) &  2k_2     \\
                                 0   &   k_1  &  -2k_2   
                                \end{array}\right], ~~~
                K_2^f=\left[\begin{array}{ccc}  -k_1 &  k_2    \\
                                 k_1 & -k_2   
                                \end{array}\right],~~~
                                K_3^f = 0,
                                $$

Note that the vector $\bm 1$ of the appropriate dimension is a left
eigenvector of each  $K_i^f$, and therefore probability is conserved
on each fast component. In view of the block structure of $K^f$ in
(\ref{examp-5}) zero is a semisimple eigenvalue of $K^f$. 

The sub-matrices of $K^s$ are 

$$
  K_1^s = \left[\begin{array}{ccc}  -2k_6 &  0 &  0     \\
                                 0 & -(k_3+k_6) &  0    \\
                                 0   &  0  &  -2k_3   
                                \end{array}\right],
K_2^s = \left[\begin{array}{cc}  -(k_4+k_5+k_6) &  0      \\
                                 0 & -(k_3+k_4+k_5)  
                                  \end{array}\right],
 K_3^s = -(2k_4+2k_5)
$$
and
$$
K_{2,1}^s=\left[\begin{array}{ccc}  2k_6 &  k_3 & 0     \\
                                 0 & k_6 & 2k_3       
                                \end{array}\right],~~~
  K_{1,2}^s=\left[\begin{array}{ccc}  k_5 & 0     \\
                                 k_4 & k_5  \\
                                 0 & k_4       
                                \end{array}\right]~~~
  $$
$$
K_{3,1}^s=\left[\begin{array}{ccc}  0 &  0 & 0  \end{array}\right], ~~~
K_{1,3}^s=\left[\begin{array}{ccc}  0 &  0 & 0    \end{array}\right]^T \quad
K_{3,2}^s=\left[\begin{array}{cc}  k_6 &  k_3  \end{array}\right],~~~ K_{2,3}^s=\left[\begin{array}{cc}  2k_5 &  2k_4 
  \end{array}\right]^T.
$$
 The matrices $L$ and $\Pi$ are given by the following. 
 $$
L=\left[\begin{array}{ccc}  L_1 &  0 &  0     \\
                                 0 & L_2 &  0     \\
                                 0   &   0  &  L_3   
                                \end{array}\right], 
    \qquad  \Pi=\left[\begin{array}{ccc}  \Pi_1 &  0 &  0     \\
                                 0 & \Pi_2 &  0     \\
                                 0   &   0  &  \Pi_3   
                                \end{array}\right],                     
                                $$
  where $L_1 = [1\ 1\ 1]$, $L_2 = [1\ 1]$, $L_3 = 1$ are basis vectors of $\cN
 (K_i^f)^T$'s. The $\Pi_i$'s are the multinomial invariant distributions on 
 the fast components \cite{Gadgil:2005:SAF}, and are found to be 
 $$  
\Pi_1 = \left[\begin{array}{ccc} \dfrac{k_2^2}{(k_1+k_2)^2} & \dfrac{2k_1k_2}{(k_1+k_2)^2} 
   & \dfrac{k_1^2}{(k_1+k_2)^2}
   \end{array}\right]^T, \quad 
   \Pi_2 = \left[\begin{array}{cc} \dfrac{k_2}{k_1+k_2} & \dfrac{k_1}{k_1+k_2}
   \end{array}\right]^T, \quad  \Pi_3 = 1.
 $$                    
0's are zero row or    column  vectors with the appropriate dimensions.

As shown above, the dynamics on the slow time scale are governed by
 (\ref{reduced}), and in this example  the transition rate
 matrix for the slow  system is given by
$$
\tilde K = LK^s\Pi=\left[\begin{array}{ccc}  L_1K_1^s\Pi_1 &  L_1K_{1,2}^s\Pi_2 &  L_1K_{1,3}^s\Pi_3     \\
                                 L_2K_{2,1}^s\Pi_1 & L_2K_{2}^s\Pi_2 &  L_2K_{2,3}^s\Pi_3     \\
                                 L_3K_{3,1}^s\Pi_1  &   L_3K_{3,2}^s\Pi_2  &  L_3K_{3}^s\Pi_3 
                                \end{array}\right]=:\left[\begin{array}{ccc}  \tilde k_1^s &  \tilde k_{1,2}^s &  \tilde k_{1,3}^s     \\
                                \tilde k_{2,1}^s & \tilde k_2^s &  \tilde k_{2,3}^s     \\
                                 \tilde k_{3,1}^s  &   \tilde k_{3,2}^s  &  \tilde k_3^s 
                                \end{array}\right]
$$ 
 where the off-diagonal transition
rates $\tilde k_{i,j}^s$  are given by     
 $$ \tilde{k}_{1,2}^s = L_1K_{1,2}^s\Pi_2,\quad \tilde{k}_{1,3}^s
=\tilde{k}_{3,1}^s= 0, \quad \tilde{k}_{2,1}^s = L_2K_{2,1}^s\Pi_1,
\quad\tilde{k}_{2,3}^s = L_2K_{2,3}^s\Pi_3,\quad\tilde{k}_{3,2}^s =
L_3K_{3,2}^s\Pi_2.  
$$
The diagonal elements are the negatives of the corresponding column sums and
therefore $\tilde K$ is the generator of a Markov chain. We will show later  -- 
in Theorem~\ref{thm:mar}  -- that 
the transition rate matrix for the reduced system is alway a Markov chain
generator.  Recall as discussed in Subsection~\ref{subsec:inv}, we have 
$$
\cE^f = \left[\begin{array}{cc} -1 & 1\\
1 & -1\\
0 & 0
\end{array}
\right]\quad
A^f = \left[\begin{array}{ccc} 1 & 1 & 0\\
0 & 0 & 1
\end{array}
\right]\quad
\tilde n = A^fn = \left[\begin{array}{c} n_1+n_2\\
n_3
\end{array}
\right],
$$
Thus $\tilde n_1 = (2,0)$, $\tilde n_2 = (1,1)$, and $\tilde n_3 = (0,2)$.
The evolution
on the slow time scale is as shown schematically below.

 \centerline{
 \includegraphics[scale=0.7]{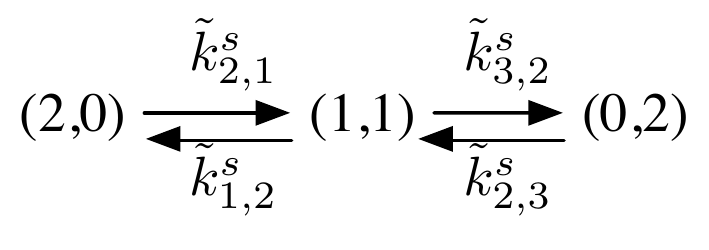}  
}
\end{Example}
The matrices for the general case in which  the total number of molecules is
$N_0$ are given in the Appendix.

\subsection{The structure of the generator of the slow dynamics}
\label{subsec:detail} 

We recall the standing assumption that there are no inputs or outputs that occur
on the fast time scale, and therefore on the fast time scale the reactions are
confined to fast simplexes as defined earlier. The slow dynamics moves the state
between fast simplexes $\Omega_f$ on a given simplex $\Omega$, and/or between
simplexes on the slow time scale. Because the transition matrix for the fast
dynamics has the block diagonal form given in (\ref{generalblock}), both $L$ and
$\Pi$ have the same block structure and it suffices to determine the null spaces
for a fixed block.

The graph $\cG_s^f$ can be decomposed into sources, internal strong components
and absorbing strong components. If all fast reactions on fast simplexes are
reversible, then all fast components are strongly connected.  If a state or
vertex is not connected to any other states, we call it an isolated state, and
it is then a (absorbing) strong component itself. Note that in the triangular reaction network as in Example~\ref{Triex}, each fast component is an absorbing strong component.

We let $S_i$ be the set of all states (nodes) in the graph of the $i^{th}$
  simplex and   denote the set of the states in the sources, internal strong components and
  absorbing strong components  by $S^{so}_i, S^{in}_i$ and $S^{ab}_i$,
  respectively. The set of all the states, $S$ can be represented by disjoint unions
$$
S= \dot \cup_i S_i = \dot \cup_i (S_i^{so} \dot \cup S_i^{in}\dot \cup
S_i^{ab}).
$$
We assume that $|S_i| = m_i$, $|S^{ab}_i|= m_{i}^{ab}$, $|S^{in}_i|= m_i^{in}$
and $|S_i^{so}|= m_i^{so}$. Thus $m_i=m_i^{ab} + m_i^{in}+m_i^{so}$. 

By following the general analysis of reaction networks developed in
\cite{Othmer:1979:GTA} we can write each block $K^f_i$ in (\ref{generalblock})
as an upper triangular block matrix
$$
K^f_i = \left[\begin{array}{ccc} K^{f,ab}_{i} &  K^{f,ab,in}_{i}  &
K^{f,ab,so}_i \\    &  K^{f,in}_{i} & K^{f,in,so}_{i} \\   &   & K^{f,so}_{i}
\end{array} \right],
$$ 
where each block   $K^{f,\al}_{i}$ is the  transition matrix between  states 
$S^\al_i$, where $\al = ab,in,so$ and  $K^{f,\al,\be}_{i}$ is the transition
matrix from the states in $S_i^{\be}$ into states in $S_i^{\al}$, where $\al,\be
= in, so$ or $ab$.  

If the $i^{th}$ simplex includes $p_i$  isolated absorbing states and $q_i$ absorbing
strong components with at least two states, then the matrix $K^{f,ab}_{i}$ can
be written as
$$K^{f,ab}_{i} = \left[\begin{array}{cccc}  {\bf 0}_{p_i} &    &  &  \\    &
K^{f,ab}_{i,1}  &  &  \\ &  & \ddots  &   \\ &   &  &   K^{f,ab}_{i,q_i}
\end{array}\right],$$ 
where  ${\bf 0}_{p_i}$ is a $p_i \times p_i$ zero matrix that reflects the transition rates
between the $p_i$ isolated absorbing states and $K^{f,ab}_{i,j}$ is an $m^{ab}_{i,j} \times
m^{ab}_{i,j}$ transition matrix describing the transitions between all states in
the $j^{th}$ absorbing strong components with at least two states. 

The following result defines an important spectral property of each $K^{f}_{i}$, which determines the structure of its left and right eigenvectors $L_i$ and $\Pi_i$ as shown in Subsection~\ref{subsec:inv}.

\begin{Theorem}\label{thm:semisimple}
Each matrix $K^{f,ab}_{i}$ has a semisimple zero eigenvalue and each $K^{f,in}_{i}$
and $K^{f,so}_{i}, \ i=1,\dots,l,$ is  nonsingular. 
\end{Theorem}
\begin{proof}The proof of this follows from the general result given in
\cite{Othmer:1979:GTA}. This theorem implies  that the zero eigenvalue of $K_{i}^f$ is semisimple.
\end{proof}

\subsection{The invariant distributions of the fast dynamics} \label{subsec:inv}

Next we consider the invariant distributions of the fast dynamics, which give
 the basis vectors in $\Pi$. These are vectors $\pi \ge 0$ such that $K^f\pi =0
 \textrm{ and } \sum_j \pi_{j} =1.$  Since every
 absorbing strong component with at least two states has a unique steady-state probability
 $\pi^{ab}_{i,j}$, $j = 1, \cdots,q_{i}$.   $\pi_{i,j}^{ab}$ is a basis for
 $\cN(K^{f,ab}_{i,j})$ since $K^{f,ab}_{i,j} \pi^{ab}_{i,j}=0$,  and dim
 $\cN(K^{f,ab}_{i,j})$=1. Notice that each absorbing strong component with
 only one state, which corresponds to each diagonal entry of the block ${\bf 0}_{p_i}$
 in $K^{f,ab}_{i}$, has steady-state probability $1$. Using these facts we define
 $$
\tilde \Pi_i =\left[ \begin{array}{cccc} I_{p_i} & & & \\ & \pi^{ab}_{i,1} &
 & \\ & & \ddots & \\ & & & \pi^{ab}_{i,q_i}
\end{array}\right],
$$ 
where $I_{p_i}$ is an $p_i\times p_i$ unit matrix and $\pi^{ab}_{i,j}$ is an
$m^{ab}_{i,j} \times 1$ vector. Since any states in sources and internal
strong components have zero probability at the steady-state, we define the 
matrices  

\begin{equation}
\Pi_i = \left[ \begin{array}{cccc}     & \tilde \Pi_i &  &  \\   {\bf 0}_{r_i}  &
{\bf 0}_{r_i}  & \cdots & {\bf 0}_{r_i} \\   {\bf 0}_{s_i}  &  {\bf 0}_{s_i}  &
\cdots & {\bf 0}_{s_i}    \end{array}\right] \qquad \qquad \Pi=
\left[ \begin{array}{ccc}   \Pi_1 &   &    \\   &  \ddots &   \\   &   & \Pi_l     
\end{array}\right],
\label{pimatrix}
\end{equation}
where ${\bf 0}_{r_i}$ and ${\bf 0}_{s_i}$ are null matrices representing the
steady-state probability of all states in internal strong components and
sources, respectively. Since  $K^f_i \Pi_i=0$ 
 it follows that $K^f \Pi=0,$ and therefore we have an orthogonal  basis for $\cN((K^f))$.

\medskip

To construct a basis of $\cN((K^f)^T)$, recall that  the multiplicity of an
eigenvalue is the same for a matrix and its adjoint, and therefore $K^{f,ab}_{i,j}$ has
exactly one zero eigenvalue and so dim$(\cN(K^{f,ab}_{i,j})^T)$ = 1 for each $j$.
The  $1 \times m^{ab}_{i,j}$  vector $\tilde L^{ab}_{i,j}
\equiv [1,\dots,1]$ is an eigenvector corresponding to the  zero eigenvalue of
$N(K^{f,ab}_{i,j})^T$, and therefore is a basis of
$\cN[(K^{f,ab}_{i,j})^T]$ for each $j$. Define
$$
\tilde L_i= \left[ \begin{array}{cccc} I_{p_i}  &   &   &   \\  &  \tilde
L^{ab}_{i,1} &   &   \\  &   & \ddots  & \\ &   &   &\tilde L^{ab}_{i,q_i}
\end{array}\right]
$$ 
and note that $\tilde L_i K^{f,ab}_{i}=0$. By defining $L_i =
\left[ \begin{array}{ccc} \tilde L_i  &  A_i   & B_i    \end{array}\right],$ 
where $A_i= -\tilde L_i K^{f,ab,in}_i (K^{f,in}_{i})^{-1}$ and $B_i=
-\big(\tilde L_iK^{f,ab,so}_{i}+A_i K^{f,in,so}_{i}\big) (K^{f,so}_{i})^{-1}$,
one can see that 
$$L_i K^f_i=  \left[ \begin{array}{ccc} \tilde L_i  &  A_i   & B_i
\end{array}\right] \left[\begin{array}{ccc} K^{f,ab}_{i} &  K^{f,ab,in}_{i}  &
K^{f,ab,so}_{i} \\  &  K^{f,in}_{i} & K^{f,in,so}_{i} \\  &   & K^{f,so}_{i}
\end{array} \right] = 0.$$ 
Thus,
$$
 L K^f = 0,
$$
where
\begin{equation}
L \equiv  \left[ \begin{array}{ccc}   L_1 &  &  \\  &  \ddots &    \\  &  &
L_l   \end{array}\right].
\label{lmatrix}
\end{equation}
One can see that all row vectors of $L_i$ consist of a basis of $\cN((K^f_i)^T)$
and it follows from the block diagonal structure of $L$ that all row vectors of
$L$ are basis vectors of $\cN((K^f)^T)$.  One can also show that the sets $\{L_i\}$
  and $\{\Pi_j\}$ are biorthogonal. 

\begin{Remark}
 When all  components are strongly connected, {\em e.g.,} all fast reactions are
 reversible,  the
 matrices $L$ and $\Pi$ reduce to 
$$
L=\left[ \begin{array}{ccc} \tilde L_1 &   &   \\   & \ddots &   \\   &   &
\tilde L_l \end{array}\right], \qquad \mbox{and} \qquad \Pi = \left[\begin{array}{ccc} \Pi_{1} & \dots  &   \\  & \ddots  &   \\   &  & \Pi_l \end{array}   \right],$$
where $\tilde{L}_i$ is a $ 1 \times m_i$ vector $[1, 1 \cdots 1]$ and $\Pi_i$ is
an $m_i \times 1$ steady-state vector of a strong component $\cL_i$.  Here $m_i$
is the number of states in $\cL_i$.
\end{Remark}

The following theorem shows that the slow evolution is always Markovian, even
when a fast simplex has multiple components, each of which may have sources
and/or internal components. 

\begin{Theorem}\label{thm:mar}
$\tilde K = LK^s \Pi$ is the generator of a  Markov chain.
\end{Theorem}
\begin{proof}
In order to prove this, we only need to show that the sum of each column of
$\tilde K$ is zero and the off-diagonal entries are nonnegative. Without loss of
generality and to avoid complicated subindices, we consider the first column of
$\tilde K$. The first column block of $\tilde K$, denoted as $\big(\tilde
K\big)_1$ is given by 
\[
\big(\tilde K\big)_1 = [L_1K_1^s\Pi_1\quad L_2K_{2,1}^s\Pi_1\quad \cdots\quad
L_nK_{n,1}^s\Pi_1]^T, 
\]
where $n$ is the number of fast components, block $L_iK_{i,1}^s\Pi_1\in {\bf
  R}_{\text{ab}_i\times \text{ab}_1}, i = 1,\ldots, n$, $\text{ab}_i$ is the
number of absorbing strong components of the $i$-th fast component. The sum of
all the rows of $\tilde K_1$ is 
\begin{align*}
&\sum_i\bigg(\sum_{k_i}(L_iK_{i,1}^s\Pi_1)^{\text{row}\ k_i}\bigg) =
\sum_i\sum_{k_i}(L_iK_{i,1}^s)^{\text{row}\ k_i}\Pi_1 
= \sum_i\bigg(\sum_{k_i}L_i^{\text{row}\ k_i}\bigg)K_{i,1}^s\Pi_1\\
&= \sum_i {\bf1}_{1\times D_i}K_{i,1}^s\Pi_1 = \bigg(\sum_i {\bf1}_{1\times
  D_i}K_{i,1}^s\bigg)\Pi_1 = {\bf 0}_{1\times D_1}\Pi_1 = {\bf 0}_{1\times
  \text{ab}_1}, 
\end{align*}
where $A^{\text{row}\ k_i}$ denotes the $k_i$-th row of matrix $A$ and  $D_i$ is the
number of nodes in the  $i$-th fast component. The  third equality follows from the fact
that the column sum of $L_i$ is always one, because each $L_i$ defines an
invariant probability and the sum must be one. The last equality follows from
the fact that the column sum of $K_{i,1}^s$ is zero. Since non-positive entries
appear only along the diagonal of $K_i^s$, the off-diagonal entries of $\tilde
K$ are nonnegative.
\end{proof}

To summarize, we now have the approximate equation for the probability distribution on the slow time
scale
\begin{equation}
\dfrac{d\tilde{p}}{dt} = LK^ s\Pi  \tilde{p} ,
\label{reducxx}
\end{equation}
 where the matrices $L$ and $\Pi$ are given above.
Clearly the utility of the reduction depends heavily on
whether the invariant distributions of the fast dynamics are easily
computed. This is possible for a network of first-order reactions, and for
closed systems the distribution is multinomial if the initial distribution is
multinomial, while for open first-order systems the distribution is a product
Poisson \cite{Gadgil:2005:SAF}. The general case of first-order reactions is
completely solved in \cite{Jahnke:2007:SCM}, where the authors show that the
solution can be represented as the convolution of multinomial and product
Poisson distributions with time-dependent parameters that evolve according to
the standard equations.

  Analytical results for bimolecular reactions are more limited. However it is
  known that the invariant distributions of a class of Markov chains that arise
  in queuing theory have a product form \cite{Boucherie:1991PFQ}, and it has
  been shown that this also holds for a zero-deficiency reaction network (which
  means that $\cR(\cE) \cap \cN({\nu}) = \phi$) under the assumption of ideal
  mass-action kinetics \cite{Anderson:2011:CTM,Melykuti:2014:EDS}. When this does not apply the
  stationary distributions frequently involve hypergeometric functions
  \cite{McQuarrie:1967:SAC}.

\subsection{ A detailed example of the reduction} 
\label{detex}
\begin{figure}[h!]
\begin{center}
{\includegraphics[scale = 0.3]{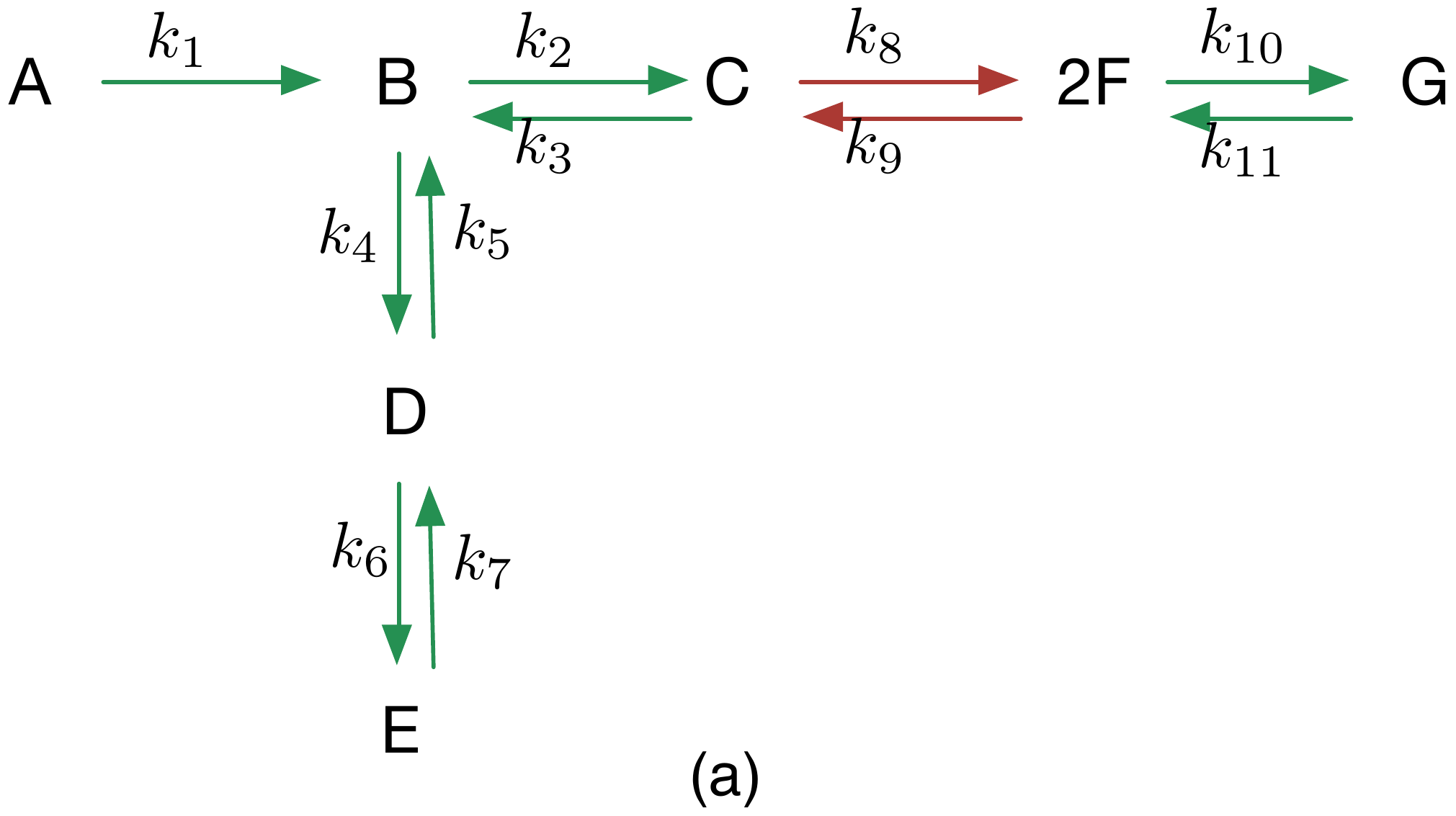}}\hspace{10pt}{\includegraphics[scale = 0.3]{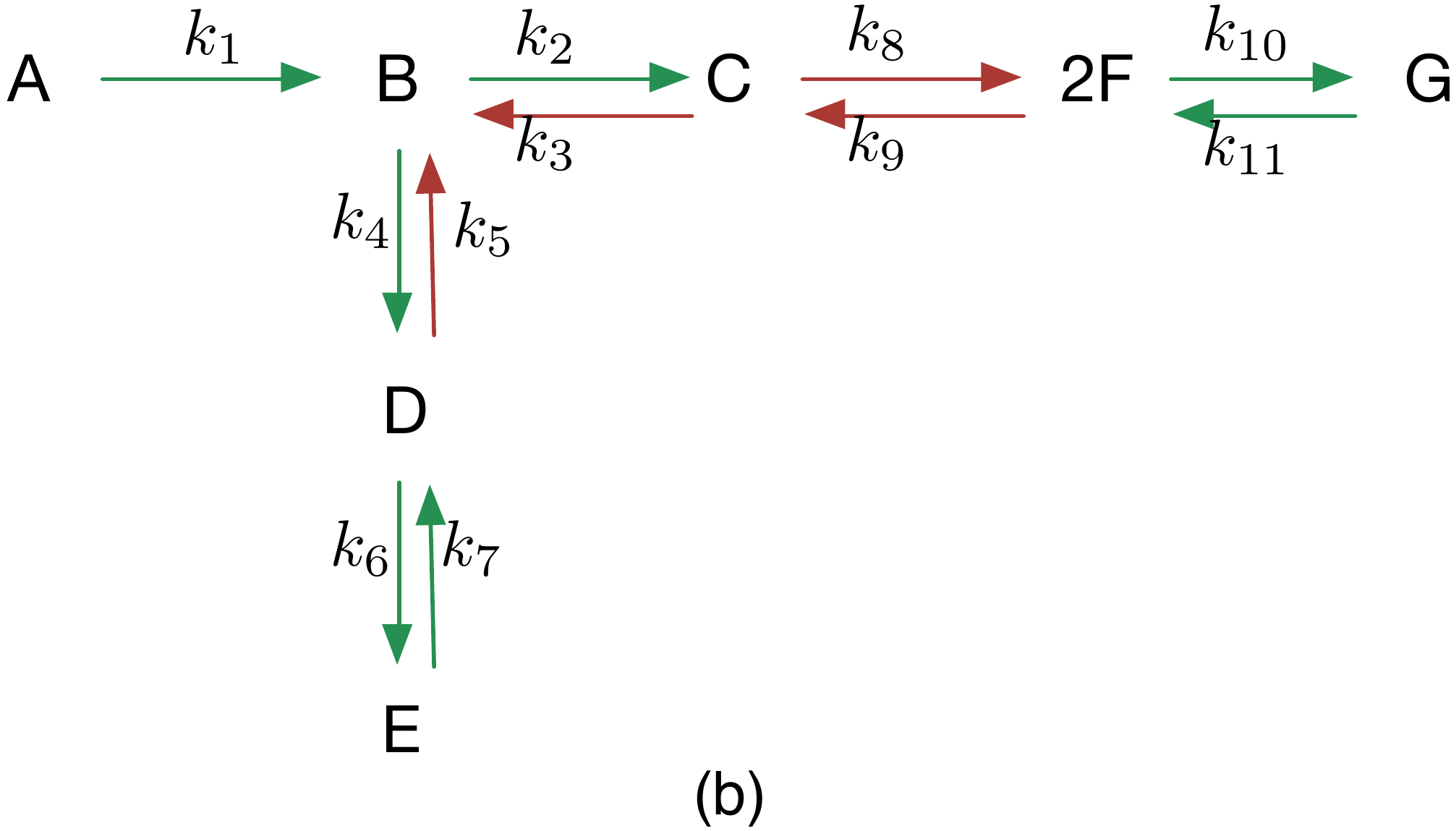}  }\\[10pt]
{\includegraphics[scale = 0.3]{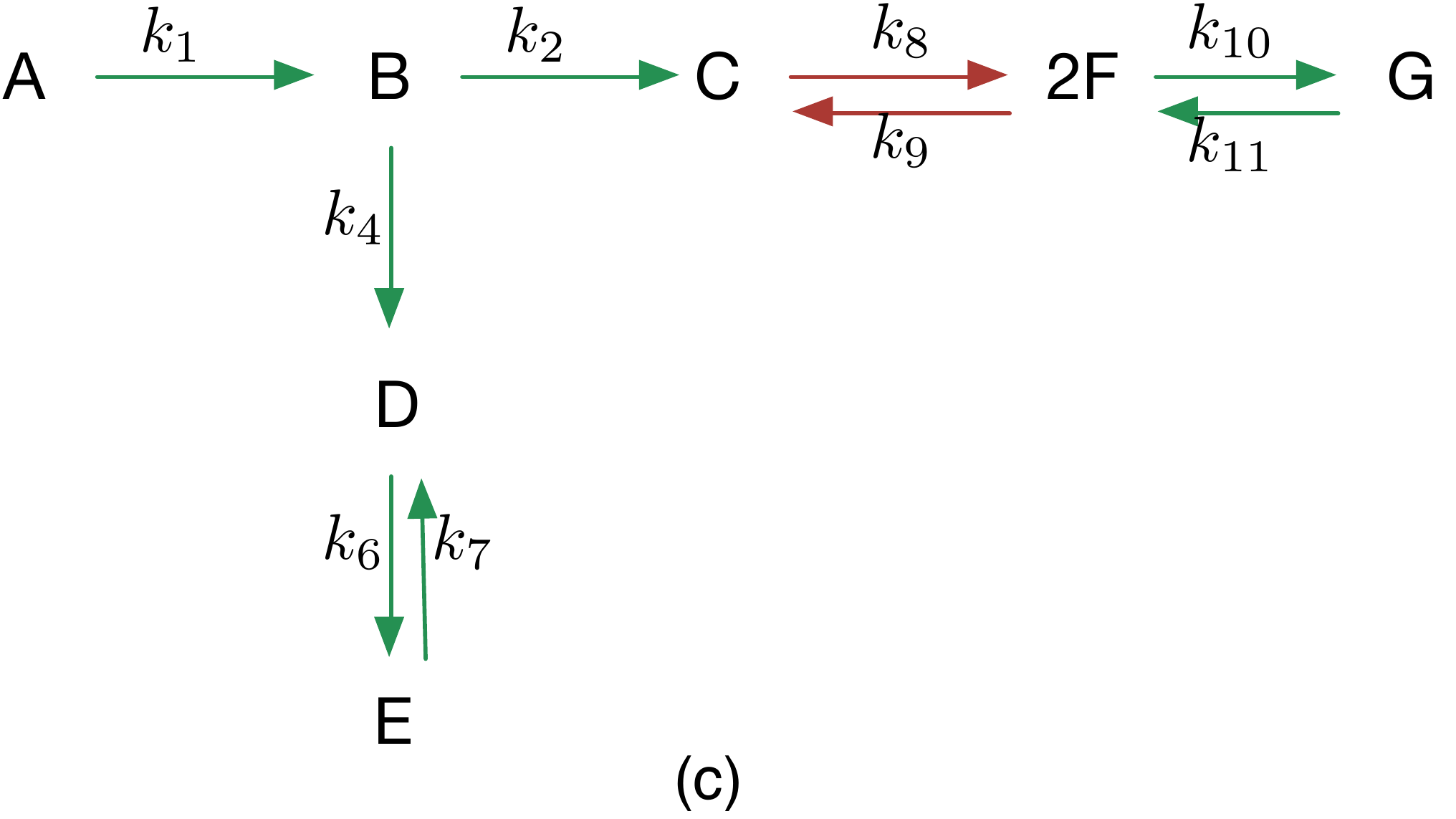}  }
\caption{\footnotesize  (a) Network 1. (b) Network 2. (c) Network 3. Red and green arrows denote slow and fast
 reactions, respectively. }
\label{fig:ab}
\end{center}
\end{figure}

\begin{figure}[h!]
\begin{center}
{\includegraphics[scale = 0.5]{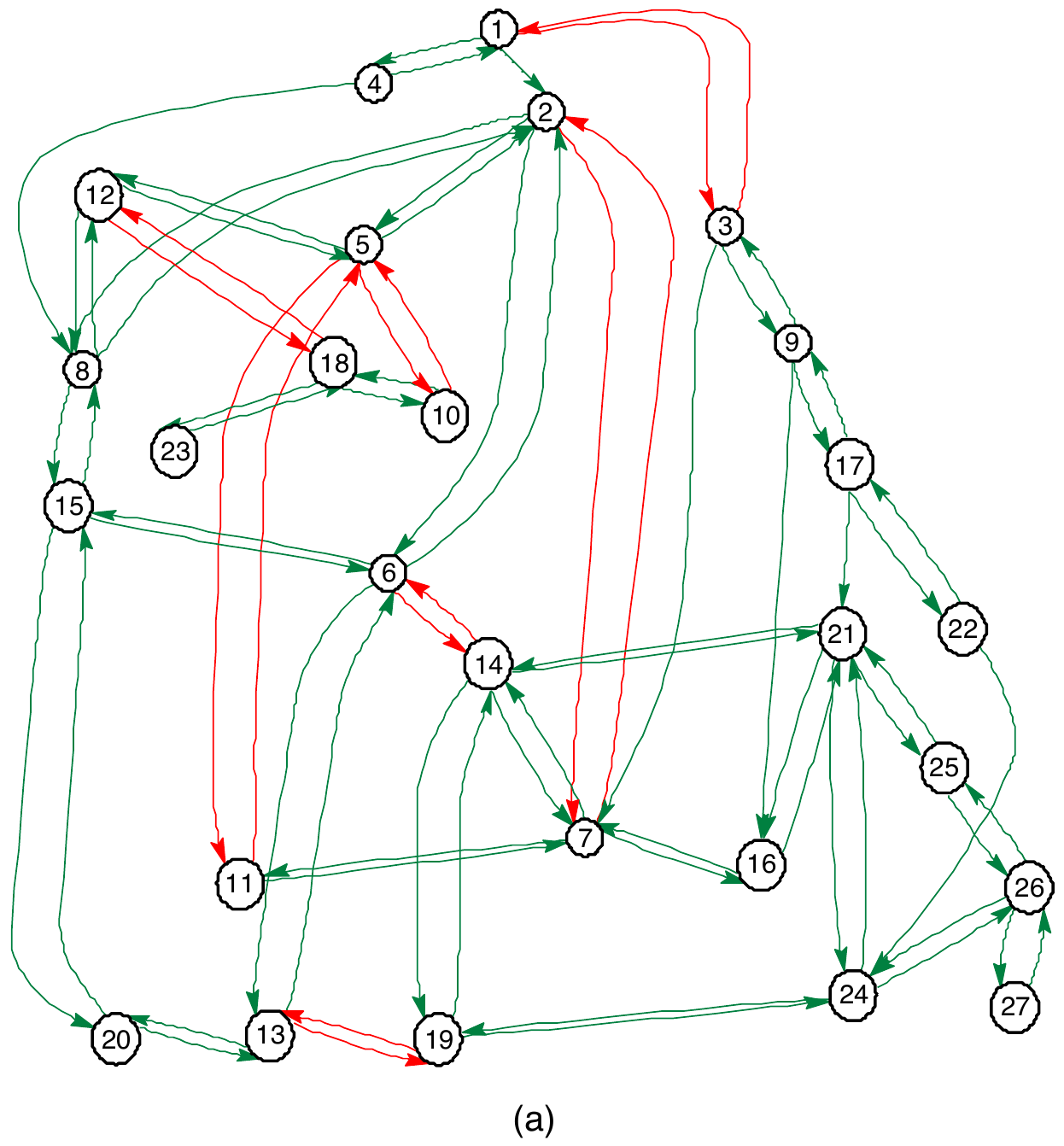}}\hspace*{10pt}{\includegraphics[scale = 0.5]{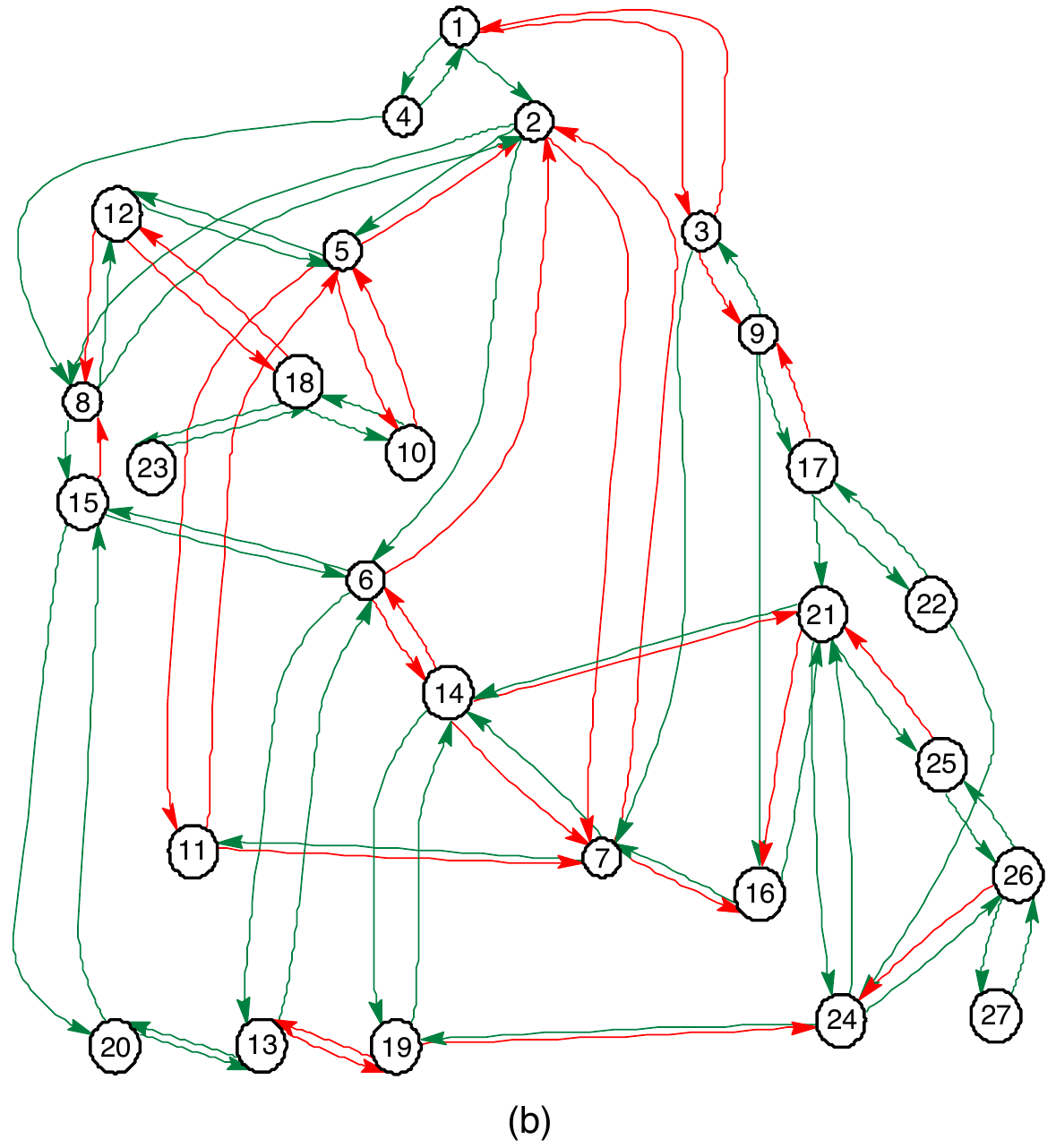}  }\\
{\includegraphics[scale = 0.5]{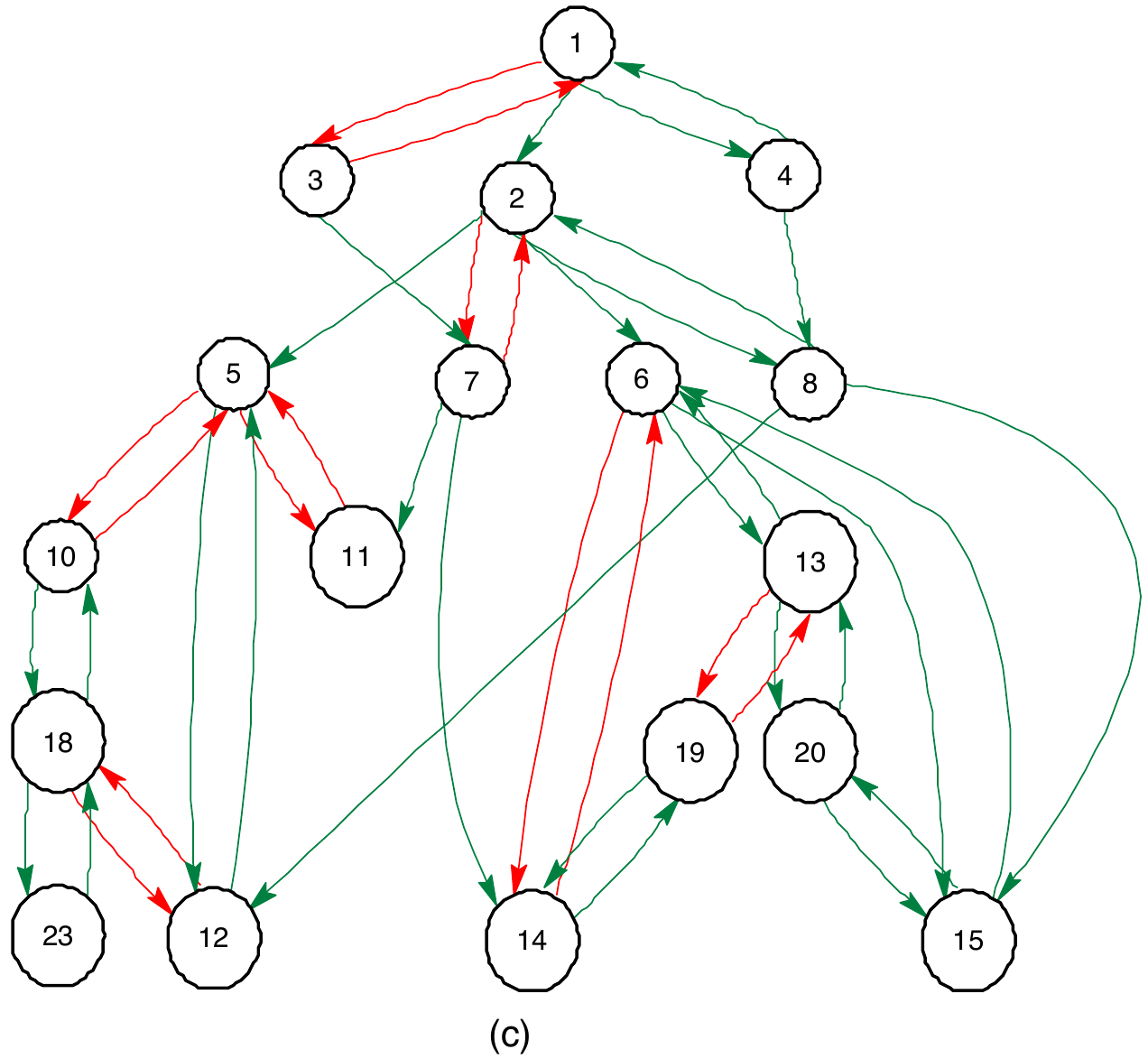}  }
\caption{\footnotesize  The graphs of the full state transition diagram
 with initial state \protect \circled{1}--$(1,0,0,0,0,2,0)$ for network 1-3, red arrows denote slow reactions. (a) Network 1. (b) Network 2. (c) Network 3. Red and green arrows denote slow and fast
 reactions, respectively. }
\label{fig:full}
\end{center}
\end{figure}

To illustrate the reduction method on general reaction networks, we consider the
three reaction networks shown in Figure~\ref{fig:ab}.  The resulting states
beginning from one molecule of A and two of F are shown in Table
\ref{state_tab}, where the states are in the order of molecular species $A, B, \cdots, G$. The state space of network 3 is a subset of that of networks 1
and 2 since reactions with rates $k_3$ and $k_5$ are missing, as seen in see
Table \ref{state_tab} (right).  The corresponding full state networks as shown in
Figure~\ref{fig:full} are generated using MATLAB.\\

\vspace*{-10pt}
\begin{table}[!htb]
\small\begin{center}
$\begin{array}{cc}
{ \begin{tabular}{ l | c | l | c | l | c}
\hline
\circled{1} & 1 0 0 0 0 2 0 &
\circled{2} & 0 1 0 0 0 2 0 &
\circled{3} & 1 0 1 0 0 0 0 \\ \hline
\circled{4} & 1 0 0 0 0 0 1 &
\circled{5} & 0 0 1 0 0 2 0 &
\circled{6} & 0 0 0 1 0 2 0 \\ \hline
\circled{7} & 0 1 1 0 0 0 0 &
\circled{8} & 0 1 0 0 0 0 1 &
\circled{9} & 1 1 0 0 0 0 0 \\ \hline
\circled{10} & 0 0 0 0 0 4 0 &
\circled{11} & 0 0 2 0 0 0 0 &
\circled{12} & 0 0 1 0 0 0 1 \\ \hline
\circled{13} & 0 0 0 0 1 2 0 &
\circled{14} & 0 0 1 1 0 0 0 &
\circled{15} & 0 0 0 1 0 0 1 \\ \hline
\circled{16} & 0 2 0 0 0 0 0 &
\circled{17} & 1 0 0 1 0 0 0 &
\circled{18} & 0 0 0 0 0 2 1 \\ \hline
\circled{19} & 0 0 1 0 1 0 0 &
\circled{20} & 0 0 0 0 1 0 1 &
\circled{21} & 0 1 0 1 0 0 0 \\ \hline
\circled{22} & 1 0 0 0 1 0 0 &
\circled{23} & 0 0 0 0 0 0 2 &
\circled{24} & 0 1 0 0 1 0 0 \\ \hline
\circled{25} & 0 0 0 2 0 0 0 &
\circled{26} & 0 0 0 1 1 0 0 &
\circled{27} & 0 0 0 0 2 0 0 \\ \hline
\end{tabular}} \quad
{ \begin{tabular}{ l | c | l | c | l | c}
\hline
\circled{1} & 1 0 0 0 0 2 0 &
\circled{2} & 0 1 0 0 0 2 0 &
\circled{3} & 1 0 1 0 0 0 0 \\ \hline
\circled{4} & 1 0 0 0 0 0 1 &
\circled{5} & 0 0 1 0 0 2 0 &
\circled{6} & 0 0 0 1 0 2 0 \\ \hline
\circled{7} & 0 1 1 0 0 0 0 &
\circled{8} & 0 1 0 0 0 0 1 &
 &  \\ \hline
\circled{10} & 0 0 0 0 0 4 0 &
\circled{11} & 0 0 2 0 0 0 0 &
 &  \\ \hline
\circled{13} & 0 0 0 0 1 2 0 &
\circled{14} & 0 0 1 1 0 0 0 &
\circled{15} & 0 0 0 1 0 0 1 \\ \hline
 &  &
 &  &
\circled{18} & 0 0 0 0 0 2 1 \\ \hline
\circled{19} & 0 0 1 0 1 0 0 &
\circled{20} & 0 0 0 0 1 0 1 &
\circled{21} & \\ \hline
\circled{22} & 1 0 0 0 1 0 0 &
\circled{23} & 0 0 0 0 0 0 2 &
\circled{24} &  \\ \hline
\end{tabular}}
\end{array}
$
\end{center}
\caption{Accessible states of networks. Left: Network 1 \& 2; Right: Network 3}
\label{state_tab}
\end{table}
\normalsize

\begin{figure}[t!]
 \begin{center}
 {\includegraphics[scale = 0.55]{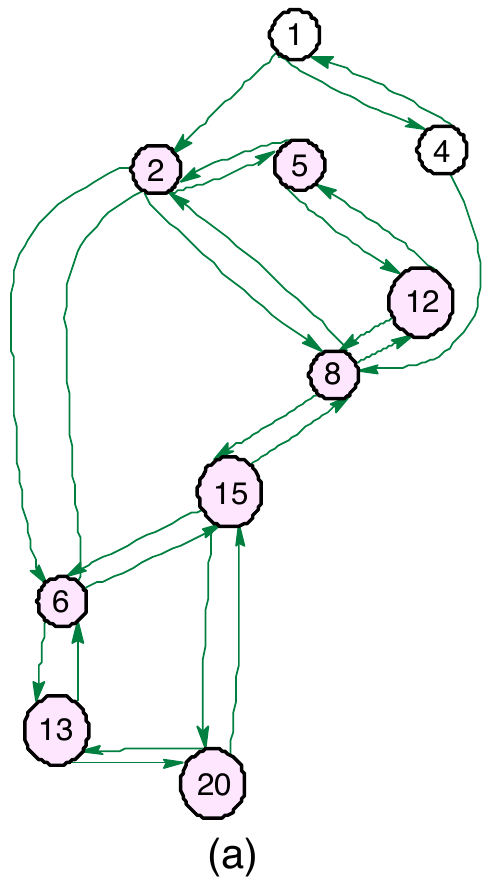}}\hspace{30pt} {\includegraphics[scale = 0.5]{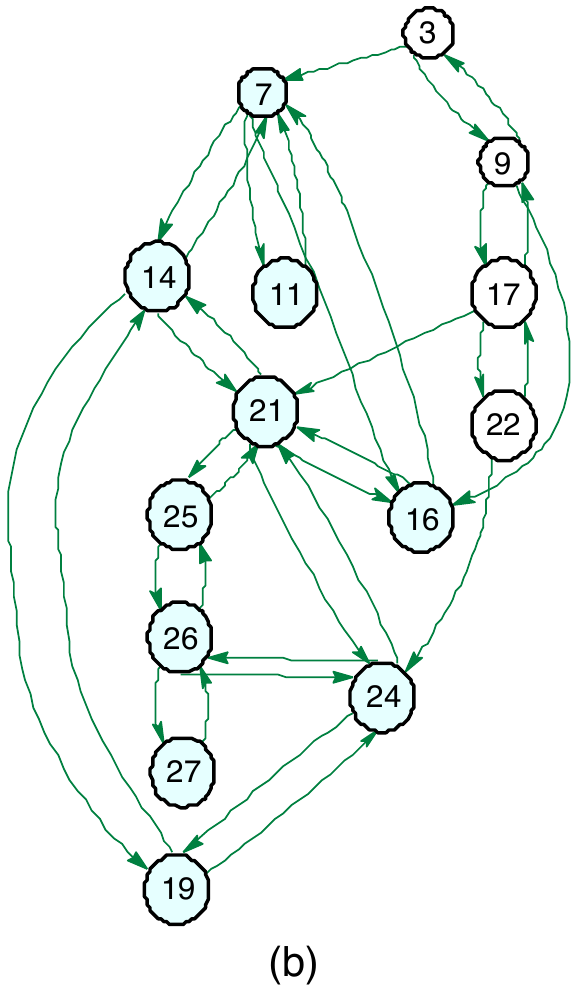}  }\hspace{30pt} {\includegraphics[scale = 0.7]{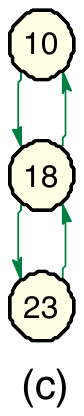}  }\\[10pt]
 {\includegraphics[scale = 0.6]{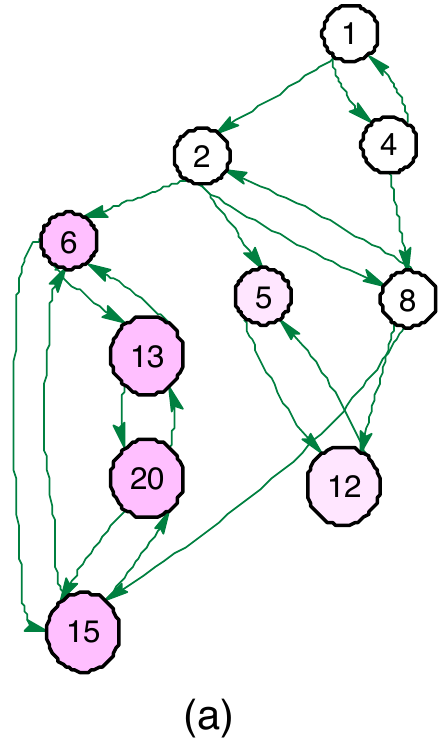}}\hspace{30pt} {\includegraphics[scale = 0.5]{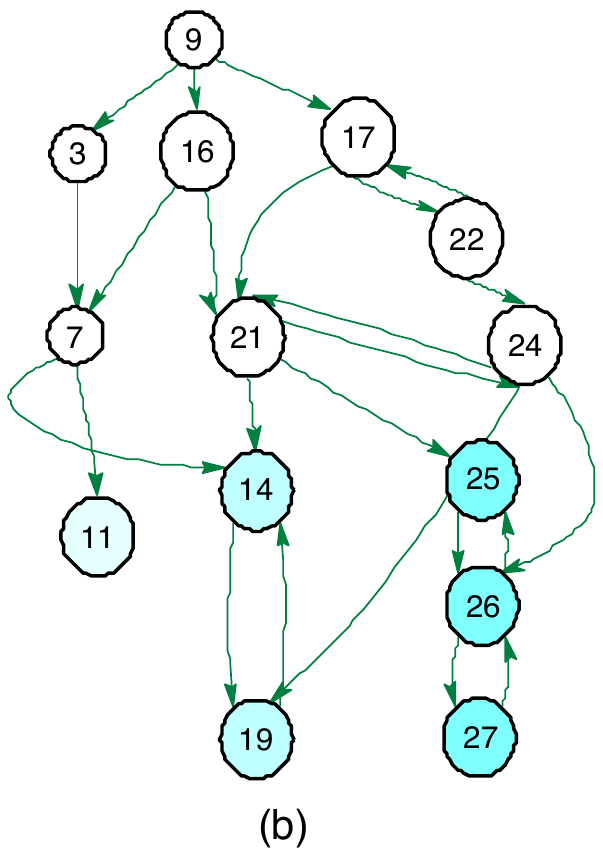}  }\hspace{30pt} {\includegraphics[scale = 0.7]{redu13.pdf}  }\\[10pt]
 {\includegraphics[scale = 0.5]{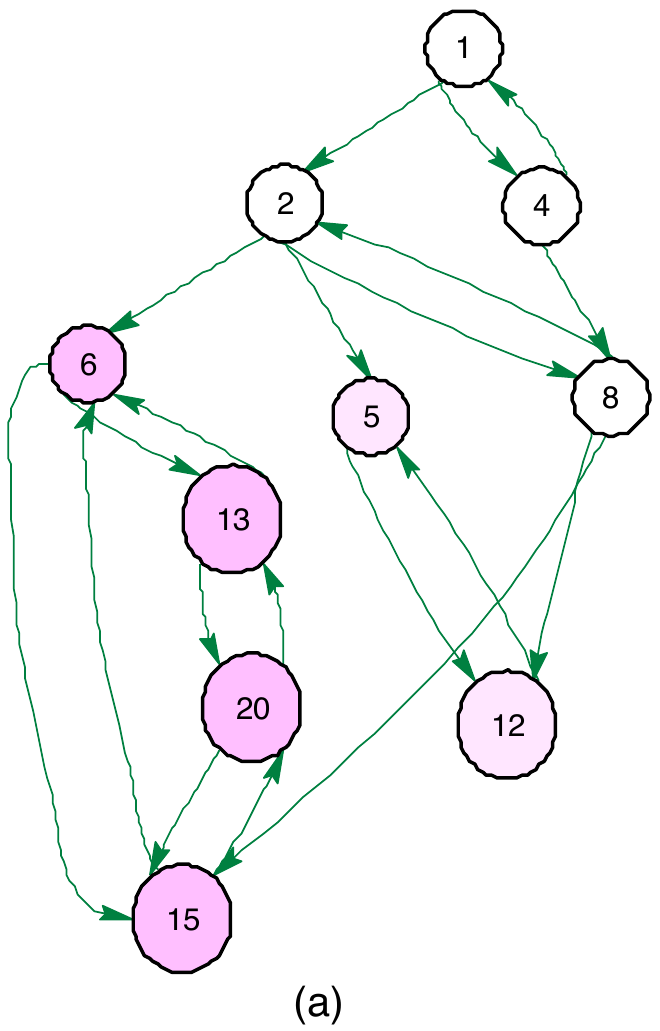}}\hspace{30pt} {\includegraphics[scale = 0.6]{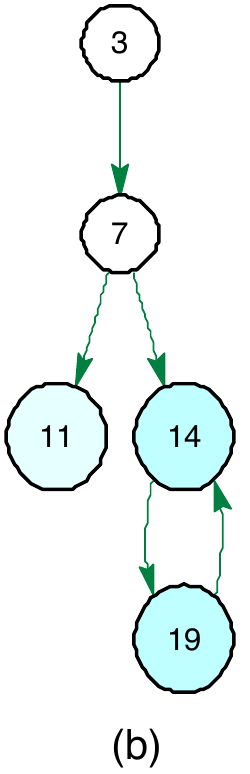}  }\hspace{30pt} {\includegraphics[scale = 0.7]{redu13.pdf}  }
 \end{center}
\caption{ The three fast components for network 1 (top), network 2 (center) and
  network 3 (bottom). The first has three fast components and three absorbing
  strong components,  the second has three  fast components
 and six absorbing strong components, and
the third has three fast components and five absorbing strong components.}
\label{reduc}
\end{figure}

By switching off the slow reactions in the graphs of the full state transition diagram shown in Figure~\ref{fig:full}
we obtain three fast components for each network as shown in Figure~\ref{reduc} (top), (center) and (bottom)
respectively. The fast components on the right for every network are strongly connected and are the same
since they are all generated by the reversible reactions $2F \leftrightarrow  G$. By highlighting the nodes appearing in
the same absorbing strong component, we see that in the state diagram of network 1, each fast component
has a unique absorbing strong component. Comparing the top and center graphs in Figure~\ref{reduc}, we see
the only difference occurs in the middle component. Network 2 has one more absorbing strong component
than network 3, and the remainder are the same. 

We perform the reduction method on the three networks and obtain the reduced
slow transition matrices for each network as follows. The detailed steps in the
reduction are shown in supplemental material available at \url{http://math.umn.edu/~othmer/Reduction.pdf}. 
One finds that the columns of the $3 \times 3$ matrix $\tilde{K}_1$, the generator of the slow dynamics
for network 1, are given by

\small
$$
 \bbordermatrix{ & (\tilde K_1)_1  &(\tilde K_1)_2]  \cr
\vspace{0.1cm}  & \frac{-k_9k_{11}k_5k_7(k_3+k_2)-k_8k_2k_5k_7(k_{10}+k_{11})-k_9k_{11}k_3k_4(k_6+k_7)}{(k_{10}+k_{11})(k_3k_4k_6+k_2k_5k_7+k_3k_4k_7+k_3k_5k_7)}& \frac{k_8k_{2}k_5k_7(k_3k_5k_7+2k_2k_5k_7+k_3k_4k_7+k_3k_4k_6)}{k_2k_5k_7(k_3k_4k_6+k_2k_5k_7+k_3k_4k_7+k_3k_5k_7)+\frac{1}{2}k_3^2(k_4k_6+k_4k_7+k_5k_7)^2} \cr
\vspace{0.1cm}  & \frac{k_9k_{11}}{k_{10}+k_{11}} & \frac{-k_8k_{2}k_5k_7(k_3k_5k_7+2k_2k_5k_7+k_3k_4k_7+k_3k_4k_6)}{k_2k_5k_7(k_3k_4k_6+k_2k_5k_7+k_3k_4k_7+k_3k_5k_7)+\frac{1}{2}k_3^2(k_4k_6+k_4k_7+k_5k_7)^2} \cr
\vspace{0.1cm} & \frac{k_8k_2k_5k_7}{k_3k_4k_6+k_2k_5k_7+k_3k_4k_7+k_3k_5k_7} & 0},
$$
\normalsize
$$
\bbordermatrix{ & (\tilde K_1)_3\cr
  \vspace{0.1cm}  & \frac{6k_9k_{11}(k_{10}+k_{11})}{3k_{10}^2+6k_{10}k_{11}+k_{11}^2} \cr
  \vspace{0.1cm}  & 0\cr
   & \frac{-6k_9k_{11}(k_{10}+k_{11})}{3k_{10}^2+6k_{10}k_{11}+k_{11}^2}}.
$$

The first three columns of $\tilde K_2$ are given by
$$
   \bbordermatrix{ & (\tilde K_2)_1 & (\tilde K_2)_2    & (\tilde K_2)_3 \cr
                  & -(k_3+k_8)-\frac{k_9k_{11}}{k_{10}+k_{11}} +
     \frac{k_3k_2}{k_2+k_4} & \frac{k_5k_2k_7}{(k_2+k_4)(k_6+k_7)}   & 2k_8 \cr 
                  & \frac{k_3k_4}{k_2+k_4} & \frac{-k_5k_2k_7}{(k_2+k_4)(k_6+k_7)}-\frac{k_9k_{11}}{k_{10}+k_{11}} &  0 \cr
                 & \frac{k_9k_{11}}{k_{10}+k_{11}} & 0 & -(k_3+2k_8)+\frac{k_3k_2}{k_2+k_4} \cr
                  & 0 & \frac{k_9k_{11}}{k_{10}+k_{11}} & \frac{k_3k_4}{k_2+k_4} \cr
                  & 0& 0& 0& \cr
                  & k_8 & 0 & 0  
                },
$$
 and the last three columns are given by
 $$
\bbordermatrix{  & (\tilde K_2)_4  & (\tilde K_2)_5 &(\tilde K_2)_6  \cr
                  & 0 & 0 & \frac{6k_9k_{11}(k_{10}+k_{11})}{3k_{10}^2+6k_{10}k_{11}+k_{11}^2}\cr
                 & k_8 & 0 & 0 \cr
                 & \frac{k_5k_2k_7}{(k_2+k_4)(k_6+k_7)} &  0 &  0\cr
                  & \frac{-k_5k_2k_7}{(k_2+k_4)(k_6+k_7)}-(k_3+k_8) + \frac{k_3k_2}{k_2+k_4} & \frac{2k_5k_2k_7}{(k_2+k_4)(k_6+k_7)} & 0\cr
                  & \frac{k_3k_4}{k_2+k_4} & \frac{-2k_5k_7}{k_6+k_7}+\frac{2k_5k_4k_7}{(k_2+k_4)(k_6+k_7)} & 0\cr
                  & 0 & 0 & \frac{-6k_9k_{11}(k_{10}+k_{11})}{3k_{10}^2+6k_{10}k_{11}+k_{11}^2} & 
                },
                $$
Finally, 
$$
  \tilde K_3
  = \bbordermatrix{ &   &   &  &   &   \cr
                  & -k_8-\frac{k_9k_{11}}{k_{10}+k_{11}} & 0   & 2k_8 & 0 & \frac{6k_9k_{11}(k_{10}+k_{11})}{3k_{10}^2+6k_{10}k_{11}+k_{11}^2}\cr
                  & 0 & -\frac{k_9k_{11}}{k_{10}+k_{11}} &  0 & k_8 & 0 \cr
                 & \frac{k_9k_{11}}{k_{10}+k_{11}} & 0 & -2k_8 &  0 &  0\cr
                  & 0 & \frac{k_9k_{11}}{k_{10}+k_{11}} & 0 & -k_8 & 0\cr
                  & k_8 & 0 & 0 & 0 & -\frac{6k_9k_{11}(k_{10}+k_{11})}{3k_{10}^2+6k_{10}k_{11}+k_{11}^2}}. $$

A number of general facts emerge from this example.

\begin{enumerate}
\item
If each fast component has a unique absorbing strong component, entries of the left eigenvector of the fast transition matrix $K^f$ are all one's , i.e. $L = [1, 1, \ldots, 1]$. In this case, there is no need to calculate $A$'s and $B$'s.

\item The column sum of $L_i$ is always one.

\item The dimension of the reduced slow transition matrix $\tilde K$ is the
  total number of absorbing strong components. Our reduction method does not
  require the uniqueness of absorbing strong component in each fast
  component. Without loss of generality, we consider the middle component of
  network 3 for example, see Figure~\ref{redu3eg}. It is obvious that the
  probability of \circled{7} choosing to transit to \circled{11} is
  $k_2/(k_2+k_4)$, and the probability of transiting to  \circled{14} is
  $k_4/(k_2+k_4)$. Conceptually we can suppose that  this 5-state component is
  comprised of  \circled{3} $\rightarrow$ \circled{7} $\rightarrow$ \circled{11}
  with probability $k_2/(k_2+k_4)$, and occurs as \circled{3} $\rightarrow$
  \circled{7} $\rightarrow$ \circled{14} $\ssang{}{}$ \circled{19} with
  probability $k_4/(k_2+k_4)$, and the   reduction method reflects this
  correctly. 

To see this, denote the probability vector corresponding to the 5-state
  component by
\[
P = [P_{11}, P_{14}, P_{19}, P_7, P_3]^T.
\]
Then the left eigenvector $L$ is found to be 
$$
L = \left[\begin{array}{ccccc}
1 & 0 & 0 & \frac{k_2}{k_2+k_4} & \frac{k_2}{k_2+k_4}\\
0 & 1 & 1 & \frac{k_4}{k_2+k_4} & \frac{k_4}{k_2+k_4}
\end{array}\right],
$$
and  the probability vector of the reduced network is given by
$$
\tilde P = LP = \left[\begin{array}{c}
P_{11} + \frac{k_2}{k_2+k_4}(P_7+P_3)\\
P_{14}+P_{19}+\frac{k_4}{k_2+k_4}(P_7+P_3)
\end{array}\right],
$$
the first row of $\tilde P$ corresponds to the probability of being trapped in 
the absorbing strong component \circled{11}, whereas the second row corresponds
to the probability of evolving to the other absorbing strong component
\circled{14}$\ssang{}{}$\circled{19}. The accumulated state associated with
\circled{3} $\rightarrow$ \circled{7} $\rightarrow$ \circled{11} is obtained by
switching off the reaction with rate $k_4$, which is $(2,0,0)$, and similarly
\circled{3} $\rightarrow$ \circled{7} $\rightarrow$ \circled{14} $\ssang{}{}$
\circled{19} gives rise to $(2,1,0)$.  

 \begin{figure}[!h]
 \begin{center}
 $\begin{array}{ccc}
 {\includegraphics[scale = 0.2]{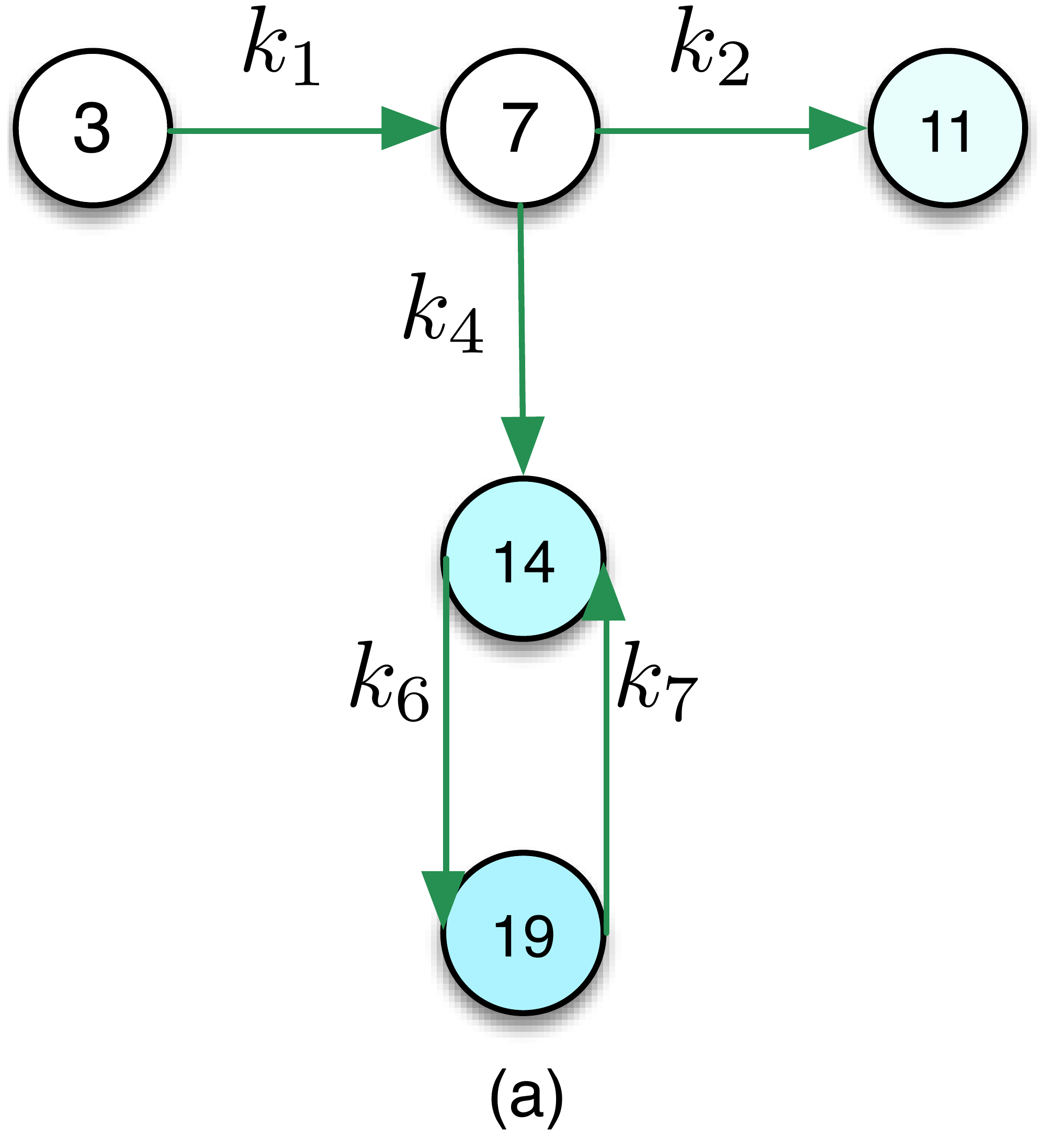}}& {\includegraphics[scale = 0.2]{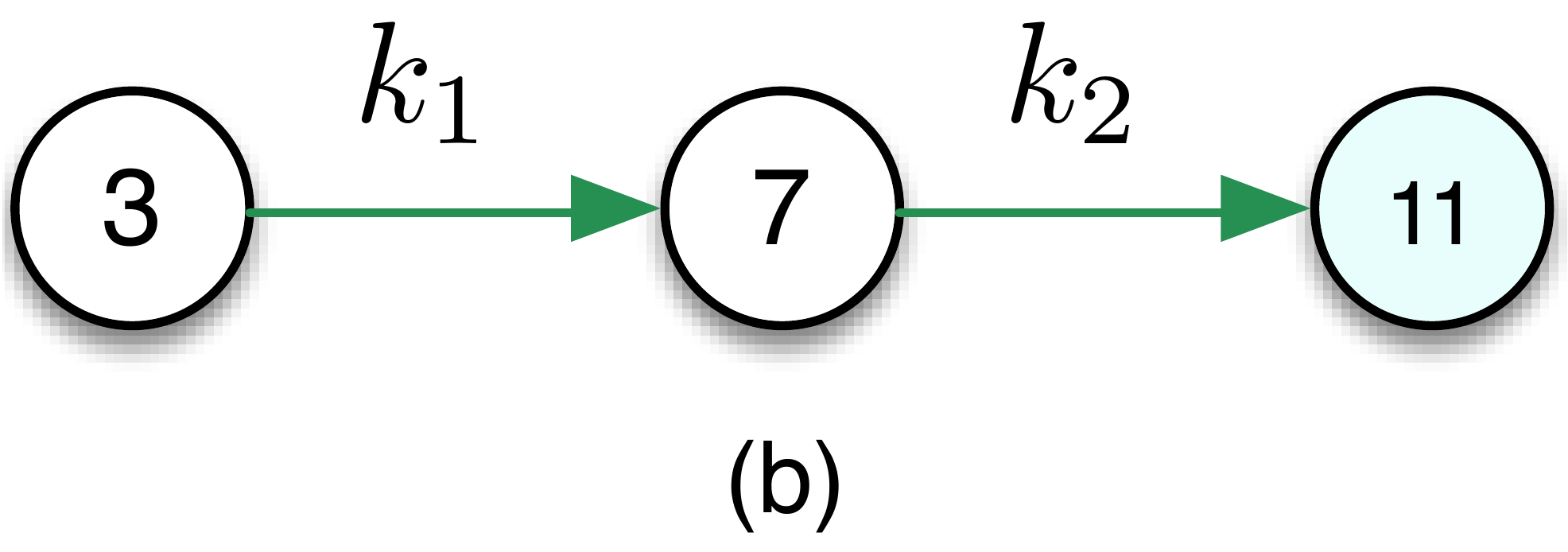}}& {\includegraphics[scale = 0.175]{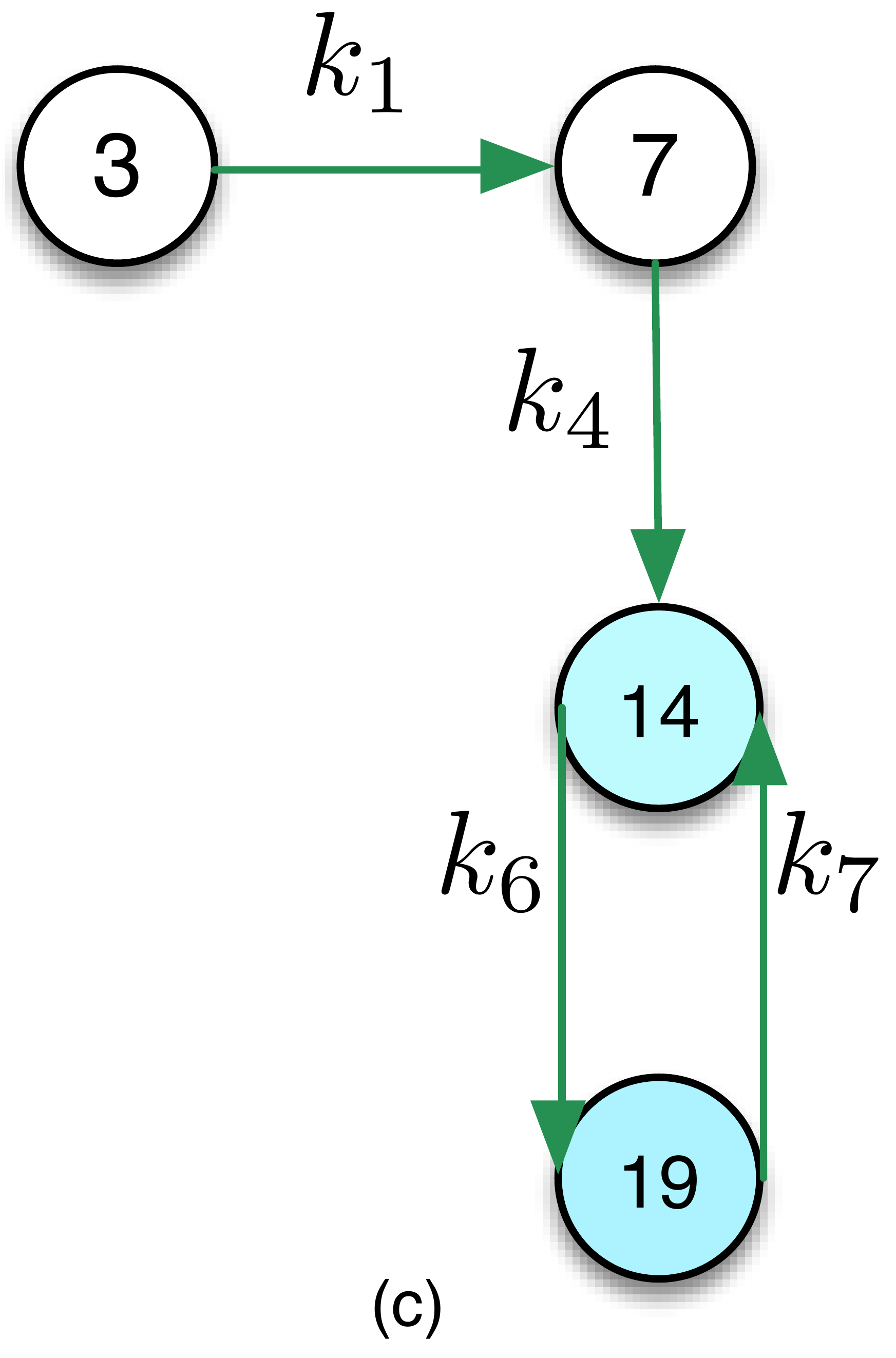}}
 \end{array}$
 \caption{\small Illustration of multiple absorbing strong components within one fast component: Network 2 middle graph}\label{redu3eg}
 \end{center}
 \end{figure}

\item
Comparing the reaction network 2 and 3 as shown in Figure~\ref{fig:ab}, we can
see network 2 has two more reverse slow reactions, i.e. reactions with rate
$k_3$ and $k_5$. Consider the reduced slow transition matrices $\tilde K_2$ and
$\tilde K_3$, if we let $k_3$ and $k_5$ be zeros in $\tilde K_2$, then $\tilde
K_2$ and $\tilde K_3$ are the same except $\tilde K_2$ has one extra row and
column of zeros, which correspond to the absorbing strong component
\circled{25}$\ssang{}{}$\circled{26}$\ssang{}{}$\circled{27} that only exists in
network 2. 

\end{enumerate}

\subsection{The mean and variance }

Often the first and second moments of the distribution of each species is
desired, and from the previous analysis one can find an approximate expression
for the  mean and variance of each species as follows.  Define a projection operator
$\cT_k$ as
$$\cT_k(n)= n_k,$$
where $n=(n_1,\dots,n_s)$.
Observing from the previous analysis that each state $S_{ij}$, the $j^{th}$
state in $\cL_i$,  has a probability distribution $\pi_{ij} \tilde p_i(t)$ at time
$t$,  one can find for each $k=1,\dots, s$ 
\begin{eqnarray}
E[n_k,t] &=& \sum_{i,j}\cT_k(n^{ij})\tilde p_i(t)\pi_{ij}, \label{remean} \\
Var[n_k,t] &=& \sum_{i,j} [\cT_k(n^{ij})]^2 \tilde p_i(t)\pi_{ij}-E[n_k,t]^2, 
\label{revar}
\end{eqnarray} 
where $n^{ij}$ denotes the vector $n$ of molecular numbers corresponding to the
state $S_{ij}$. Thus we can obtain formal solutions (\ref{remean}) and
(\ref{revar}) for mean and variance of each species on the slow dynamics. 

In linear models, it is possible to obtain explicit expressions for mean and
variance \cite{Gadgil:2005:SAF}. In that case one  can consider  transitions between
species as random walks of a molecule. This fact leads to a simple stochastic
algorithm for computation and  especially, if each fast component $C_i$ in the
reaction network is strongly connected, the quasi-steady-state distribution of
the fast component is multinomial \cite{Gadgil:2005:SAF}. Thus one can compute
the conditional expectation $E[\cR^s(n)|\tilde n]$ as follows; The
quasi-steady-state probability vector $\pi_i$ of a component $C_i$ can be
computed by using 
$$K^i_f \pi_i =0, \sum_{j=1}^{m_i} \pi_{ij}=1,$$
where $\pi_{ij}$ is a steady-state probability of $j^{th}$ species $\cM_{ij}$ in the $i^{th}$ component.  
If the reactant of a slow reaction  is $\cM_{ij}$, then one can obtain the
transition rate for the reaction  
$$E[\cR^s(n) \mid\tilde n] =  K^s \pi_{ij} \tilde n,$$
where $k_s$ is the transition rate constant of the reaction $R^s$.

Further insight into the structure of the slow dynamics is gotten from the fact
that  entry of $(LK^s \Pi)$ is a conditional expectation  when all fast
components are strongly connected, as proven next. 

\begin{Theorem}
Suppose that fast components are strongly connected.  Then the approximate
transition rate $\tilde \cR_{\ell}$ for the transition rate of all slow reactions from $j^{th}$ to $i^{th}$ fast component  is  
$$
\tilde \cR_{\ell} = \sum_{\ell} c_{\ell} E[h_{\ell}(n)\mid \tilde n]= \sum_{\ell} E[\cR^s_{\ell}(n) \mid \tilde n],
$$ 
where $c_{\ell}$ is the transition rate constant of the $\ell^{th}$ slow reaction which transits states from $j^{th}$ to $i^{th}$ fast component and
$h_{\ell}$ is the combinatorial number of reactants of the $\ell^{th}$ slow reaction. 
\label{scrate}
\end{Theorem}
\begin{proof} 
First note that by an occurrence of the $\ell^{th}$ slow reaction, a state in the
fast component $\cL_f^j$ is transformed into a state in the fast component
$\cL_f^i$. When a fast component is strong the variable $\tilde n =A_f n$ is uniquely defined in each
component, we suppose that $\tilde n= \tilde n^i$ for $\cL_f^i$ and $\tilde n =
\tilde n^j$ for $\cL_f^j$.  Let $n^{i,k}$ be the vector of molecular numbers that
denote the state $S_{ik}$, the $k^{th}$ state in the component $\cL_f^i$. 

Since $L K^s \Pi$ is a Markov chain generator, it follows from the reduced
equation (\ref{reduced}) that the approximate transition rate for a transitions
from the $j^{th}$ to the $i^{th}$ fast component is given by the $(i,j)^{th}$
entry of $ L K^s \Pi $, which is computed as follows.
 
\begin{eqnarray*}
(L K^s \Pi)_{ij} &=& (diag[L_1,\dots,L_{n_f}] K^s  diag[\Pi_i,\Pi_2,\dots, \Pi_{n_f}
                ])_{ij}\\
                &=& [ 1 \ 1 \ \dots \ 1] \left[ \begin{array}{ccc}
                (K^s_{ij})_{11} & \dots & (K^s_{ij})_{1m_j} \\  (K^s_{ij})_{21}
                & \dots & (K^s_{ij})_{2m_j}  \\ \vdots & \vdots & \vdots \\
                (K^s_{ij})_{m_i1} & \dots & (K^s_{ij})_{m_im_j}
                \end{array}\right] \left[\begin{array}{c} \pi_{j1} \\  \pi_{j2}
                \\ \vdots \\ \pi_{jm_j} \end{array}\right] \\ [10pt]
                &=& \sum_q \sum_p (K^s_{ij})_{pq} \pi_{jq}   =  \sum_q c_{\ell} h_{\ell}(n^{j,q}) P(n=n^{j,q} \mid \tilde n = \tilde n^j)\\
                &=&  c_{\ell} E[h_{\ell}(n) \mid \tilde n=\tilde n^j]\
                =  E[\cR^s_{\ell}(n) \mid \tilde n=\tilde n^j]
\end{eqnarray*}
where the $(i,j)^{th}$ block matrix of $K^s$, $K^s_{ij}$, is a matrix of
transition rates from the states of $\cL_f^j$ to those of $\cL_f^i$ and so
$(K^s_{ij})_{pq}$ is the rate of transitions from the $q^{th}$ state in $\cL_f^j$
into the $p^{th}$ state in $\cL_f^i$. In the fourth equality, we made use of the
fact  that  $\sum_p (K^s_{ij})_{pq}= \sum_{\ell}c_{\ell} h_{\ell}(n^{j,q})$  and  $\pi_{jq}$ is the
conditional probability $P(n=n^{j,q}\mid\tilde n = \tilde n^j)$.    
\end{proof}

\begin{Remark}
The result is true in general but the proof is somewhat more involved since the
structure of the $L_i$ and the $\Pi_j$ is more complicated. 
\end{Remark}

When  $\tilde n$ defined earlier is uniquely defined in each  
discrete simplex $\cL_f^i$, we have the relationship between the full and
reduced systems  shown in the following table. 

\begin{center}
\begin{tabular}{|c|c|c|}
\hline
& \textbf {Original system}    & \textbf{Reduced system} \\
\hline
\textrm{variable} &  $n$   &   $\tilde n = A^f n$  \\
\textrm{stoichiometry} &  $\nu\cE^s $     &  $\tilde{\nu\cE}^s =  A^f\nu\cE^s$ \\
\textrm{transition rate }  & $\cR^s(n)$   & $\tilde \cR^s(n)= E[\sum \cR^s | \tilde n]$ \\
\hline
\end{tabular}
\end{center}
Moreover, under the assumptions in Theorem \ref{scrate} we can obtain an
approximate chemical master equation on the slow time scale when the fast
components are strongly connected.  If the original master equation is given by
\begin{eqnarray*}
\dfrac{d}{dt}P(n,t) &= & \sum_{\ell}\dfrac{1}{\ep} \big[\cR^f_{\ell}(n-\nu\cE^f_{(\ell)})\cdot
P(n-\nu\cE^f_{(\ell)},t)-\cR^f_{\ell}(n)\cdot
P(n,t)\big]\\ 
&& + \sum_{k} \big[\cR^s_{k}(n-\nu\cE^s_{(k)})\cdot
P(n-\nu\cE^s_{(k)},t)- \cR^s_{k}(n)\cdot
P(n,t)\big],
\end{eqnarray*}
then the master equation in the reduced system is approximated (with error less
than $O(\ep)$) 
\begin{equation}
\dfrac{d\tilde p(\tilde n,t)}{dt} = \sum_k \tilde R^s_k(\tilde n-\tilde{\nu\cE^s}_{(k)}) \tilde p(\tilde n-\tilde{\nu\cE^s}_{(k)},t) -  \tilde R^s_k(\tilde n) \tilde p(\tilde n,t), \label{redmaster}
\end{equation}
where $\tilde n= A^f n$, $\tilde {\nu\cE^s}_{(k)} = A^f \nu\cE^s_{(k)}$ and
$\tilde R^s_k(\tilde n)=  E[\sum\cR^s_k(n)|\tilde n]$. In a later section we obtain
a modified stochastic simulation algorithm from (\ref{redmaster}). If there are
sources or internal strong components the master equation for the slow system
contains additional terms.

\section{A stochastic simulation algorithm for  the slow  dynamics}
 
\subsection{ A stochastic simulation algorithm based on the QSS approximation}

Preparation:

\begin{itemize}

\item Identify fast and slow reactions.

\item For a given initial state $n(0)$, switch off the slow reactions and generate
  the transition matrix $K^f$ of the fast component that $n(0)$ lies in by implementing Algorithm~\ref{alg}. 

\item Identify the fast simplex that $n(0)$ lies in to determine $A^{f,ab}$.

\item Define $\tilde n\equiv A^{f,ab} n$, $\widetilde{\nu\cE^s} \equiv A^{f,ab} \nu\cE^s$.

\item Denote the fast simplex that $n(0)$ lies in to be  $\tilde n(0) = A^{f,ab}n(0)\equiv\tilde n_j$, i.e., $jth$ fast simplex and its target fast simplexes through slow reactions, denoted as $\tilde n_i$(could be more than one).

\item Identify the fast transition matrices within the $ith$ and $jth$ fast simplexes, respectively, denoted as $K^{f,ab_i},K^{f,ab_j}$. Identify slow transition matrices from $jth$ to $ith$ fast simplexes, denoted as $K_{i,j}^s$.

\item Compute $L_i$ and $\Pi_j$ such that $L_iK^{f,ab_i} = 0$ and $K^{f,ab_j}\Pi_j = 0$.

\end{itemize}

Step 1. {\bf Scheme for a fast simplex} : Computing slow reaction rates from the
QSS state of fast a simplex 
\begin{enumerate} 
\item Compute the slow transition rate $\tilde \cR^s(\tilde n_j) \equiv \tilde
  K_{i,j}^s = L_iK_{i,j}^s\Pi_j$ from the $jth$ fast simplex to all its target simplexes.

\item  Compute $\tilde R^s_{tot}= \sum_i \tilde\cR^s_i(\tilde n_j)$.  
\end{enumerate}

\medskip

Step 2. {\bf Scheme for simulation of slow reactions} : Simulating the slow reactions with the reaction rates obtained in Step 1
\begin{enumerate}
\item Generate two random numbers $r_1$ and $r_2$ from the uniform distribution
  on $(0,1)$.

\item Set $\tau = -\dfrac{\log(r_1)}{\tilde R^s_{tot}} $ and choose $k$ such that \\
$ \sum_{i=1}^{k-1} \tilde \cR^s_i(\tilde n_j) < r_2 \tilde R^s_{tot} \le \sum_{i=1}^{k} \tilde \cR^s_i(\tilde n_j)$.
\end{enumerate}

\medskip

Step 3. {\bf Update}

Let $t \leftarrow t+\tau$ and  $\tilde n \leftarrow \tilde n + \widetilde {\nu\cE^s_k}$. Go to Step 1.

\begin{Remark}

\item[1.] The main idea of the modified Gillespie algorithm is to perform the
  traditional algorithm on a reduced rather than the original reaction network.
  For instance, in Example \ref{Triex},
$$
A^f = \left[\begin{array}{ccc} 1 & 1 & 0\\
0 & 0 & 1
\end{array}
\right]\quad
\tilde n(0) = A^fn(0) = (2,0)\quad
\widetilde {\nu\cE^s} = A^f\nu\cE^s = \left[\begin{array}{cc} -1 & 1\\
1 & -1
\end{array}
\right],
$$
and $\widetilde {\nu\cE^s}$ provides the stoichiometric matrix for a  reduced
reaction network.  The reduced reaction network can be written as
$$
\mathcal M_1 \rightleftharpoons \mathcal M_2,
$$
where $\mathcal M_1, \mathcal M_2$ are what we call ``pseudo species" in the reduced network, and the new initial state $\tilde n(0) = (2,0)$ corresponds to only two molecules of species $\mathcal M_1$. Each pseudo species represents one fast (absorbing) strong component, in this case
$$
\mathcal M_1: A\rightleftharpoons B\quad \mathcal M_2: C.
$$
The slow reaction rates in the reduced network are computed using  $LK^s\Pi$ as before except $K^f$ and $K^s$ are for the original reaction diagram rather than state diagram, i.e.
$$
K_1^f = \bbordermatrix{ & A & B  \cr
                A & -k_1 & k_2 \cr
                B & k_1& -k_2}
                \quad
K_2^f = \bbordermatrix{& C \cr
C & 0}\quad
K_{2,1}^s = \bbordermatrix{& A & B \cr
C & k_6 & k_3}\quad
K_{1,2}^s = \bbordermatrix{& C \cr
A & k_5 \cr
B & k_4}
$$
After finding the forward reaction rate from $\mathcal M_1$ to $\mathcal M_2$ to be $(k_1k_3+k_2k_6)/(k_4+k_5)$ and the backward rate $k_4+k_5$, one can initialize the Gillespie algorithm as usual.

\end{Remark}
 
 Example~\ref{Triex} is the special case in which each fast component is
 strongly connected, i.e. only fast absorbing strong components or sinks
 exist. In general, one should identify sink, internal, and source strong
 components in the reaction graph. Once they are identified, the graph becomes a
 tree. Then one can use for example depth-first search to identify the fast
 discrete simplexes as defined in \ref{fastsimplex} in order to find $A^f$.

The algorithm entails computation of the matrix of the slow dynamics, and whether
one does this once for all initially, which may be advantageous when doing
multiple realizations, or on the fly when only a few realizations are desired,
is a matter of choice. In single realizations the entire state space generated
by Algorithm 1 may not be explored and 'on-the-fly' computation of the slow
transition rates may be advantageous.

\section{Applications} 

In this section we analyze three examples: a bacterial moter model, an enzyme-inhibitor model and a model from PFK system. These examples illustrate the reduction when
there are relatively few states in the system. In these cases one can find the
approximate probability distribution of the slow variables by solving the
reduced matrix equation ${d\tilde p}/{dt} = LK^s \Pi \tilde p$ directly. Of
course the main step is to find the matrix in the reduced equation
$$
\dfrac{d\tilde p}{dt} = LK^s \Pi \tilde p.
$$

\begin{example}[A bacterial moter model]

This  example arises as a model for control of the rotational bias in the
flagellar motor of {\em E. coli} \cite{Othmer:2005:ACR}. At the base of the
motor are sites at which the protein $CheY_{p}$ can bind, and the occupancy of
the sites biases the probability of switching the direction of rotation of the
motor. Here we consider the following  scheme for binding $CheY_p$ (represented by $Y$) to these
sites. 

\small
\begin{displaymath}
 \begin{array}{cccccccccccccc}
 \vspace{.01in}
 & nk_{1}Y& &  (n-1)k_{1}Y& & \qquad &  2k_{1}Y& &  k_{1}Y \\[-2pt]
 CW_0 & \harr & CW_1 & \harr & CW_2 & \cdots & \harr & CW_{n-1} & \harr & CW_{n} \\
 [-5pt]
 &  k_{-1}& & 2k_{-1}& \qquad & & (n-1)k_{-1}&
 & nk_{-1} \\
 \alpha_{0}\duarr\beta_{0} &  & \alpha_{1} \duarr \beta_{1} &
  & \alpha_{2} \duarr \beta_{2} & \cdots &   & \alpha_{n-1} \duarr
 \beta_{n-1} &  & \alpha_{n} \duarr \beta_{n} \\
 & nk_{3}Y& & (n-1)k_{3}Y& & \qquad & 2k_{3}Y& & k_{3}Y \\[-2pt]
 \vspace{.01cm}
 CCW_{0} & \harr & CCW_{1} & \harr & CCW_{2} & \cdots & \harr & CCW_{n-1} & \harr & CCW_{n} \\
 \vspace{.015cm}
 & k_{-3}& & 2k_{-3}& \qquad & & (n-1)k_{-3}& & nk_{-3} \\
 \end{array}
 \end{displaymath}

 %

 %
 %

\normalsize

Here $CW_k$ and $CCW_k$ represent clockwise and counterclockwise flagellar rotation, respectively. We assume that the horizontal transitions are fast, while the vertical
transitions are slow, which  leads to  two strongly connected fast
components comprising the horizontal steps.

Thus the two $(n+1)\times (n+1)$ fast reaction rate matrices are given by
\small
$$
K^f_1= \left[ \begin{array}{ccccc} -nk_1Y & k_{-1} &  0 & \dots &  \\ nk_1Y &
-(k_{-1} +(n-1)k_{1}Y) & \ddots &   &  \\  & (n-1)k_{1}Y & \ddots & \ddots &  \\
&   & \ddots & \ddots & nk_{-1}\\  & & &k_1Y & -nk_{-1} \end{array} \right]
\quad
K^f_2=  \left[ \begin{array}{ccccc} -nk_3Y & k_{-3} &  0 & \dots &  \\ nk_3Y &
-(k_{-3} +(n-1)k_{3}Y) & \ddots &   &  \\  & (n-1)k_{3}Y & \ddots & \ddots &  \\
&   & \ddots & \ddots & nk_{-3}\\  & & &k_3Y & -nk_{-3} \end{array} \right].
$$ 
\normalsize
Moreover, the slow reaction rate matrix  is  
$$K^s = \left[ \begin{array}{cc} -A & B \\ A  & -B \end{array}\right],$$
where 
$$
A =diag(\alpha_0,\dots,\alpha_n) \textrm{ and } B=diag(\beta_0,\dots,\beta_n).
$$
Let $\Pi=[ \Pi_1 \ | \ \Pi_2]$ and let $\Pi_i,i=1,2$ be stationary distributions
of $K^f_i,i=1,2$ respectively , {\it i.e}, $\Pi_i$ is the right eigenvector of
$K^f_i$ corresponding to the zero eigenvalue with $\sum_j \Pi_{ij} =1$, where 
$\Pi_{ij}$ denotes the $j^{th}$ entry of the vector $\Pi_i$. Then the
reduced equation  can be written
$$
\frac{d\tilde p}{dt}= LK^s\Pi \tilde p,
$$
where 
\begin{eqnarray*}
LK^s\Pi &=& \left[ \begin{array}{cc} -\sum_{i=1}^{n+1} \alpha_{i-1} \Pi_{1i} &
\sum_{i=1}^{n+1} \beta_{i-1} \Pi_{2i}\\  \sum_{i=1}^{n+1} \alpha_{i-1} \Pi_{1i}
& -\sum_{i=1}^{n+1} \beta_{i-1} \Pi_{2i} \end{array}\right] 
\end{eqnarray*}
Thus we reduce the system into two-state model,
$$CCW \ \ \ssang{k^+}{k^-} \ \ CW,$$
where $k^+= \sum_{i=1}^{n+1} \alpha_{i-1} \Pi_{1i}$ and $k^-=\sum_{i=1}^{n+1}
\beta_{i-1} \Pi_{2i}$. Here we note that $k^{\pm}$ are functions of
$k_i,n,\alpha,\beta$ and $Y$. These parameters are reported in the experimental
literature. 
 
Figure~\ref{flagellum} illustrate the evolution of probabilities of $CW$ and $CCW$ obtained from
the full and the reduced equations. 

\begin{figure}[h!] 
\begin{center}
\begin{tabular}{cc}
{\includegraphics[width=2.5in]{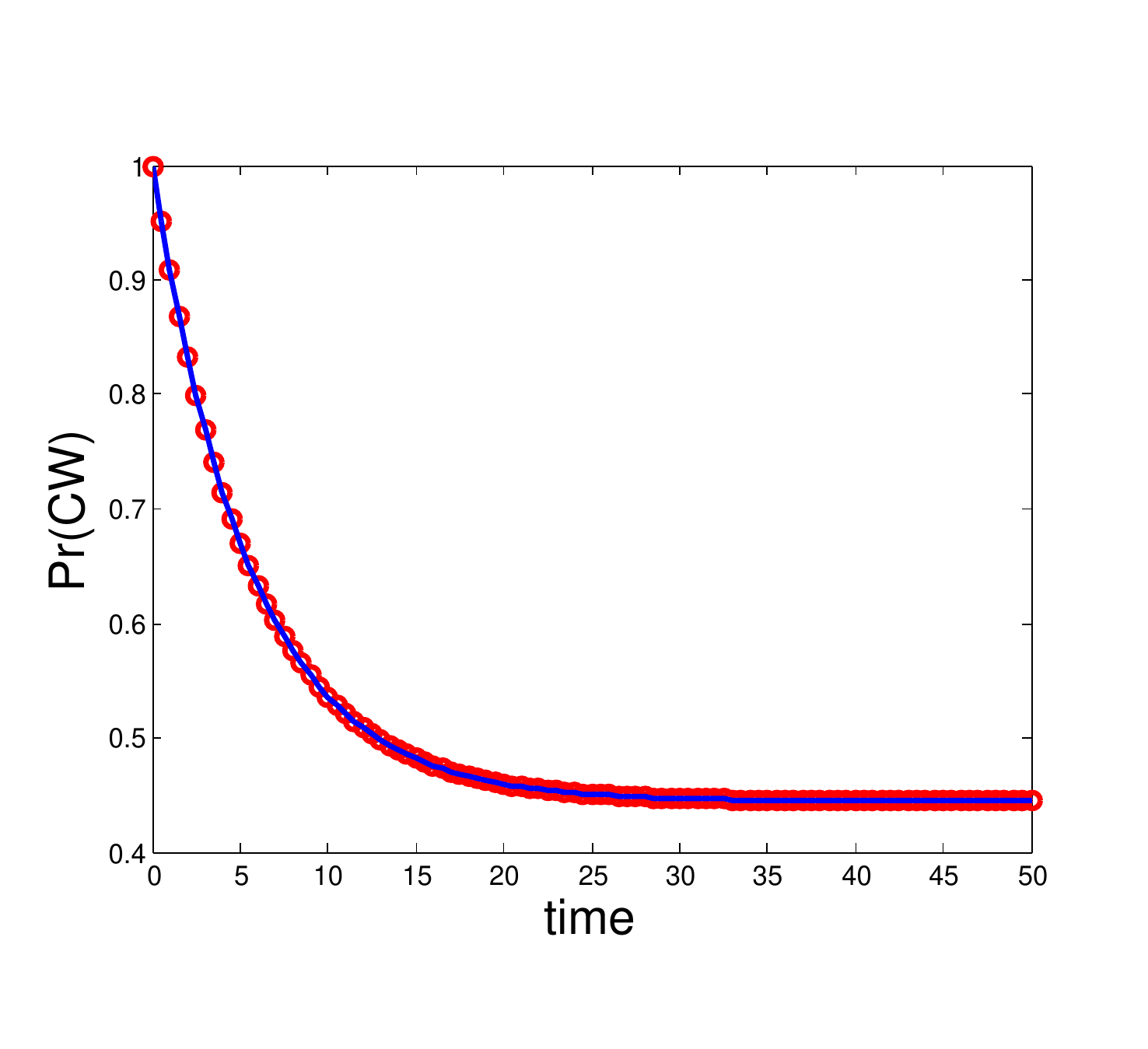}}  & {\includegraphics[width=2.5in]{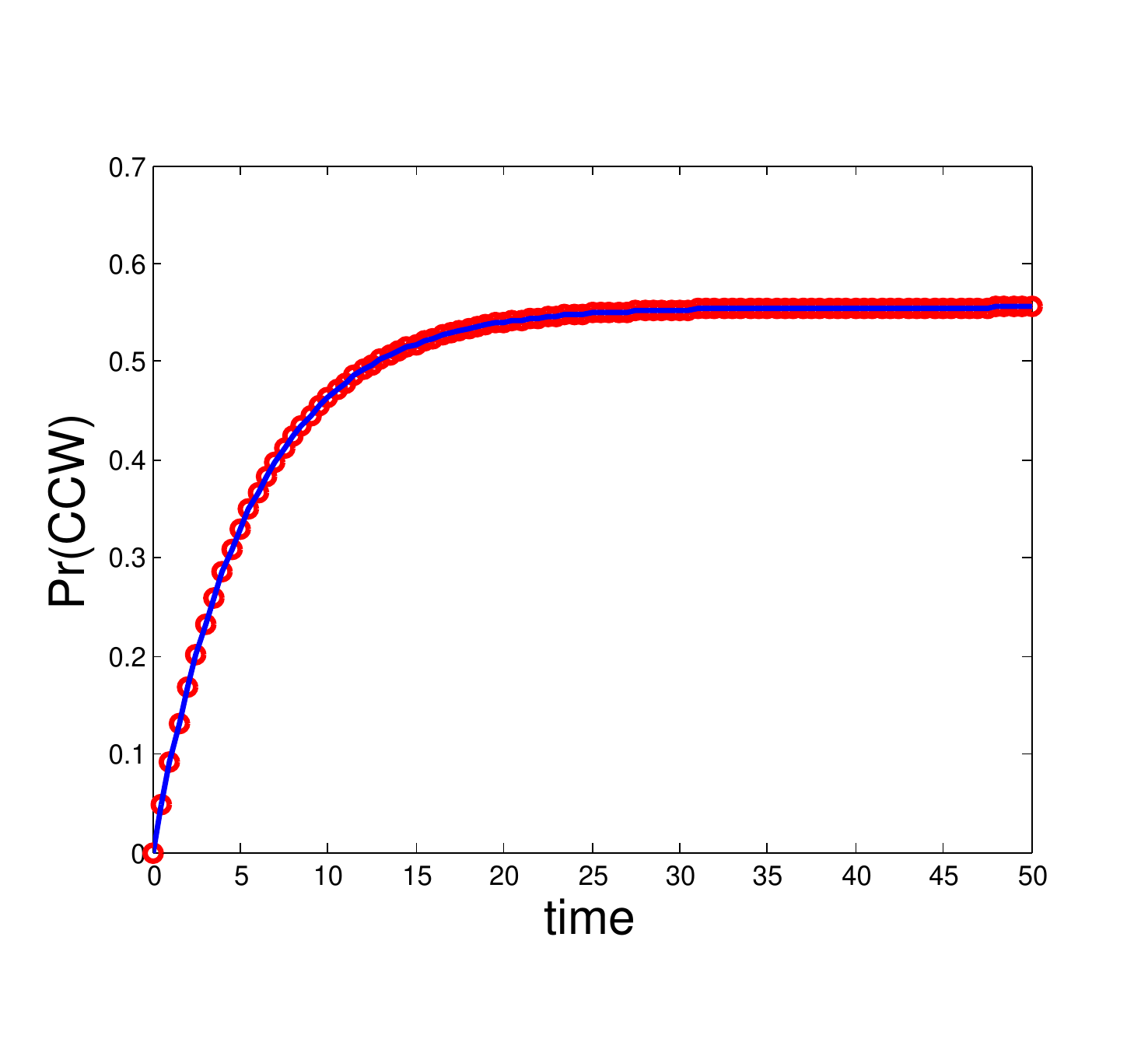}} \\
{\includegraphics[width=2.5in]{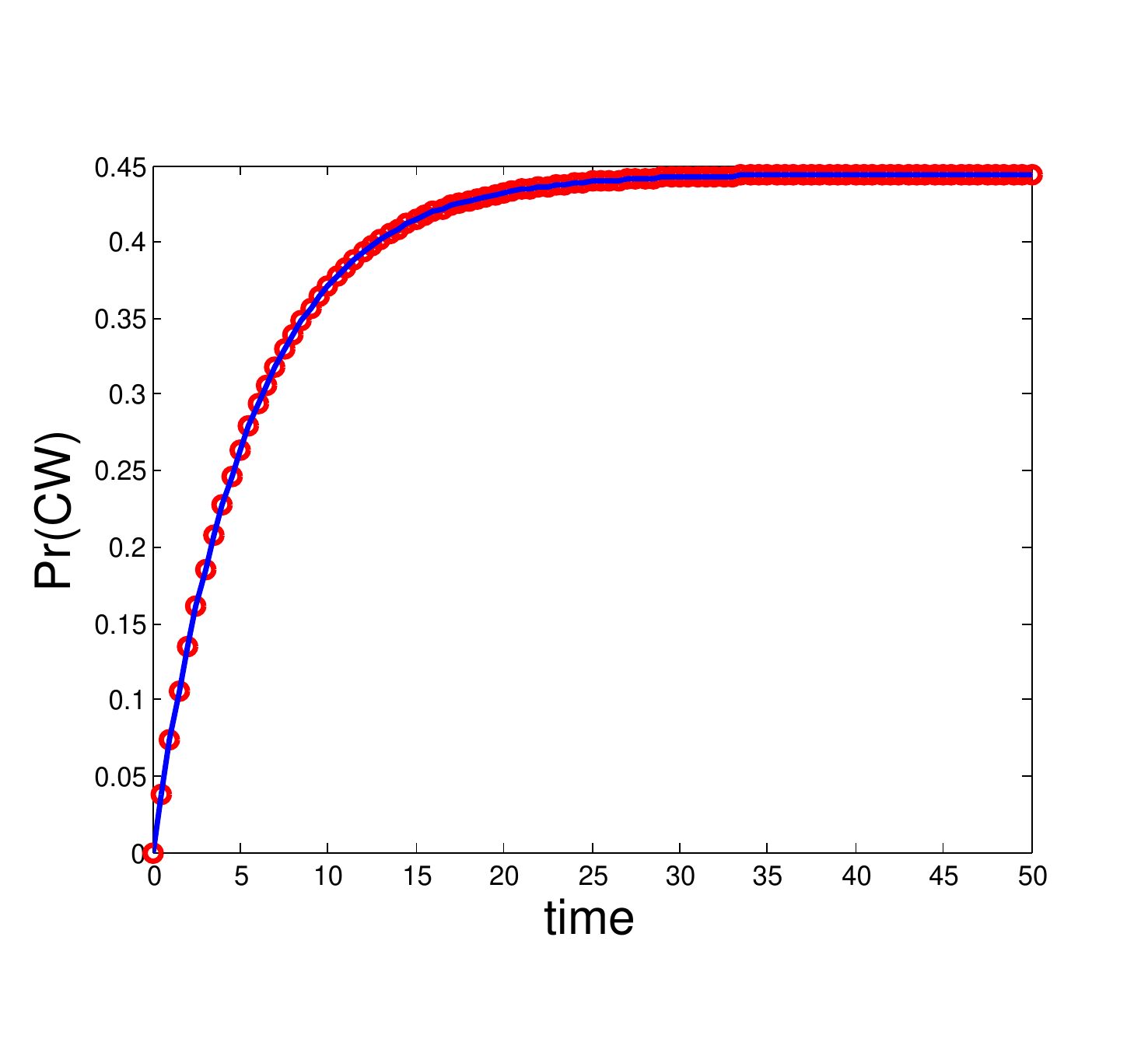}}  & {\includegraphics[width=2.5in]{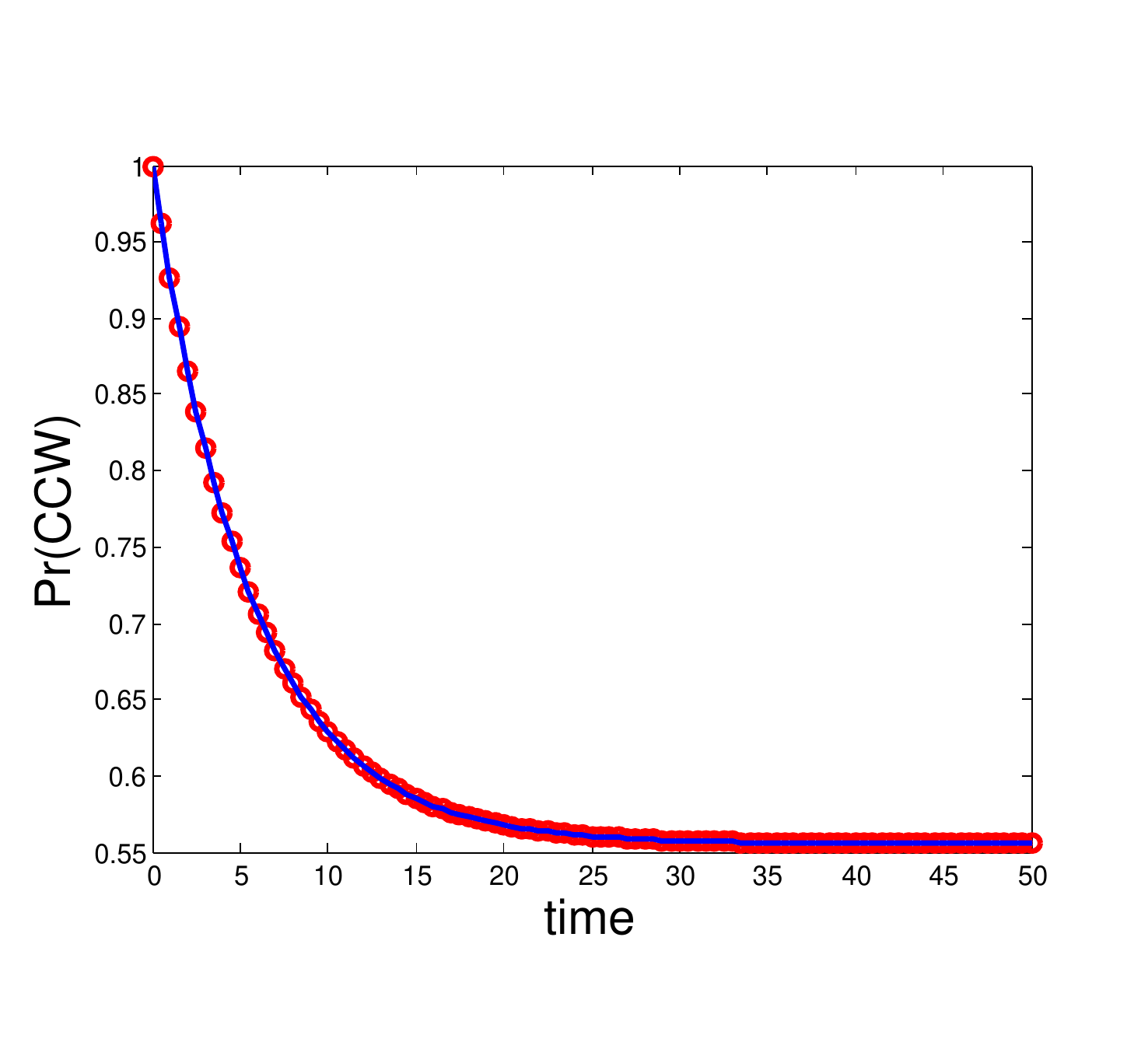}} 
\end{tabular}
\end{center}
\caption{\footnotesize Shown here are the time evolutions of the probability that the motor rotates clockwise ($Pr(CW) = 1$) and counterclockwise ($Pr(CCW) = 1$) with different initial conditions using $k_1=2, k_{-1}=1, k_3=2, k_{-3}=1$, $\alpha_i=0.1, \beta_i=0.08, Y=100$. In the upper two figures, the initial condition is
that the motor is rotating clockwise with no $Y$ bound. In the lower two figures, the initial condition is
that the motor is rotating counterclockwise with no $Y$ bound. Blue lines are generated by computing the full system ($Pr(CW) = \sum Pr(CW_i)$, $Pr(CCW) = \sum Pr(CCW_i)$), whereas red circles are for the reduced system.} 
\label{flagellum}
\end{figure}
\end{example}

\begin{example}[An enzyme inhibition model]

We consider an enzyme-substrate reaction system with a competitive inhibition 
\begin{eqnarray}
E+S \ \ssang{k_1}{k_2} \  ES
\stackrel{k_3}{\rightarrow} E+P, \ \ E+I \ \ssang{k_4}{k_5} \  EI,   \label{inhibit}
\end{eqnarray}
where $E, S, I$ and $P$ denote enzyme, substrate, inhibitor and product, respectively. 
Let $n_i(t),i=1,\dots,6$ denote the number of molecules of $E,S,ES,I,EI$ and $P$, respectively.
We assume that the two reversible reactions $E+S \ \ssang{}{} \  ES, \ \ E+I \ \ssang{}{} \  EI $ are much faster
than the irreversible reaction $ES {\rightarrow} E+P$. After finding 
$$A^f= \left[\begin{array}{cccccc} 
 1 & 0 & 1 & 0 & 1 & 0 \\
 0 & 1 & 1 & 0 & 0 & 0 \\
 0 & 0 & 0 & 1 & 1 & 0 \\
 0 & 0 & 0 & 0 & 0 & 1 
\end{array} \right], $$
one can identify the slow variable
$$\tilde n = A^f n = (n_1 + n_3 + n_5, n_2+ n_3, n_4 + n_5, n_6)^T.$$
Since the fast subsystem $E+S \ \ssang{}{} \  ES, \ \ E+I \ \ssang{}{} \  EI $ has 
a deficiency of zero and is weakly reversible, one can use the result by \cite{Anderson:2011:CTM} for finding the equilibrium probability of the fast subsystem;   
For the convenience of computation, we assume that $k_1=k_2$ and $k_3=k_4$. 
If the deterministic equilibrium values of $E,S,ES,I$ and $EI$ for the fast dynamics are denoted by $c_1, c_2, c_3, c_4$ and $c_5$, respectively,  
one can find the only one solution
\begin{eqnarray*}
c_3 &=& \frac{\alpha_2 (1+\alpha_1+\alpha_2+\alpha_3 +\sqrt{4\alpha_1 + (1-\alpha_1+\alpha_2+\alpha_3)^2)}}{2
(\alpha_2+\alpha_3)}, \\
 c_4 &=& \frac{\alpha_3 (1+\alpha_1+\alpha_2+\alpha_3 +\sqrt{4\alpha_1 + (1-\alpha_1+\alpha_2+\alpha_3)^2)}}{2
(\alpha_2+\alpha_3)}, 
\end{eqnarray*}
and $$c_1=\alpha_1 - c_3 -c_5, \ c_2=\alpha_2 - c_3, \ c_4 = \alpha_3-c_5,$$
where $\alpha_1 = c_1 + c_3 + c_5,  \alpha_2 = c_2 + c_3,$ and $\alpha_3 = c_4 + c_5$ are conserved quantities in the fast subsystem. Using the result of \cite{Anderson:2010:PFS}, one finds the equilibrium probability of the fast subsystem as 
$$p(n_1,\dots,n_5) =  M \prod_{i=1}^5 \frac{c_i^{n_i}}{n_i !}, n_i = 0, 1, 2,\dots$$ 
where $M$ is the normalizing constant for $\sum_n p(n) =1$ and $n_1,\dots, n_5$ satisfy the conserved quantities $n_1+n_3+n_5$, $n_2+n_3$ and $n_4+n_5$.  
The following figure shows the simulation results for the slow variable $\tilde
n_4$ (the number of product molecules) obtained from the approximate algorithm and the
exact algorithm. 
\begin{figure}[h!]
\begin{center}
{\includegraphics[width=2.5in]{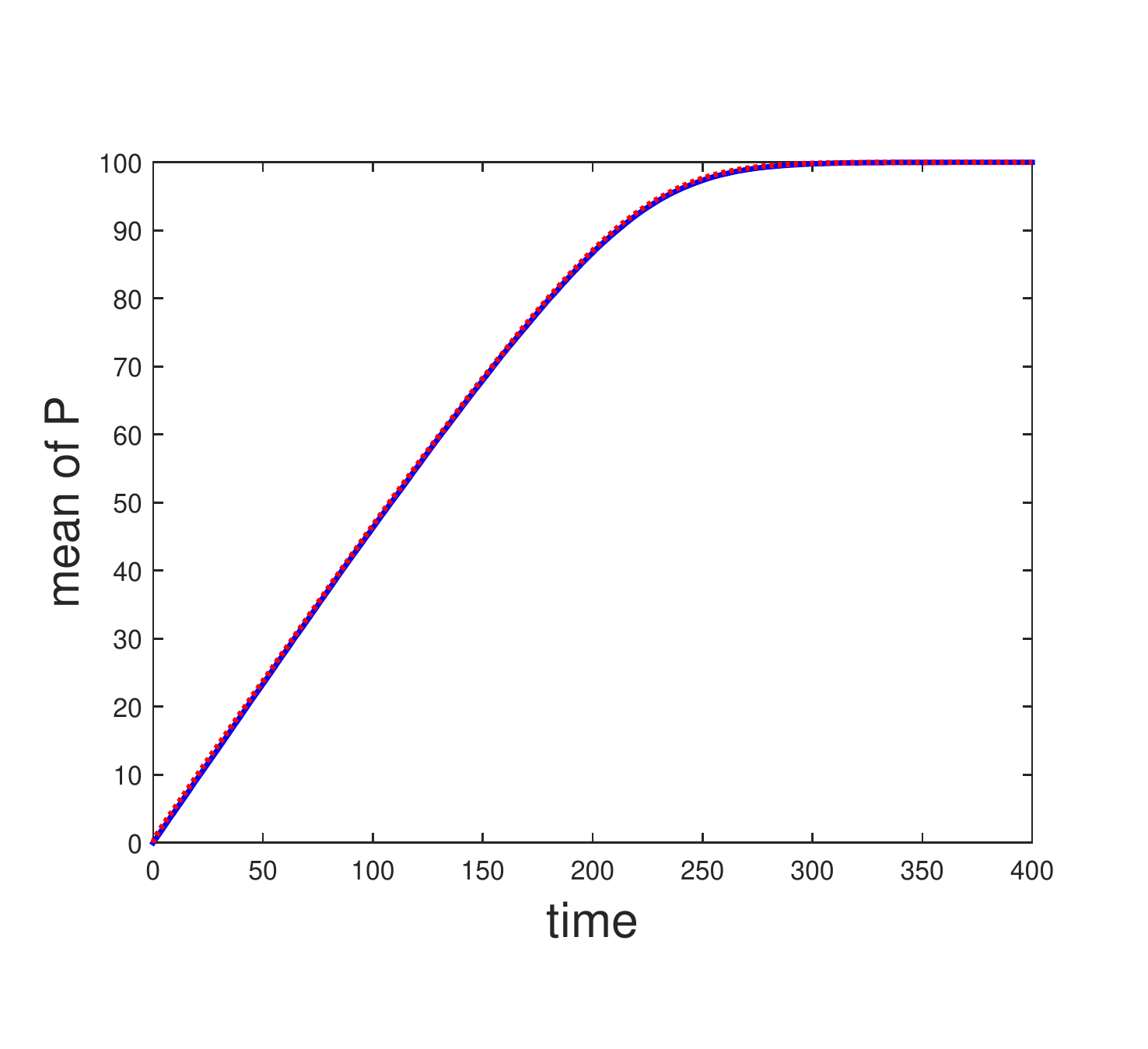}}  {\includegraphics[width=2.5in]{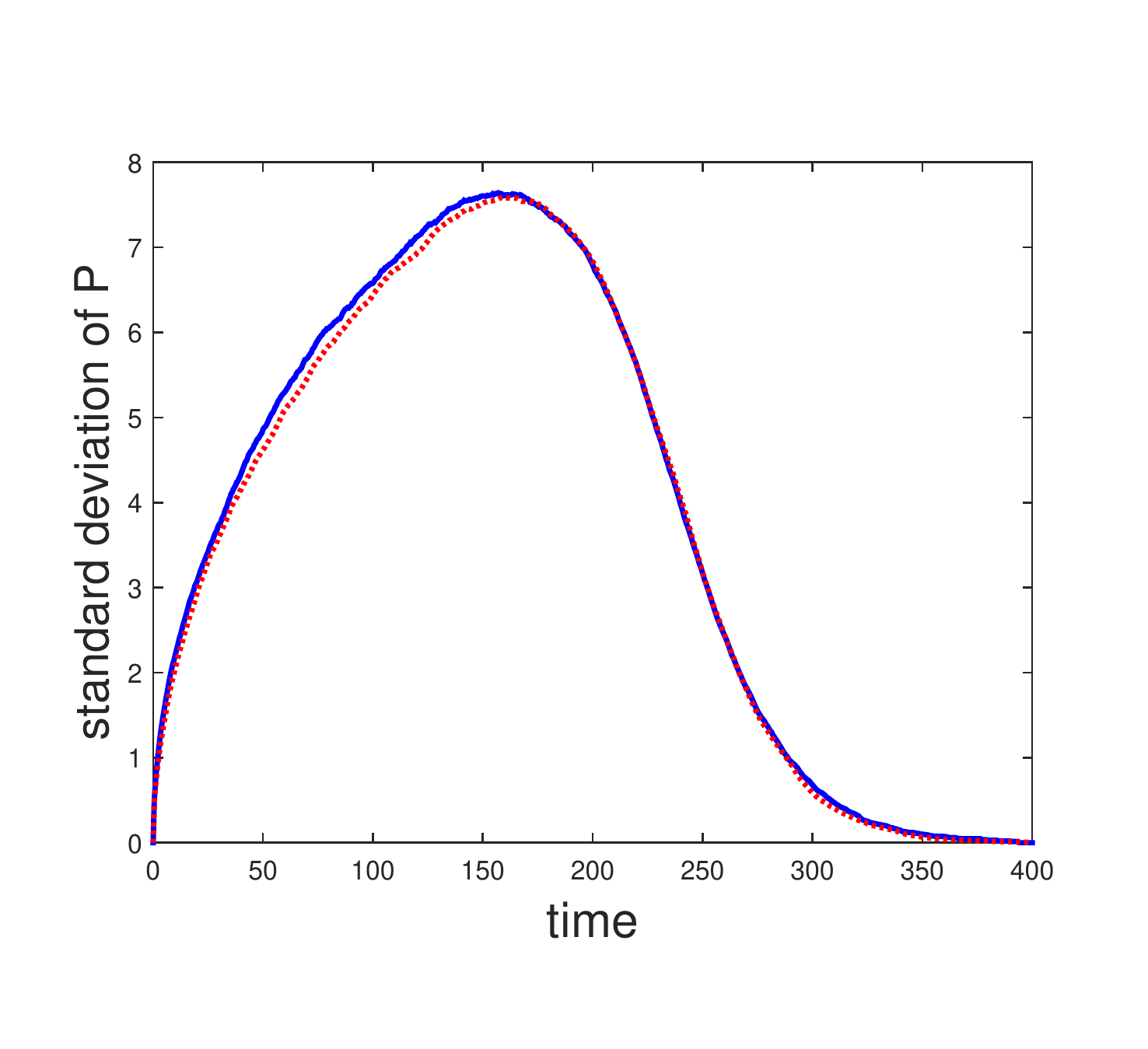}} \\
\label{detfig1}
  \caption{Enzyme-substrate model with an inhibitor. 
Comparison of approximate stochastic simulation algorithm (red dotted) to exact stochastic simulation algorithm (blue solid). Evolution of means and standard deviations of numbers of the product $P$ 
when the initial condition $(E, S, ES, I, EI, P) = (5, 100, 0, 5, 0, 0)$ and 
$k_1=k_2=k_4=k_5=10, k_3=0.1$.
The results are based on 5000 realizations. Concerning the relative CPU time for
one realization of the stochastic simulation, if the approximate algorithm takes
1 second,  the exact algorithm takes about 3.6 seconds, when a quad-core machine
with Windows 8.1 and MATLAB 2014 is used.} 
\end{center}
\end{figure}
\end{example}

\begin{example}[A model for the PFK system]

We consider a reaction network from a model of the  PFK  step in glycolysis \cite{Othmer:1978:ECD}
\begin{eqnarray*}
A_1+E_1 \ \ssang{k_1}{k_{-1}} \  & E_1 A_1 & \stackrel{k_2}{\to}  E_1+A_2 \\
A_1+E_1^\ast \ \ssang{k_3}{k_{-3}} \  & E_1^\ast A_1 & \stackrel{k_4}{\to}  E_1^\ast +A_2 \\
A_2+E_2 \ \ssang{k_5}{k_{-5}} \  & E_2 A_2 & \stackrel{k_6}{\to}  E_2+\textrm{\ Product}. \\
\end{eqnarray*}
Here $A_1, A_2, E_1, E_1^\ast$ and $E_2$ denote  $F6P, ADP$, the low activity
and activated forms of free PFK and the enzyme for the ADP sink reaction,
respectively, and  $E_1A_1$, $E_1^\ast A_1$ and $E_2A_2$ represent
enzyme-substrate complexes.  

If we assume three binding/unbinding reactions are much faster than others,
 there are  two different fast  subsystems $C_1$ and $C_2$ in the reaction network,
$$C_1 \ \ : \ \ A_1+E_1 \ \ssang{k_1}{k_{-1}} \   E_1 A_1, \ \ \ A_1+E_1^\ast \ \ssang{k_3}{k_{-3}} \  E_1^\ast A_1$$
and
$$C_2 \ \ : \ \ A_2+E_2 \ \ssang{k_5}{k_{-5}} \  E_2 A_2.$$
First we can find stationary distribution of $C_2$ using the hypergeometric
functions \cite{McQuarrie:1967:SAC}. 
Let initial numbers of $A_2$ and $E_2$ be $a_0$ and $b_0$, respectively and define $Q=\frac{k_{-1}}{k_1}$.  If $b_0 \ge a_0$, then
the stationary marginal distribution of $A_2$ is given by
$$P_{A_2}(k)= D \frac{Q^k}{k!} \frac{(a_0+c_0)(a_0+c_0-1)\cdots (a_0+c_0 - k+1)}{(b_0-a_0+1)(b_0-a_0+2) \cdots (b_0-a_0+k)},$$ 
where $k=1,\dots,a_0+c_0$ and $D$ is a normalization constant.

Similarly, if $a_0 > b_0$, then
the stationary marginal distribution of $E_2$
$$P_{E_2}(k)= D \frac{Q^k}{k!} \frac{(b_0+c_0)(b_0+c_0-1)\cdots (b_0+c_0 - k+1)}{(a_0-b_0+1)(a_0-b_0+2) \cdots (a_0-b_0+k)},$$ 
where $k=1,\dots,b_0+c_0$

Since the fast subsystem $C_1$ has a deficiency of zero and is weakly reversible, we can find the equilibrium probability similar to the enzyme inhibitor model;
If  we assume that $k_1=k_2$ and $k_3=k_4$ and the deterministic equilibrium values of $A_1,E_1, E_1A_1,E_1^\ast, E_1^\ast A_1$ are denoted by $c_i, i=1,...,5$, respectively,   then
we can find the equilibrium probability of the fast subsystem $C_1$  
$$p(n_1,\dots,n_5) =  M \prod_{i=1}^5 \frac{c_i^{n_i}}{n_i !}, n_i = 0, 1, 2,\dots$$ 
where $n_i,i=1,...,5$ are the number of $A_1,E_1,E_1 A_1, E_1^*,E_1^*A_1$, subject to the conserved quantities $n_1+n_3+n_5$, $n_2+n_3$ and $n_4+n_5$,
$M$ is the normalizing constant for $\sum_n p(n) =1$ and $n_1,\dots, n_5$ satisfy the conserved
  quantities $n_1+n_3+n_5$, $n_2+n_3$ and $n_4+n_5$.   
\begin{eqnarray*}
c_3 &=& \frac{\alpha_2 (1+\alpha_1+\alpha_2+\alpha_3 +\sqrt{4\alpha_1 + (1-\alpha_1+\alpha_2+\alpha_3)^2)}}{2
(\alpha_2+\alpha_3)}, \\
 c_4 &=& \frac{\alpha_3 (1+\alpha_1+\alpha_2+\alpha_3 +\sqrt{4\alpha_1 + (1-\alpha_1+\alpha_2+\alpha_3)^2)}}{2
(\alpha_2+\alpha_3)}, 
\end{eqnarray*}
$$c_1=\alpha_1 - c_3 -c_5, \ c_2=\alpha_2 - c_3, \ c_4 = \alpha_3-c_5,$$
and $\alpha_1 = c_1 + c_3 + c_5,  \alpha_2 = c_2 + c_3$ and $\alpha_3 = c_4 +
c_5$ are conserved quantities in the fast subsystem.  

The following figure shows the simulation results for the number of product, which is a slow variable.
\begin{figure}[h!]
\begin{center}
\begin{tabular}{cc}
{\includegraphics[width=3.in]{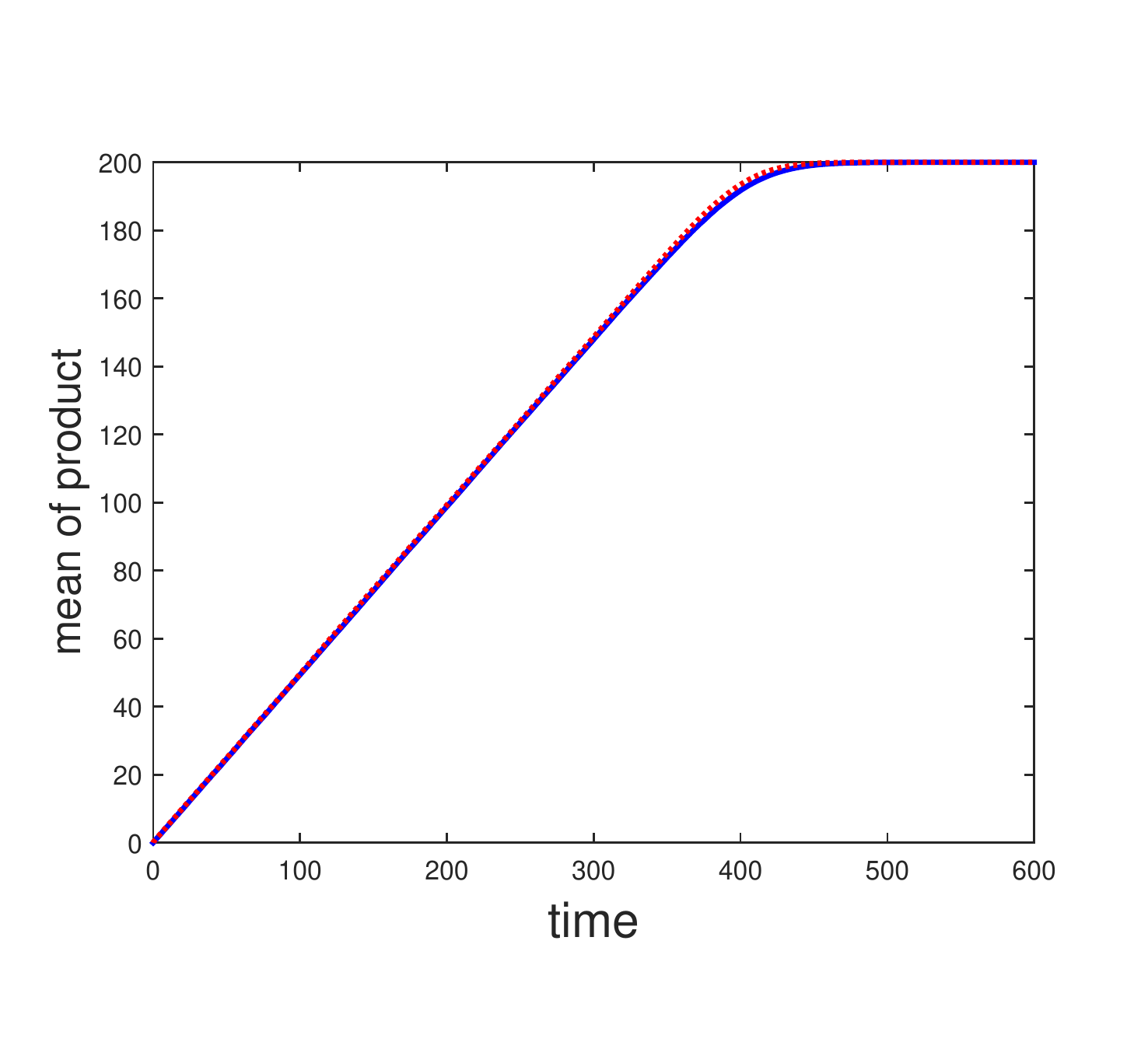}}  &
{\includegraphics[width=3.in]{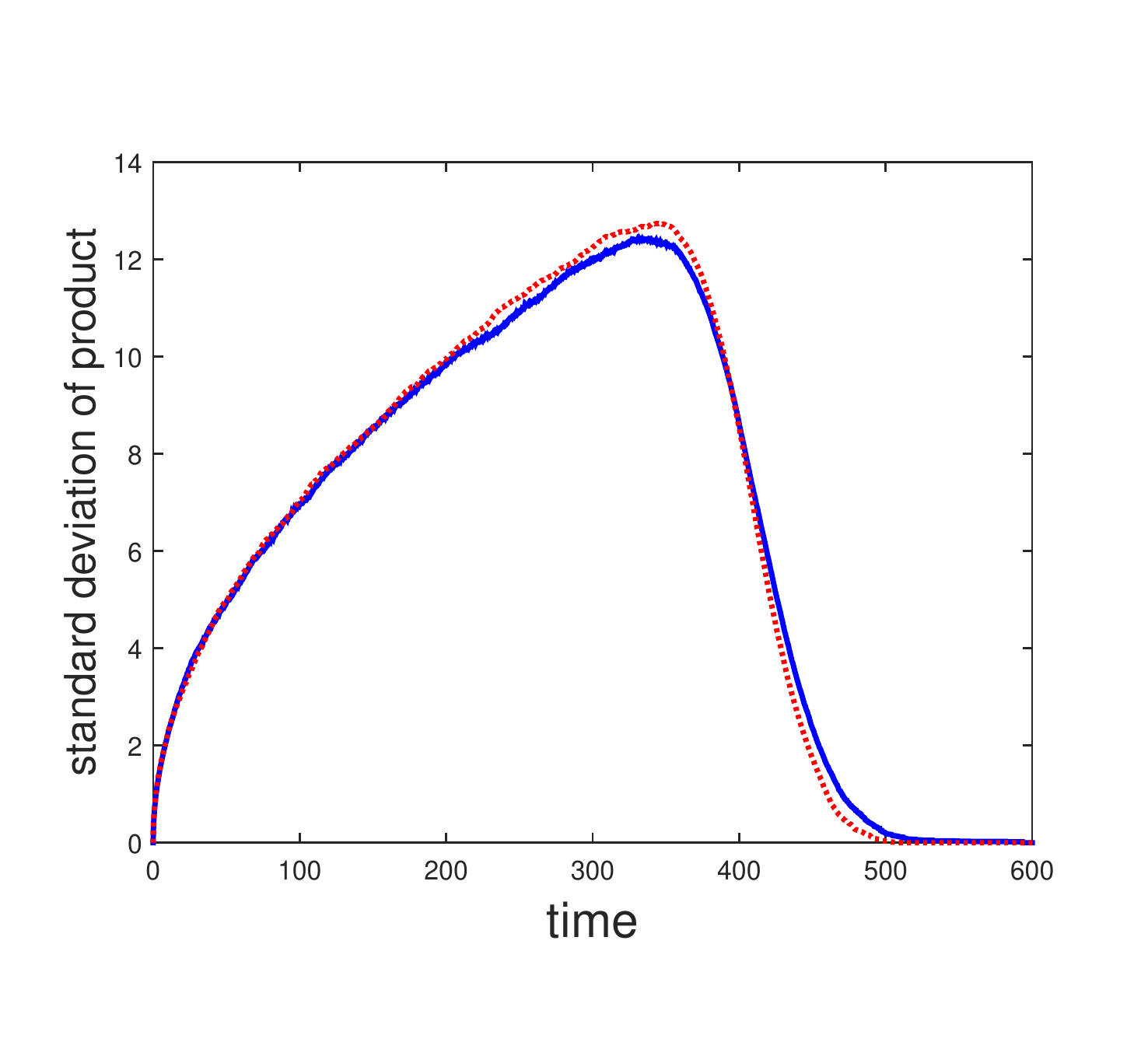}} 
\label{detfig2} \\
\end{tabular}
 \caption{ 
A comparison of the approximate stochastic simulation algorithm (red dotted) to
the exact stochastic simulation algorithm (blue solid) for the model of PFK reaction system. Evolution of means and standard deviations of numbers of the product when the initial condition $(A_1, E_1, E_1 A_1, E_1^*, E_1^* A_1, A_2, E_2, E_2 A_2, P) = (100, 5, 0, 5, 0, 100, 5, 0, 0)$ and $k_2=k_4=k_{6}=0.1,  k_1=k_{-1}=k_3=k_{-3}=k_5=k_{-5}=10$.
The results are based on 5000 realizations. Concerning the relative CPU time for
one realization of the stochastic simulation, if the approximate algorithm takes
1 second,  the exact algorithm takes about 2.9 seconds, when a quad-core machine
with Windows 8.1 and MATLAB 2014 is used. } 
\end{center}
\end{figure}

\end{example}

\section{Conclusion}

We developed a reduction method for stochastic biochemical reaction networks
with coupled fast and slow reactions, and formulated an associated extension of
the Gillespie method. The reduced equation for the time-dependent probability
distribution on the slow time scale involves a projection in state space, and we
show that the generator of the reduced system is also Markovian.  We found that
the transition rate for a reaction on the slow time scale is approximated by the
expectation of the original transition rate conditional on the invariant distribution of the fast dynamics. Throughout the reduction the invariants of fast subsystems play a
significant role, similar to what was found earlier in the analogous reduction
of deterministic systems \cite{Lee:2009:MTS}.  When the stationary distribution
of the fast dynamics can be computed, these can be used directly in a modified
stochastic simulation algorithm.  Several biological examples, including a model
of motor behavior for a single flagellum, an enzyme-inhibitor model and a model
for the PFK system, illustrate the numerical accuracy of the approximation for
two lowest moments.

\section{Appendix}

\subsection{The general case in Example \ref{Triex}}
\label{subsec:reduc}

When there are $N_0$ molecules, there are $N_0+1$ fast components and the
transition matrix $K$ is block  
tridiagonal, i.e.

$$ K = \left[\begin{array}{cccccc} \dfrac{1}{\epsilon}K_1^f+K_1^s & K_{1,2}^s & & & &\\
K_{2,1}^s & \dfrac{1}{\epsilon}K_2^f+K_2^s & K_{2,3}^s & & &\\
 & K_{3,2}^s & \ddots & \ddots & & \\
 & & \ddots & \ddots & \ddots & \\
 & & & \ddots & \ddots & K_{N_0,N_0+1}^s\\
 & & & & K_{N_0+1,N_0}^s & \dfrac{1}{\epsilon}K_{N_0+1}^f+K_{N_0+1}^s
\end{array}
\right],
  $$

where  the fast  blocks are given by 
\footnotesize
$$ 
K_1^f = \left[\begin{array}{cccccc} -N_0k_1 & k_2 & & & &\\
N_0k_1 & -[(N_0-1)k_1+k_2] & 2k_2 & & &\\
 & (N_0-1)k_1 & -[(N_0-2)k_1+2k_2] & 3k_2 &  &  \\
 & & \ddots & \ddots & \ddots & \\
 & & & \ddots & \ddots & N_0k_2\\
 & & & & k_1 & -N_0k_2
\end{array}
\right]_
{(N_0+1)\times(N_0+1)}\\
$$
$$
K_2^f = \left[\begin{array}{cccccc} -(N_0-1)k_1 & k_2 & & & &\\
(N_0-1)k_1 & -[(N_0-2)k_1+k_2] & 2k_2 & & &\\
 & (N_0-2)k_1 & -[(N_0-3)k_1+2k_2] & 3k_2 &  &  \\
 & & \ddots & \ddots & \ddots & \\
 & & & \ddots & \ddots & (N_0-1)k_2\\
 & & & & k_1 & -(N_0-1)k_2
\end{array}
\right]_
{N_0\times N_0}
\cdots\\
$$
$$
K_{N_0}^f = \left[\begin{array}{cc} -k_1 & k_2\\
 k_1 & -k_2
  \end{array}
  \right],
  K_{N_{0}+1}^f = 0.
$$

\normalsize
The slow off-diagonal blocks are given by 
\footnotesize
$$
K_{1,2}^s=\left[\begin{array}{cccc} k_5 & & & \\
k_4 & k_5 & & \\
& \ddots & \ddots & \\
& & \ddots & k_5\\
& & & k_4
\end{array}\right]_{(N_0+1)\times N_0},
K_{2,3}^s=\left[\begin{array}{cccc} 2k_5 & & & \\
2k_4 & 2k_5 & & \\
& \ddots & \ddots & \\
& & \ddots & 2k_5\\
& & & 2k_4
\end{array}\right]_{N_0\times (N_0-1)},
\cdots,\\
K_{N_{0},N_{0}+1}^s = \left[\begin{array}{c}N_0k_5\\
N_0k_4
\end{array}\right].
$$
\normalsize
For the lower diagonal blocks,
\footnotesize

$$
K_{2,1}^s=\left[\begin{array}{ccccc} N_0k_6 & k_3 & & & \\
 & (N_0-1)k_6 & 2k_3 & &\\
&  & \ddots & \ddots &\\
& & & k_6 & N_0k_3
\end{array}\right]_{N_0\times (N_0+1)} \quad
K_{3,2}^s=\left[\begin{array}{ccccc} (N_0-1)k_6 & k_3 & & & \\
 & (N_0-2)k_6 & 2k_3 & &\\
&  & \ddots & \ddots &\\
& & & k_6 & (N_0-1)k_3
\end{array}\right]_{(N_0-1)\times N_0},
$$
$\cdots, \qquad 
K_{N_0+1,N_0}^s = \left[\begin{array}{cc}k_6 & k_3
\end{array}\right]
$. 
\normalsize
Finally, the $K_i^s$ along the diagonal are diagonal matrices of the same
 dimension as $K_i^f$, and the $(j,j)$-th entry of $K_i^s$ is the negative sum
 of rates leaving $j$-th node of the $i$-th fast component.  

In this case the transition
 rate is given by 
$$
\tilde k_{i,j}^s=L_iK_{i,j}^s\Pi_j.
$$

\subsection{Moment equations of the invariant distributions} 

The low-order moments of the distributions for the fast systems play a role in
the QSS reduction in Section \ref{master_reduc}, and here we consider the
low-order moment equations.
\begin{Theorem}
\label{steadyEQ}
Let $r$ be the total number of the reactions in the system. Then the invariant
(steady-state) distribution  of $P(n,t)$, which we denote $P(n)$, satisfies 
\begin{eqnarray}
\sum_{\ell=1}^{r}  \cR_{\ell}(n-\nu\cE_{(\ell)} P(n-\nu\cE_{(\ell)}) = \sum_{\ell=1}^{r}  \cR_{\ell}(n) P(n) \label{QEprob}
\end{eqnarray}
and 
and the first two moment equations  lead to 
\begin{eqnarray}
&& \hspace{-0.2cm} E[  \nu\cE  \cR(n)]=0, \label{QEmean}\\
&&  \hspace{-0.2cm} \sum_{i=1}^{r}\big[  \nu\cE_{(i)} \otimes E[n \cR_i(n)] +  E[n \cR_i(n)] \otimes \nu\cE_{(i)} +  \nu\cE_{(i)} \otimes \nu\cE_{(i)} E[ \cR_i(n)] \big] = 0. \nonumber \\
\label{QEsecond}
\end{eqnarray}
\end{Theorem}
\begin{proof}   At the steady-state  
\begin{eqnarray}
\sum_{\ell=1}^{r} \cR_{\ell}(n-\nu\cE_{(\ell)}) P(n-\nu\cE_{(\ell)}) = \sum_{\ell=1}^{r} \cR_{\ell}(n) P(n) \label{quasieq}
\end{eqnarray}
By multiplying by $n$ and summing over all the values of $n\in\cL(n_0)$, we obtain
$$\sum_{n}\sum_{\ell=1}^{r} n \cR_{\ell}(n-\nu\cE_{(\ell)}) P(n-\nu\cE_{(\ell)}) = \sum_{n}\sum_{\ell=1}^{r} n \cR_{\ell}(n) P(n).$$
Using the transformation $ n-\nu\cE_{(\ell)} \to n$ on the left side,
we obtain
$$\sum_{n}\sum_{\ell=1}^{r} (n+\nu\cE_{(\ell)}) \cR_{\ell}(n) P(n) = \sum_{n}\sum_{\ell=1}^{r} n \cR_{\ell}(n) P(n).$$
By subtracting the right side from the left one,
$$\sum_{\ell=1}^{r} \nu\cE_{(\ell)} \sum_{n}  \cR_{\ell}(n) P(n) = \sum_{\ell=1}^{r} \nu\cE_{(\ell)} E[\cR_{\ell}(n)]= 0.$$
Thus we conclude that
$$\nu\cE E[ \cR(n)]= 0.$$
If the deficiency $\delta \equiv  \rho (\cE) -  \rho (\nu\cE)$ is zero, then
$E[ \cR(n)]$ is a cycle in the graph \cite{Othmer:1979:GTA}.\footnote{The reader
can show that the presence of inputs or outputs of the form given in Table
\ref{rxntypes} does not alter the deficiency.}

At the next order  we multiply equation (\ref{quasieq}) by a
tensor product $n\otimes n$ and sum over $n$. Then by a
similar argument, we obtain 
\begin{eqnarray*}
 && \sum_{i=1}^{r}\big[ \nu\cE_{(i)} \otimes \sum_{n}  n \cR_i(n) p(n) +  \sum_{n}  n \cR_i(n) p(n)\otimes \nu\cE_{(i)}  +  \nu\cE_{(i)} \otimes \nu\cE_{(i)} \sum_{n}  \cR_i(n) p(n) \big]\\
&=& \sum_{i=1}^{r}\big[  \nu\cE_{(i)} \otimes E[n \cR_i(n)] +  E[n \cR_i(n)] \otimes \nu\cE_{(i)} +  \nu\cE_{(i)} \otimes \nu\cE_{(i)} E[ \cR_i(n)] \big]\\
&=& 0.
\end{eqnarray*} 
\end{proof}
If  $\delta =0$, the two lowest moment equations can be  simplified to 
\begin{eqnarray}
&& E[ \cR(n)]=0 \\
&&  \sum_{i=1}^{r}\big[  \nu\cE_{(i)} \otimes E[n \cR_i(n)] +  E[n \cR_i(n)] \otimes \nu\cE_{(i)} \big] = 0.
\end{eqnarray}
When all reactions are linear the problem is much simpler, and the evolution
equations for the first and second moments can be written {\em explicitly in
terms of those moments} \cite{Gadgil:2005:SAF}. 

As a consequence of  Theorem \ref{steadyEQ}, similar equations can be obtained
for the quasi-steady-state of the probability distribution for the fast
subsystem in a two-time scale stochastic network. We first define the
expectation of a function $f(n)$ over a discrete reaction simplex  $\cL_f$ for
the fast subsystem as follows. 
$$E_{\cL_f}[f(n)] \equiv \sum_{n\in \cL_f} f(n) p(n).$$
\begin{Corollary}
Let $r_f$ be the total number of fast reactions.
Then at the steady-state of the fast subsystem,
the governing equation is given by 
\begin{eqnarray}
\sum_{i=1}^{r_f}  \cR^f_i(n-\nu\cE^f_{(i)}) P(n-\nu\cE^f_{(i)}) =
\sum_{i=1}^{r_f}  \cR^f_i(n) P(n) 
\label{fastequipro} 
\end{eqnarray}
and  for each discrete reaction simplex $\cL_f$,
\begin{eqnarray}
&& E_{\cL_f} [  \cR^f(n)]=0, \\
&&  \sum_{i=1}^{r_f}\big[  \nu\cE^f_{(i)} \otimes E_{\cL_f}[n \cR^f_{i}(n)] + E_{\cL_f}[n \cR^f_i(n)] \otimes \nu\cE^f_{(i)}   \big] = 0.
\end{eqnarray}
\end{Corollary}

\begin{proof}
The `state-wise' form  of the master  equation (\ref{prob_evol1})  can be
written 
\begin{eqnarray}
\dfrac{d}{dt}P(n,t) &= & \sum_{\ell}\dfrac{1}{\ep} \big[\cR^f_{\ell}(n-\nu\cE^f_{(\ell)})\cdot P(n-\nu\cE^f_{(\ell)},t)-\cR^f_{\ell}(n)\cdot P(n,t)\big] \nonumber\\ 
\label{twotimemaster} \\[-15pt]
&& + \sum_{k} \big[\cR^s_{k}(n-\nu\cE^s_{(k)})\cdot P(n-\nu\cE^s_{(k)},t)- \cR^s_{k}(n)\cdot
P(n,t)\big], \nonumber
\end{eqnarray}
where $\cR^f$ and $\cR^s$ are the transition rates of fast and slow reactions,
respectively and $\cE^f$ and $\cE^s$ are incidence matrices for fast and slow
reactions, respectively. 

 In the previous theorem, substitute $\cE^f, \cR^f$ and
  $E_{\cL_f}[\cdot]$   into  $\cE, \cR$ and $E[\cdot]$ and use the  full rank
  assumption on $ \nu\cE^f$.
 \end{proof}
%


\begin{thebibliography}{}

\bibitem[Anderson {\em et~al.}, 2010]{Anderson:2010:PFS}Anderson, D.~F.,
 Craciun, G., \& Kurtz, T.~G. (2010) Product-form stationary distributions for deficiency zero chemical reaction networks. {\em Bull. Math. Biol.} {\bf 72}:1947--1970
 
\bibitem[Anderson \& Kurtz, 2011]{Anderson:2011:CTM}
 Anderson, D.~F. \& Kurtz, T.~G. (2011) Continuous time markov chain models for
  chemical reaction networks In: {\em Design and analysis of biomolecular
  circuits} pp. 3--42 Springer

\bibitem[Aris, 1965]{Aris:1965:PRA}
 Aris, R. (1965) Prolegomena to the rational analysis of chemical reactions.
  {\em Arch. Rat. Mech. Anal.} {\bf 19} (2):81--99

\bibitem[Boucherie \& Dijk, 1991]{Boucherie:1991PFQ}
 Boucherie, R.~J. \& Dijk, N. M.~V. (1991) Product forms for queueing networks
  with state-dependent multiple job transitions. {\em Advances in Applied
  Probability, } {\bf }:152--187

\bibitem[Bundschuh {\em et~al.}, 2003]{Bundschuh:2003:FSV}
 Bundschuh, R., Hayot, F., \& Jayaprakash, C. (2003) Fluctuations and slow
  variables in genetic networks. {\em Biophy. J.} {\bf 84}:1606

\bibitem[Campbell \& Meyer, 1991]{Campbell:1991:GIL}
 Campbell, S.~L. \& Meyer, C.~P. (1991) {\em Generalized Inverses of Linear
  Transformations} Dover

\bibitem[Cao {\em et~al.}, 2005]{Cao:2005:SSS}
 Cao, Y., Gillespie, D.~T., \& Petzold, L.~R. (2005) The slow-scale stochastic
  simulation algorithm. {\em J. Chem. Phys.} {\bf 122}:014116--

\bibitem[Chen, 1971]{Chen:1971:AGT}
 Chen, W.~K. (1971) {\em Applied Graph Theory} Amsterdam: North-Holland

\bibitem[Chevalier \& EI-Samad, 2009]{Chevalier:2009}
 Chevalier, M.~W. \& EI-Samad, H. (2009) A rigorous framework for multiscale
  simulation of stochastic cellular networks. {\em The Journal of Chemical
  Physics, } {\bf 131} (5):054102--
  
\bibitem[Cotter, 2015]{Cotter:2015:ARX} Cotter, Simon. (2015) Constrained Approximation of Effective Generators for Multiscale Stochastic Reaction Networks and Application to Conditioned Path Sampling.{\em arXiv preprint arXiv:1506.02446}.


\bibitem[Deuflhard {\em et~al.}, 2008]{Deuflhard:2008:ADG}
 Deuflhard, P., Huisinga, W., Jahnke, T., \& Wulkow, M. (2008) Adaptive
  discrete galerkin methods applied to the chemical master equation. {\em SIAM
  J. SCI. COMPUT.} {\bf 30} (6):2990--3011

\bibitem[E {\em et~al.}, 2005]{E:2005:NSS}
 E, W., Liu, D., \& Vanden-Eijnden, E. (2005) Nested stochastic simulation
  algorithm for chemical kinetic systems with disparate rates. {\em J. Chem.
  Phys.} {\bf 123}:194107

\bibitem[Gadgil {\em et~al.}, 2005]{Gadgil:2005:SAF}
 Gadgil, C., Lee, C.~H., \& Othmer, H.~G. (2005) A stochastic analysis of
  first-order reaction networks. {\em Bull. Math. Biol.} {\bf 67}:901--946

\bibitem[Gillespie, 2007]{Gillespie:2007:SSC}
 Gillespie, D.~T. (2007) Stochastic simulation of chemical kinetics. {\em Annu.
  Rev. Phys. Chem.} {\bf 58}:35--55

\bibitem[Goutsias, 2005]{Goutsias:2005:QAF}
 Goutsias, J. (2005) Quasiequilibrium approximation of fast reaction kinetics
  in stochastic biochemical systems. {\em J. Chem. Phys.} {\bf 122}:184102--

\bibitem[Goutsias \& Jenkinson, 2013]{Goutsias:2013}
 Goutsias, J. \& Jenkinson, G. (2013) Markovian dynamics on complex reaction
  networks. {\em Physics Reports, } {\bf 529}:199--264

\bibitem[Haseltine \& Rawlings, 2002]{Haseltine:2002:ASC}
 Haseltine, E.~L. \& Rawlings, J.~B. (2002) Approximate simulation of coupled
  fast and slow reactions for stochastic chemical kinetics. {\em Journal of
  Chemica Physics, } {\bf 117} (15):6959--6969

\bibitem[Hellander \& L{\"o}tstedt, 2007]{Hellander:2007}
 Hellander, A. \& L{\"o}tstedt, P. (2007) Hybrid method for the chemical master
  equation. {\em Journal of Computational Physics, } {\bf 227}:100--122

\bibitem[Horn \& Jackson, 1972]{Horn:1972:GMA}
 Horn, F. \& Jackson, R. (1972) General mass action kinetics. {\em Arch. Rat.
  Mech. Anal.} {\bf 48}:81

\bibitem[Hu {\em et~al.}, 2013]{Hu:2013:SAR}
 Hu, J., Kang, H.-W., \& Othmer, H.~G. (2013) Stochastic analysis of
  reaction--diffusion processes. {\em Bulletin of Mathematical Biology, } {\bf
  }:1--41

\bibitem[Huang \& Liu, 2014]{Huang:2014}
 Huang, C. \& Liu, D. (2014) Strong convergence and speed up of nested
  stochastic simulation algorithm. {\em Commun. Comput. Phys.} {\bf 15}
  (4):1207--1236

\bibitem[Jahnke \& Huisinga, 2007]{Jahnke:2007:SCM}
 Jahnke, T. \& Huisinga, W. (2007) Solving the chemical master equation for
  monomolecular reaction systems analytically. {\em J. Math. Biol.} {\bf
  54}:1--26

\bibitem[Janssen, 1989a]{Janssen:1989:EFVII}
 Janssen, J. A.~M. (1989a) The elimination of fast variables in complex
  chemical reactions. {II}.. {\em J. Stat. Phys.} {\bf 57}:171--186

\bibitem[Janssen, 1989b]{Janssen:1989:EFVIII}
 Janssen, J. A.~M. (1989b) The elimination of fast variables in complex
  chemical reactions. {III}. mesoscopic level(irreducible case). {\em J. Stat.
  Phys.} {\bf 57}:187--198

\bibitem[Kato, 1966]{Kato:1966:PTL}
 Kato, T. (1966) Perturbation {Theory} for {Linear Operators}. {\em
  Spring{\-}er-Ver{\-}lag, } {\bf }

\bibitem[Kazeev {\em et~al.}, 2014]{Kazeev:2014:DSC}
 Kazeev, V., Khammash, M., Nip, M., \& Schwab, C. (2014) Direct solution of the
  chemical master equation using quantized tensor trains. {\em PLoS Comput
  Biol, } {\bf }:e1003359

\bibitem[Kim {\em et~al.}, 2014]{Kim:2014}
 Kim, J.~K., Josic, K., \& Bennett, M.~R. (2014) The validity of
  quasi-steady-state approximations in discrete stochastic simulations. {\em
  Biophysical Journal, } {\bf 107}:783--793

\bibitem[Lee \& Lui, 2009]{Lee:2009}
 Lee, C.~H. \& Lui, R. (2009) A reduction method for multiple time scale
  stochastic reaction networks. {\em J Math Chem, } {\bf 46}:1292--1321

\bibitem[Lee \& Othmer, 2009]{Lee:2009:MTS}
 Lee, C.~H. \& Othmer, H.~G. (2009) {A multi-time-scale analysis of chemical
  reaction networks: I. Deterministic systems}. {\em Journal of Mathematical
  Biology, } {\bf 60} (3):387--450

\bibitem[Mastny {\em et~al.}, 2007]{Mastny:2007}
 Mastny, E.~A., Haseltine, E.~L., \& Rawlings, J.~B. (2007) Two classes of
  quasi-steady-state model reductions for stochastic kinetics. {\em The Journal
  of Chemical Physics, } {\bf 127} (9):094106--

\bibitem[McQuarrie, 1967]{McQuarrie:1967:SAC}
 McQuarrie, D.~A. (1967) Stochastic approach to chemical kinetics. {\em J
  Applied Probability, } {\bf 4} (3):413--478

\bibitem[M{{\'e}}lyk{{\'u}}ti {\em et~al.}, 2014]{Melykuti:2014:EDS}
 M{{\'e}}lyk{{\'u}}ti, B., Hespanha, J.~P., \& Khammash, M. (2014) Equilibrium
  distributions of simple biochemical reaction systems for time-scale
  separation in stochastic reaction networks. {\em Journal of The Royal Society
  Interface, } {\bf 11} (97)

\bibitem[Menz {\em et~al.}, 2012]{Menz:2012:HSD}
 Menz, S., Latorre, J.~C., Schutte, C., \& Huisinga, W. (2012) Hybrid
  stochastic-deterministic solution of the chemical master equation. {\em
  Multiscal Model Simul.} {\bf 10} (4):1232--1262
  
  
\bibitem[Norris, 1998]{Norris:1998}
Norris, J.~R. (1998) Markov chains. {\em University of Cambridge}

\bibitem[Othmer, 1976]{Othmer:1976:NEC}
 Othmer, H.~G. (1976) Nonuniqueness of equilibria in closed reacting systems.
  {\em Chemical engineering science, } {\bf 31}:993--1003

\bibitem[Othmer, 1979]{Othmer:1979:GTA}
 Othmer, H.~G. (1979) A graph-theoretic analysis of chemical reaction networks
  Lecture Notes, Rutgers University -- available at 
  \url{http://math.umn.edu/~othmer/graphrt.pdf}

\bibitem[Othmer, 2005]{Othmer:2005:ACR}
 Othmer, H.~G. (2005) Analysis of complex reaction networks Lecture Notes,
  University of Minnesota

\bibitem[Othmer \& Aldridge, 1978]{Othmer:1978:ECD}
 Othmer, H.~G. \& Aldridge, J.~A. (1978) The effects of cell density and
  metabolite flux on cellular dynamics. {\em J. Math. Biol.} {\bf 5}:169--200

\bibitem[Peles {\em et~al.}, 2006]{Samad:2006:SHS}
 Peles, S., Munsky, B., \& Khammash, M. (2006) Reduction and solution of the
  chemical master equation using time scale separation and finite state
  projection. {\em J. Chem. Phys.} {\bf 125}:204104--

\bibitem[Rao \& Arkin, 2003]{Rao:2003:SCK}
 Rao, C.~V. \& Arkin, A.~P. (2003) Stochastic chemical kinetics and the
  quasi-steady state assumption: {A}pplication to the gillespie algorithm. {\em
  J. Chem. Phys.} {\bf 118} (11):4999--5010

\bibitem[Salis \& Kaznessis, 2005]{Salis:2005:AHS}
 Salis, H. \& Kaznessis, Y. (2005) Accurate hybrid stochastic simulation of a
  system of coupled chemical or biochemical reactions. {\em J. Chem. Phys.}
  {\bf }

\bibitem[Schnell \& Turner, 2004]{Schnell:2004:RKI}
 Schnell, S. \& Turner, T. (2004) Reaction kinetics in intracellular
  environments with macromolecular crowding: simulations and rate laws. {\em
  Progress in biophysics and molecular biology, } {\bf 85} (2):235--260
  
  \bibitem[Smith {\em et~al.}, 2015]{Smith:2015:ARX} Smith, S., Cianci, C., \& Grima, R. (2015) Model reduction for stochastic chemical systems with abundant species. {\em arXiv preprint arXiv:1510.03172}.

\bibitem[Srivastava {\em et~al.}, 2011]{Srivastava:2011:SQS}
 Srivastava, R., Haseltine, E.~L., Mastny, E., \& Rawlings, J.~B. (2011) The
  stochastic quasi-steady-state assumption: reducing the model but not the
  noise. {\em The Journal of Chemical Physics, } {\bf 134} (15):154109

\bibitem[Thomas {\em et~al.}, 2011]{Thomas:2011:CLS}
 Thomas, P., Straube, A.~V., \& Grima, R. (2011) Communication: Limitations of
  the stochastic quasi-steady-state approximation in open biochemical reaction
  networks. {\em The Journal of chemical physics, } {\bf 135} (18):181103

\bibitem[Wilhelm, 2009]{Wilhelm:2009:SCR}
 Wilhelm, T. (2009) The smallest chemical reaction system with bistability.
  {\em BMC systems biology, } {\bf 3} (1):90

\end{thebibliography}

\end{document}